\documentclass[a4paper, 11pt, oneside]{memoir}
\usepackage{LaTeXLoader}
\usepackage{stackengine}\stackMath
\usepackage{LaTeXSizes}
\usepackage{TimNewLists} 
\usepackage{TimNewBib}
\usepackage{TimNewTitlesSections} 
\usepackage{LaTeXFonts}
\usepackage{TimNewTheorem}
\usepackage{LTstyle}
\usepackage{TimNewStocks}
\counterwithout{equation}{chapter}
\bibliography{SetTheory.bib}
\author{Tim Button}

\begin{document}\midsloppy
\pagestyle{nicelypage}
\chapter[Part 1]{Level Theory, Part \thechapter
	\chapsubhead{Axiomatizing the bare idea of a cumulative hierarchy of sets}}\label{pt:lt}

\noindent\textcolor{blue}{This document contains preprints of Level Theory, Parts 1--3. All three papers are forthcoming at \emph{Bulletin of Symbolic Logic}.}

\begin{quote}
\textbf{Abstract.} The following bare-bones story introduces the idea of a cumulative hierarchy of pure sets: `Sets are arranged in stages. Every set is found at some stage. At any stage S: for any sets found before S, we find a set whose members are exactly those sets. We find nothing else at S'. Surprisingly, this story already guarantees that the sets are arranged in well-ordered levels, and suffices for quasi-categoricity. I show this by presenting Level Theory, a simplification of set theories due to Scott, Montague, Derrick, and Potter.
\end{quote}

\begin{epigraph}
	{What we shall try to do here is to axiomatize the types in as simple a way as possible so that everyone can agree that the idea is natural.}
	{\textcite[208]{Scott:AST}}
\end{epigraph}

\noindent 
The following bare-bones story introduces the idea of a cumulative hierarchy of pure sets:\footnote{See e.g.\ \textcites[323]{Shoenfield:AST}. I have modified Shoenfield's story in two ways. First: Shoenfield speaks of sets as `formed' at stages; I avoid this way of speaking, to avoid begging the question against platonists. Second: Shoenfield speaks of forming `collections consisting of sets' into sets; I simply speak plurally. Note that the Basic Story takes no stance on whether sets `depend' upon their members in anything other than an heuristic sense (cf.\ \cites{Incurvati:MAS}[51--69]{Incurvati:CS}).}
\begin{storytime}\textbf{The Basic Story.}	
	Sets are arranged in stages. Every set is found at some stage. At any stage $\stage{s}$: for any sets found before $\stage{s}$, we find a set whose members are exactly those sets. We find nothing else at $\stage{s}$.
\end{storytime}\noindent
This story says nothing at all about the height of any hierarchy, and apparently says almost nothing about the order-type of the stages. It lays down nothing more than the \emph{bare idea} of a pure cumulative hierarchy. Surprisingly, though, this bare idea already guarantees that the sets are arranged in well-ordered levels. Indeed, this bare idea is quasi-categorical. Otherwise put: the Basic Story pins down any cumulative hierarchy completely, modulo that hierarchy's height, on which the Story takes no stance. The aim of this paper is to show all of this.

I begin by axiomatizing the Basic Story in the most obvious way possible, obtaining Stage Theory, \ST. It is clear that any pure cumulative hierarchy satisfies \ST. Unfortunately, \ST has multiple primitives. To overcome this, I develop Level Theory, \LT. Its only primitive is $\in$, but \LT and \ST say exactly the same things about sets (see \S\S\ref{s:1:st}--\ref{s:1:LTST}). As such, any cumulative hierarchy satisfies \LT. Moreover, \LT proves that the levels are well-ordered, and \LT is quasi-categorical (see \S\S\ref{s:1:wodiscussion}--\ref{s:1:quasicat}). 

My theory \LT builds on work by Dana Scott, Richard Montague, John Derrick, and Michael Potter. I discuss their theories in \S\ref{s:1:history}, but I wish to be very clear at the outset: \LT is significantly technically simpler than its predecessors, but it owes everything to them. 

This paper is the first in a triptych. In Part \ref{pt:pst}, I explore potentialism, by considering a tensed variation of the Basic Story. In Part \ref{pt:blt}, I modify the Story again, to provide every set with a complement. Part \ref{pt:pst} presuppose Part \ref{pt:lt}, but Parts \ref{pt:lt} and \ref{pt:blt} can be read in isolation.

\setcounter{section}{-1}
\section{Preliminaries}\label{s:1:prelim}
I use second-order logic throughout. Mostly, though, my use of second-order logic is just for convenience. Except when discussing quasi-categoricity (see \S\ref{s:1:quasicat}), any second-order claim can be replaced with a first-order schema in the obvious way. In using second-order logic, I assume the Comprehension scheme, $\exists F\forall x(F(x) \liff \phi)$, for any $\phi$ not containing `$F$'. 

For readability, I concatenate infix conjunctions, writing things like $a \subseteq r \in s \in t$ for $a \subseteq r \land r \in s \land s \in t$. I also use some simple abbreviations (where $\Psi$ can be any predicate whose only free variable is $x$, and $\lhd$ can be any infix predicate):
\begin{align*}
	(\forall x : \Psi)\phi&\coloneq \forall x(\Psi(x) \lonlyif \phi) & 	(\forall x \lhd y)\phi&\coloneq\forall x(x \lhd y \lonlyif \phi)\\
	(\exists x : \Psi)\phi&\coloneq\exists x(\Psi(x) \land \phi)& 
	(\exists x \lhd y)\phi&\coloneq\exists x(x \lhd y \land \phi)
\end{align*}
When I announce a result or definition, I list in brackets the axioms I am assuming.

\section{Stage Theory}\label{s:1:st}
The Basic Story, which introduces the bare idea of a cumulative hierarchy, mentions sets and stages. To begin, then, I will present a theory which quantifies distinctly over both sorts of entities. (It is a simple modification of Boolos's \cite*{Boolos:IA} theory; see \S\S\ref{s:1:history:scott}--\ref{s:1:history:boolos}.)

Stage Theory, \ST, has two distinct sorts of variable, for \emph{sets} (lower-case italic) and for \textbf{stages} (lower-case bold). It has three primitive predicates:
\begin{listbullet}
	\item[$\in$:] a relation between sets; read `$a \in b$' as `$a$ is in $b$'
	\item[$<$:] a relation between stages; read `$\stage{r} < \stage{s}$' as `$\stage{r}$ is before $\stage{s}$'
	\item[{$\foundat$}:] a relation between a set and a stage; read `$a \foundat \stage{s}$' as `$a$ is found at $\stage{s}$'
\end{listbullet}
For brevity, I write $a \foundby \stage{s}$ for $\exists \stage{r}(a \foundat \stage{r} < \stage{s})$, i.e.\ $a$ is found before $\stage{s}$. Then \ST has five axioms:\footnote{\label{fn:lt:cheap}Classical logic yields a `cheap' proof of the existence of a stage and an empty set: by classical logic, there is some object, $a$; by \ref{st:stage} we have some \stage{s} such that $a \foundat \stage{s}$; with $F(x)$ given by $x \neq x$, \ref{st:spec} yields a set, $\emptyset$, such that $\forall x\ x\notin \emptyset$. Those who find such proofs \emph{too} cheap can adopt a free logic and then add explicit existence axioms; I will retain classical logic.}
\begin{listaxiom}
	\labitem{Extensionality}{ext} $\forall a \forall b( \forall x(x \in a \liff x \in b) \lonlyif a = b)$
	\labitem{Order}{st:ord} $\forall \stage{r} \forall \stage{s}\forall \stage{t}(\stage{r} < \stage{s} < \stage{t} \lonlyif \stage{r} < \stage{t})$
	\labitem{Staging}{st:stage} 
	$\forall a \exists \stage{s}\ \, a \foundat \stage{s}$
	\labitem{Priority}{st:pri} 
	$\forall \stage{s} (\forall a \foundat \stage{s})(\forall x \in a)x \foundby \stage{s}$
	\labitem{Specification}{st:spec} 
	$\forall F \forall \stage{s}((\forall x : F)x \foundby \stage{s} \lonlyif (\exists a \foundat \stage{s})\forall x(F(x) \liff x \in a))$
\end{listaxiom}
The first two axioms make implicit assumptions explicit: whilst I did not mention \ref{ext} in the Basic Story of a cumulative hierarchy, I take it as analytic that sets are extensional;\footnote{For brevity, I am considering hierarchies of pure sets; I revisit this in  \S\S\ref{s:1:app:ur:conventional}--\ref{s:urelements:abs}.} similarly, \ref{st:ord} records the analytic fact that `before' is a transitive relation. The remaining three axioms can then be read off the Basic Story directly: \ref{st:stage} says that every set is found at some stage; \ref{st:pri} says that a set's members are found before it; and \ref{st:spec} says that, if we find every $F$ before $\stage{s}$, then we find the set of $F$s at $\stage{s}$. So all of \ST's axioms are obviously true of the Basic Story. Otherwise put: any cumulative hierarchy obviously satisfies \ST.\footnote{\label{fn:lt:nonnull}Or, given footnote \ref{fn:lt:cheap}: any \emph{non-null} hierarchy satisfies \ST.}

This is \ST's chief virtue. Its chief drawback is that it contains multiple primitives. To see why this is a defect, suppose that we were forced to axiomatize the bare idea of a cumulative hierarchy using something like \ST's two-sorted logic. In that case, our grasp of the (cumulative iterative) notion of \emph{set} would unavoidably depend upon a concept which we had not rendered set-theoretically, namely, \emph{stage of a hierarchy}. And that would somewhat undercut the commonplace ambition, that our notion of \emph{set} might serve as a certain kind of autonomous foundation for mathematics.

\section{Level Theory}\label{s:1:lt}
To overcome this problem, I {present} Level Theory, \LT. This theory's only primitive is $\in$, but it makes exactly the same claims about sets as \ST does. I begin with a definition, due to Scott and Montague (see \S\ref{s:1:history:sm}), which forms the linchpin of this paper:\footnote{\textcites[Definition 22.4, p.161]{MontagueScottTarski}[214]{Scott:AST}. They used the `$\pot$' symbol, but not the name `potentiation'.}
\begin{define}\label{def:pot}
	For any $a$, let $a$'s  \emph{potentiation} be $\pot{a} \coloneq \Setabs{x}{\exists c(x \subseteq c \in a)}$, if it exists.\footnote{\label{fn:lt:potfudge}By the notational conventions, $\pot{a} = \Setabs{x}{(\exists c \in a)x \subseteq c} = \Setabs{x}{\exists c(x \subseteq c \land c \in a)}$. We do not initially assume that $\pot{a}$ exists for every $a$; instead, we initially treat every expression of the form `$b = \pot{a}$' as shorthand for `$\forall x(x \in b \liff \exists c(x \subseteq c \in a))$', and must double-check whether $\pot{a}$ exists. Ultimately, though, \LT proves that $\pot{a}$ exists for every $a$ (Lemma \ref{lem:lt:levof}.\ref{levofexists}).}
\end{define}\noindent
The name \emph{potentiation} emphasises the conceptual connection with powersets; note that $\pot{\{a\}} = \powerset a$.\footnote{NB: by design, \LT does not prove that every set has a powerset; for that, we have \LT + \ref{lt:cre} (see \S\ref{s:1:ltsubzf}).} The next two definitions employ this notion of potentiation (and thereby simplify definitions due to Derrick and Potter; see \S\ref{s:1:history:dp}):\footnote{\textcite[41]{Potter:STP}.}
\begin{define}
	\label{def:history}\label{def:level}
	Say that $h$ is a \emph{history}, written $\histpred(h)$, iff $(\forall x \in h) x = \pot(x \cap h)$. Say that $s$ is a \emph{level}, written $\levpred(s)$, iff $(\exists h : \histpred) s = \pot{h}$.
\end{define}\noindent
The intuitive idea behind this definition is that a history is an initial sequence of levels, and that the levels go proxy for stages. It is not obvious that this will work as described; indeed, the next two sections are dedicated to establishing this fact. But, using the notion of a level, \LT has just three axioms:\footnote{For ultra-economy, we can replace \ref{sep}+\ref{lt:strat} with $\forall F(\exists a \forall x(F(x) \liff x \in a) \liff (\exists s : \levpred)(\forall x : F)x \in s)$. We can read this as: \emph{a property determines a set iff its instances are all in some level} (cf.\ \cite[183, Definition 8.9]{ButtonWalsh:PMT}). As in footnote \ref{fn:lt:cheap}, above, the use of classical logic offers a `cheap' proof of the existence of $\emptyset$. Moreover, \LT has a model whose \emph{only} denizen is $\emptyset$.} 
\begin{listaxiom}
	\labitem{\ref{ext}}{lt:ext} $\forall a \forall b (\forall x(x \in a \liff x \in b) \lonlyif a =b)$
	\labitem{Separation}{sep} $\forall F \forall a \exists b \forall x(x \in b \liff (F(x) \land x \in a))$
	\labitem{Stratification}{lt:strat} $\forall a(\exists s : \levpred)a \subseteq s$
\end{listaxiom}

\section{The well-ordering of the levels}\label{s:1:LTwo}
In \S\ref{s:1:LTST}, I will show that \LT makes exactly the same claims about sets as \ST does. First, I must develop the elements of set theory within \LT. To do so, I need some more definitions:
\begin{define}\label{def:potent}Say that $a$ is \emph{transitive} iff $(\forall x \in a)x \subseteq a$. 
	Say that $a$ is \emph{potent} iff $\forall x(\exists c(x \subseteq c \in a) \lonlyif x \in a)$.
\end{define}\noindent
\emph{Transitivity} is completely familiar. \emph{Potency} is discussed in a few places, albeit with no standard name.\footnote{\textcites[19]{Potter:S} uses `hereditary'; \textcites[78]{Doets:RRP}[193]{ButtonWalsh:PMT} use `supertransitive'. \textcites[487]{Mathias:SMZST}[208]{Burgess:EPU} use `supertransitive' for sets which are both transitive and potent;  \textcite[1046]{LevyVaught:PPR} use `supercomplete' for such sets.} As my choice of name suggests, though, there is a tight link between the operation of potentiation (see Definition \ref{def:pot}) and the property of potency:
\begin{lem}\label{lem:null:potent} 
	If $\pot{a}$ exists, then $\pot{a}$ is potent.
\end{lem}
\begin{lem}[\ref{ext}]\label{lem:ext:idempotent}
	$a$ is potent iff $a = \pot{a}$.
\end{lem}\noindent
Recall the conventions: Lemma \ref{lem:null:potent} follows from the definitions alone, but Lemma \ref{lem:ext:idempotent} requires \ref{ext}. I leave the trivial proofs to the reader. 

My aim now is to prove several results about levels, in the sense of Definition \ref{def:level}.\footnote{The next few results simplify \textcite[41--6]{Potter:STP}. Lemma \ref{lem:es:regularity} is inspired by Potter's Proposition 3.6.4; Lemma \ref{lem:es:levhist} by Potter's Proposition 3.4.1; Lemma \ref{lem:es:acc} by Potter's Proposition 3.6.8; and Lemma \ref{lem:es:comparability} by Potter's Proposition 3.6.11.} These results do not need \LT's \ref{lt:strat} axiom, since any sets which were not subsets of levels would be irrelevant.\footnote{Cf.\ \textcite[211n.1]{Scott:AST}.} 
\begin{lem}\label{lem:es:order}Every level is potent and transitive.
\end{lem}
\begin{proof}
	Fix a level, $s$, so $s = \pot{h}$, for some history $h$. Potency follows by Lemma \ref{lem:null:potent}. For transitivity, fix $a \in s = \pot{h}$; so $a \subseteq c \in h$ for some $c$, and $c = \pot(c \cap h)$ as $h$ is a history; so $a \subseteq \pot(c \cap h) \subseteq \pot{h} = s$.
\end{proof}
\begin{lem}[\ref{ext}, \ref{sep}]\label{lem:es:regularity}
	If every $F$ is potent and something is $F$, then there is an $\in$-minimal $F$. Formally: $\forall F((\exists x F(x) \land (\forall x : F)x\text{ is potent}) \lonlyif (\exists a : F)(\forall x : F)x \notin a)$.
\end{lem}
\begin{proof}
	Let $F$ be as described and let $u$ be an $F$. Using \ref{sep} twice, let:
	\begin{align*}
		c &= \Setabs{x \in u}{(\forall y : F)x \in y} = \Setabs{x}{(\forall y : F)x \in y}\\
		d &= \Setabs{x \in c}{x \notin x}
	\end{align*}
	Clearly $d \notin c$, since otherwise $d \in d \liff d \notin d$; so there is some $a$ which is $F$ with $d \notin a$. Now if $x$ is $F$, then $d \subseteq c \subseteq x$, but $d \notin a$ and $a$ is potent, so $x \notin a$.
\end{proof}
\begin{lem}[\ref{ext}, \ref{sep}]\label{lem:es:induction}
	If some level is $F$, then there is an $\in$-minimal level which is $F$. Formally: $\forall F((\exists s : \levpred)F(s) \lonlyif (\exists s: \levpred)(F(s) \land (\forall r : \levpred)(F(r) \lonlyif r \notin s)))$
\end{lem}
\begin{proof}
	All levels are potent, by Lemma \ref{lem:es:order}; now use Lemma \ref{lem:es:regularity}. 
\end{proof}
\begin{lem}[\ref{ext}, \ref{sep}]\label{lem:es:levhist}Every member of a history is a level.
\end{lem}
\begin{proof}
	For reductio, let $h$ be a history with some non-level in it. Since $c = \pot{(c \cap h)}$ for all $c \in h$, every member of $h$ is potent by Lemma \ref{lem:null:potent}. Using Lemma \ref{lem:es:regularity}, let $a$ be an $\in$-minimal non-level in $h$. Now $a = \pot{(a \cap h)}$; and $a \cap h = \Setabs{x \in a}{x \in h}$ exists by \ref{sep}. So, to obtain our desired contradiction, it suffices to show that $a \cap h$ is a history. Fix $b \in a \cap h$. So $b$ is a level, by choice of $a$, and $b = \pot{(b \cap h)}$ as $b \in h$. If $x \in b$, then $x \subseteq b$, since $b$ is transitive by Lemma \ref{lem:es:order}; so $x \in a$, since $a$ is potent as above; hence, $b \subseteq a$. So $b = \pot{(b \cap h)} = \pot{(b \cap (a \cap h))}$. Generalising, $a \cap h$ is a history. 
\end{proof}
\begin{lem}[\ref{ext}, \ref{sep}]\label{lem:es:acc} $s = \pot{\Setabs{r \in s}{\levpred(r)}}$, for any level $s$.
\end{lem}
\begin{proof}
	Let $s$ be a level. If $a \subseteq r \in s$, then $a \in s$, as $s$ is potent by Lemma \ref{lem:es:order}. If $a \in s$, then as $s = \pot{h}$ for some history $h$, we have $a \subseteq r \in h \subseteq \pot{h} = s$ for some $r$, and $r$ is a level by Lemma \ref{lem:es:levhist}. 
\end{proof}
\begin{lem}[\ref{ext}, \ref{sep}]\label{lem:es:comparability}
	All levels are comparable.\footnote{Say that $x$ is \emph{comparable} with $y$ iff $x \in y \lor x = y \lor y \in x$} Formally:  
	$(\forall s: \levpred)(\forall t : \levpred)(s \in t \lor s = t \lor t \in s)$
\end{lem}
\begin{proof}
	For reductio, suppose that some levels are incomparable. By Lemma \ref{lem:es:induction}, there is an $\in$-minimal level, $s$, which is incomparable with some level; and by Lemma \ref{lem:es:induction} again, there is an $\in$-minimal level, $t$, which is incomparable with $s$. To complete the reductio, I will show that $s = t$. 
	
	To show that $s \subseteq t$, fix $a \in s$. So $a \subseteq r \in s$ for some level $r$, by Lemma \ref{lem:es:acc}. Now $r$ is comparable with $t$, by choice of $s$. But if either $r = t$ or $t \in r$, then $t \in s$ as $s$ is transitive, contradicting our assumption; so $r \in t$. Now $a \subseteq r \in t$, so that $a \in t$ as $t$ is potent. Generalising, $s \subseteq t$.
	
	Exactly similar reasoning, based on the choice of $t$, shows that $t \subseteq s$. So $t = s$.
\end{proof} \noindent
Rolling Lemmas \ref{lem:es:induction} and \ref{lem:es:comparability} together, we obtain the fundamental theorem of level theory:
\begin{thm}[\ref{ext}, \ref{sep}]\label{thm:es:wo}
	The levels are well-ordered by membership.
\end{thm}\noindent 
Combining this result with \ref{lt:strat}, we obtain a powerful tool, which intuitively allows us to consider the level at which a set is first found:
\begin{define}[\LT]\label{def:levof}
	Let $\levof{a}$ be the $\in$-least level with $a$ as a subset; i.e., $a \subseteq \levof{a}$ and $\lnot (\exists s : \levpred)a \subseteq s \in \levof{a}$.
\end{define}
\begin{lem}[\LT]\label{lem:lt:levof}
	For all sets $a, b$, and all levels $r, s$: 
	\setcounter{ncounts}{0}
	\begin{listn}
		\item\label{levofexists} {$\levof{a}$ and $\pot{a}$ both exist, and $\pot{a} \subseteq \levof{a}$}
		\item\label{levofnotin} $a \notin \levof{a}$
		\item\label{levofquick} $r\subseteq s$ iff $s\notin r$
		\item\label{levofidem} $s = \levof{s}$
		\item\label{levofsubs} if $b \subseteq a$, then $\levof{b} \subseteq \levof{a}$
		\item\label{levofin} if $b \in a$, then $\levof{b} \in \levof{a}$
		{\item\label{levofmin} $\levof{a} = \pot{\Setabs{\levof{x}}{x\in a}}$
			\item\label{levoffatten} if every member of $a$ is a level, then $\pot{a} = \levof{a}$}
	\end{listn}
\end{lem}
\begin{proof}
	\emphref{levofexists} {$\levof{a}$ exists by \ref{lt:strat} and Theorem \ref{thm:es:wo}. Now if $x \subseteq c \in a \subseteq \levof{a}$, then $x \in \levof{a}$ since $\levof{a}$ is potent; so $\pot{a} \subseteq \levof{a}$ exists by \ref{sep}.}
	
	\emphref{levofnotin} There is no level $t$ with $a \subseteq t \in \levof{a}$, so $a \notin \levof{a}$ by Lemma \ref{lem:es:acc}. 
	
	\emphref{levofquick} If $r \subseteq s$ then $s \notin r$ by the well-ordering of levels. Conversely, if $s \notin r$, then either $r \in s$ or $r = s$ by comparability; and $r \subseteq s$ either way, as $s$ is transitive.
	
	\emphref{levofidem}  By \eqref{levofnotin}, $s \notin \levof{s}$. By  \eqref{levofquick}, $\levof{s} \notin s$. So $s = \levof{s}$, by comparability.
	
	\emphref{levofsubs} If $b \subseteq a$ then $b \subseteq \levof{a}$. So $\levof{a} \notin \levof{b}$, by definition of $\levof{b}$, so $\levof{b} \subseteq \levof{a}$ by \eqref{levofquick}.
	
	\emphref{levofin} If $b \in a$ then $b \in \levof{a}$. By \eqref{levofnotin}, $b \notin \levof{b}$; so $\levof{a} \nsubseteq \levof{b}$, and hence $\levof{b} \in \levof{a}$ by \eqref{levofquick}.
	
	\emphref{levofmin} Let $k = \Setabs{\levof{x}}{x \in a}$. If $c \in \pot{k}$ then $c \subseteq \levof{x}$ for some $x \in a$; now $\levof{x} \in \levof{a}$ by \eqref{levofin}, so $c \in \levof{a}$. Conversely, if $c \in \levof{a}$ then $c \subseteq r \in \levof{a}$ for some level $r$ by Lemma \ref{lem:es:acc}; since $a \nsubseteq r$ by definition of $\levof{a}$, there is some $x \in a \setminus r$; now $\levof{x} \notin r$ as $r$ is potent, so that $r \subseteq \levof{x}$ by \eqref{levofquick} and hence $c \subseteq \levof{x}$; so $c \in \pot{k}$.
	
	\emphref{levoffatten} In this case, $a = \Setabs{\levof{x}}{x \in a}$ by \eqref{levofidem}, so $\levof{a} = \pot{a}$ by \eqref{levofmin}.
\end{proof}

\section{The set-theoretic equivalence of \ST and \LT}\label{s:1:LTST}
Having explained how to work within \LT, I will now make good on my earlier promise, and show that \LT and \ST make exactly the same claims about sets. More precisely, I will prove the following:
\begin{thm}\label{thm:LTST}
	$\ST \proves \phi$ iff $\LT \proves \phi$, for any {\LT-sentence $\phi$}.
\end{thm}\noindent 
To show that \ST says no more about sets than \LT does, I define a translation, $*$, from \ST-formulas into \LT-formulas. In effect, $*$ treats stages as levels, ordered by membership. Specifically, its non-trivial actions are as follows:\footnote{So the other clauses are:
			$(\lnot \phi)^* \coloneq \lnot \phi^*$; 
			$(\phi \land \psi)^* \coloneq (\phi^* \land \psi^*)$; 
			$(\forall x\phi)^* \coloneq \forall x \phi^*$; 
			$(\forall F \phi)^* \coloneq \forall F \phi^*$; and $\alpha^* \coloneq \alpha$ for all atomic formulas $\alpha$ which are not of the forms mentioned in the main text.}
	\begin{align*}
		(\stage{s} < \stage{t})^* &\coloneq \stage{s} \in \stage{t} & 
		(x \foundat \stage{s})^* &\coloneq x \subseteq \stage{s}&(\forall \stage{s} \phi)^* &\coloneq (\forall \stage{s} : \levpred)(\phi^*)
\end{align*}\noindent
After translation, we treat all first-order variables---whether bold or italic---as being of the same sort. Fairly trivially, for any \LT-sentence $\phi$, if $\ST \proves \phi$ then $\ST^* \proves \phi$. The left-to-right half of Theorem \ref{thm:LTST} now follows from this simple observation, together with the fact that $* : \ST \functionto \LT$ is an interpretation:
\begin{lem}[\LT]\label{lem:lt:sttrans} {$\ST^*$ holds.}
\end{lem}
\begin{proof}	
	\ref{ext}$^*$ is \ref{ext}. \ref{st:stage}$^*$ is \ref{lt:strat}. 
	\ref{st:ord}$^*$ holds by Lemma \ref{lem:es:order}.  
	Note that Lemma \ref{lem:es:acc} allows us to simplify $(x \foundby \stage{s})^*$, i.e.\ $(\exists \stage{r}(x \foundat \stage{r} < \stage{s}))^*$, to $(x \in \stage{s})$. Now \ref{st:pri}$^*$ holds trivially. And \ref{st:spec}$^*$ holds as if $(\forall x : F)x \in \stage{s}$, then $\Setabs{x}{F(x)} \subseteq \stage{s}$ by \ref{sep}.\footnote{Note that the $*$-translation of any \ST-Comprehension instance is an \LT-Comprehension instance.}
\end{proof}\noindent
 To obtain the right-to-left half of Theorem \ref{thm:LTST}, I must first prove some quick results in \ST:
\begin{lem}[\ST]\label{lem:st:sep} \ref{sep} holds.
\end{lem}
\begin{proof}
	By \ref{st:stage}, $a \foundat \stage{s}$ for some $\stage{s}$. By \ref{st:pri}, $(\forall x \in a)x \foundby \stage{s}$. Now use \ref{st:spec}.
\end{proof}
\begin{lem}[\ST]\label{lem:st:conversepri}$\forall \stage{s}\forall a(a \foundat \stage{s} \liff (\forall x \in a)x \foundby \stage{s})$
\end{lem} 
\begin{proof}
	Left-to-right is \ref{st:pri}. For right-to-left, suppose $(\forall x \in a)x \foundby \stage{s}$; then $\Setabs{x}{x \in a} = a \foundat \stage{s}$ by \ref{ext} and \ref{st:spec}.
\end{proof}
\noindent
I next introduce \emph{slices}. These will turn out to be levels, in the sense of Definition \ref{def:history}. Here is the definition of a slice, and some elementary results concerning slices:
\begin{define}\label{def:slice}
	For each $\stage{s}$, let $\slice{s} = \Setabs{x}{x \foundby \stage{s}}$, if it exists. Say that $a$ is a \emph{slice} iff $a = \slice{s}$ for some $\stage{s}$.
\end{define}        
\begin{lem}[\ST]\label{lem:st:slice} For any $\stage{s}$:
	\begin{listn-0}
		\item\label{slice:exists} $\slice{s}$ exists
		\item\label{slice:foundup} {$\forall \stage{r}\forall a(a \foundat \stage{r}  \leq \stage{s} \lonlyif a \foundat \stage{s})$}
		\item\label{slice:foundat} $\forall a(a \subseteq \slice{s} \liff a \foundat \stage{s})$
		\item\label{slice:trans} $\slice{s}$ is transitive
		\item\label{slice:acc} $\slice{s} = \pot{\Setabs{\slice{r}}{\slice{r} \in \slice{s}}}$
	\end{listn-0}
\end{lem}
\begin{proof}
	\emphref{slice:exists} By \ref{st:spec} and \ref{ext}.
	
	\emphref{slice:foundup} Let $a \foundat \stage{r} \leq \stage{s}$. Now $(\forall x \in a)x \foundby \stage{r}$ by \ref{st:pri}, so $(\forall x \in a)x \foundby \stage{s}$ by \ref{st:ord}, and $a \foundat \stage{s}$ by Lemma \ref{lem:st:conversepri}.
	
	\emphref{slice:foundat} $a \subseteq \slice{s}$ iff $(\forall x \in a)x \in \slice{s}$ iff $(\forall x \in a)x \foundby \stage{s}$ iff $a \foundat \stage{s}$ by Lemma \ref{lem:st:conversepri}.
	
	\emphref{slice:trans} {If} $a \in \slice{s}$, {then} $a \foundat \stage{r} < \stage{s}$ for some $\stage{r}${; hence $a \foundat \stage{s}$ and $a \subseteq \slice{s}$ by \eqref{slice:foundup}--\eqref{slice:foundat}.}
	
	\emphref{slice:acc} If $a \in \slice{s}$, then $a \foundat \stage{r} < \stage{s}$ for some $\stage{r}${; hence} $a \subseteq \slice{r} \foundat \stage{r} < \stage{s}$ by \eqref{slice:foundat}{,} so $a \subseteq \slice{r} \in \slice{s}$. If $a \subseteq \slice{r} \in \slice{s}$, then $a \subseteq \slice{r} \foundat \stage{t} < \stage{s}$ for some $\stage{t}$; now $a \subseteq \slice{r} \subseteq \slice{t}$ by \eqref{slice:foundat}, so $a \foundat \stage{t}$ by \eqref{slice:foundat}, i.e.\ $a \in \slice{s}$.	
\end{proof}\noindent
It is now easy to show that $\in$ well-orders the slices: just transcribe the proofs of Lemmas \ref{lem:es:induction} and \ref{lem:es:comparability} within \ST, replacing `levels' with `slices', noting that \ST proves \ref{sep} (see Lemma \ref{lem:st:sep}), and replacing appeal to Lemmas \ref{lem:es:order} and \ref{lem:es:acc} with Lemma \ref{lem:st:slice}.\ref{slice:trans}--\ref{slice:acc}. We can then go on to prove that the levels are the slices. 
\begin{lem}[\ST]\label{lem:st:levelsslices} $s$ is a level iff $s$ is a slice.
\end{lem}
\begin{proof}
	For induction on slices, suppose: $(\forall \slice{q} \in \slice{s})(\forall a \subseteq \slice{q})(a\text{  is a slice} \liff \levpred(a))$.  I will show that $(\forall a \subseteq \slice{s})(a\text{  is a slice} \liff \levpred(a))$. The result will follow by \ref{st:stage} and Lemma \ref{lem:st:slice}.\ref{slice:foundat}.
	
	First, fix a level $r \subseteq \slice{s}$. Let $h = \Setabs{q \in r}{\levpred(q)}$; so $r = \pot{h}$ by Lemma \ref{lem:es:acc}. (Note that \ST proves all of Lemmas \ref{lem:null:potent}--\ref{lem:es:comparability}, verbatim, since \ST proves \ref{sep}.) Fix $a \in r$; so $a \in \slice{s}$, so $a \subseteq \slice{q} \in \slice{s}$ for some $\slice{q}$ by Lemma \ref{lem:st:slice}.\ref{slice:acc}; hence, by the induction hypothesis, $a$ is a slice iff $a$ is a level. So $h = \Setabs{\slice{q}}{\slice{q} \in r}$. Noting that $h \subseteq \slice{s}$, let $\slice{t}$ be the $\in$-least slice such that $h \subseteq \slice{t}$. Since $r$ is transitive and the slices are well-ordered, $h = \Setabs{\slice{q}}{\slice{q} \in \slice{t}}$. So $r = \pot{h} = \slice{t}$ by Lemma \ref{lem:st:slice}.\ref{slice:acc}, i.e.\ $r$ is a slice.	
	
	Next, fix $\slice{r} \subseteq \slice{s}$. Let $h = \Setabs{\slice{q}}{\slice{q} \in \slice{r}}$; so $\slice{r} = \pot{h}$ by Lemma \ref{lem:st:slice}.\ref{slice:acc}; and $h = \Setabs{q \in \slice{r}}{\levpred(q)}$ by the induction hypothesis. Fix $q \in h$; since $\slice{r}$ is transitive, $q \cap h = \Setabs{p \in q}{\levpred(p)}$, so that $q = \pot{(q \cap h)}$ by Lemma \ref{lem:es:acc}. So $h$ is a history and $\slice{r} = \pot{h}$ is a level.
\end{proof}\noindent
This allows us to prove the last axiom of \LT within \ST:
\begin{lem}[\ST]\label{lem:st:strat}
	\ref{lt:strat} holds.
\end{lem}
\begin{proof}
	Fix $a$; by \ref{st:stage}, $a \foundat \stage{s}$ for some $\stage{s}$, i.e.\ $a \subseteq \slice{s}$ by Lemma \ref{lem:st:slice}.\ref{slice:foundat}, and $\slice{s}$ is a level by Lemma \ref{lem:st:levelsslices}. 
\end{proof}\noindent
So $\ST \proves \LT$, completing the proof of Theorem \ref{thm:LTST}.

\section{The inevitability of well-ordering}\label{s:1:wodiscussion}
A simple argument now establishes that \LT axiomatizes the bare idea of a cumulative hierarchy of sets:
\begin{listl-0}
	\item\label{arg:st:sound} Any cumulative hierarchy of sets satisfies \ST (see \S\ref{s:1:st}).
	\item\label{arg:lt:st} \LT is set-theoretically equivalent to {\ST} (see Theorem \ref{thm:LTST}).
	\item\label{arg:lt:sound} So: any cumulative hierarchy of sets satisfies \LT (from \eqref{arg:st:sound} and \eqref{arg:lt:st}).
\end{listl-0}
Otherwise put: \LT is true of the Basic Story I told at the start of this paper, and which I repeat here for ease of reference: \emph{Sets are arranged in stages. Every set is found at some stage. At any stage $\stage{s}$: for any sets found before $\stage{s}$, we find a set whose members are exactly those sets. We find nothing else at $\stage{s}$.}

In fact, \eqref{arg:lt:sound} takes on an even deeper significance when we reflect on just \emph{how} bare-bones this Basic Story is. The Story says that some stages are `before' others, and we can safely assume that `before' is a transitive relation on stages (hence \ST's \ref{st:ord} axiom).\footnote{In similar spirit, \textcite[323]{Shoenfield:AST} says: `We should certainly expect \emph{before} to be a partial ordering of the stages; and this is the only fact about this relation which we need for our axioms.' But Shoenfield obtains well-ordering by arguing for Foundation using a proof due to Scott (see \S\ref{s:1:history:scott}) and then using Replacement to define the $V_\alpha$s; \LT, of course, does not include Replacement (see \S\ref{s:1:ltsubzf}).} But it is not obvious, for example, that it would be inconsistent to augment the Story by saying \emph{for every stage there is an earlier stage}, or \emph{between any two stages there is another stage}. This might prompt us to start entertaining cumulative hierarchies which are ordered like the integers, or the rationals, or more exotically still. A very simple argument, however, puts an abrupt end to such speculation:
\begin{listl}
	\item\label{arg:lt:wo} \LT proves the well-ordering of the levels (see Theorem \ref{thm:es:wo}). 
	\item So: any cumulative hierarchy of sets has well-ordered levels (from \eqref{arg:lt:sound} and \eqref{arg:lt:wo}).
\end{listl}
Scott was the first to prove a well-ordering result from a similarly spartan starting point {(see \S\ref{s:1:history:sm})}, and he put the point beautifully: `This at first surprising result shows how little choice there is in setting up the type hierarchy.'\footnote{\textcite[210]{Scott:AST}.} Scott's deep observation deserves to be much more widely known. 

The connection between \ST and \LT also helps to demystify the definition of \emph{level}. Working in \ST, suppose that $h$ is an initial sequence of slices; if $\slice{s} \in h$, then $\slice{s} \cap h$ is the set of all slices less than $\slice{s}$, so that $\slice{s} = \pot(\slice{s} \cap h)$ by Lemma \ref{lem:st:slice}.\ref{slice:acc}. These observations motivate Definition \ref{def:history}. We say that $h$ is a history iff $(\forall x \in h)x = \pot{(x \cap h)}$, in the hope that, so defined, a history will be an initial sequence of slices; if it is, then the next slice in the sequence is the potentiation of that history, by Lemma \ref{lem:st:slice}.\ref{slice:acc}; and this is how we define levels.

\section{The quasi-categoricity of \LT}\label{s:1:quasicat}
We just saw that every cumulative hierarchy of sets has well-ordered levels. In fact, we can push this point further. By design, \LT says nothing about the height of any hierarchy. But, as I will show in this section, \LT is quasi-categorical. Informally, we can spell out \LT's quasi-categoricity as follows: 
\begin{listl}
	\item\label{arg:lt:quasi} Any two hierarchies satisfying \LT are structurally identical for so far as they both run, but one may be taller than the other.
\end{listl}
Since every cumulative hierarchy satisfies \LT, we obtain:
\begin{listl}
	\item\label{arg:lt:height} Any two cumulative hierarchies are structurally identical for so far as they both run, but one may be taller than the other (from \eqref{arg:lt:sound} and \eqref{arg:lt:quasi}).
\end{listl} 
So, echoing Scott: when we set up a cumulative hierarchy, our only choice is how tall to make it.

It just remains to establish \eqref{arg:lt:quasi}, i.e.\ to show that \LT is quasi-categorical. In fact, there are {at least} two ways to explicate the informal idea of quasi-categoricity, and \LT is quasi-categorical on both explications. (Note that both ways make essential use of second-order logic; this is the only section of the paper where my use of second-order logic is not merely for convenience.)

The first notion of quasi-categoricity is familiar from Zermelo. Working in some (set-theoretic) model theory, we define the $V_\alpha$s as usual:
\begin{align*}
	V_0 & = \emptyset; & V_{\alpha + 1} & = \powerset V_\alpha; & V_{\alpha} &= \bigcup_{\beta \in \alpha} V_\beta \text{ when $\alpha$ is a limit}
\end{align*}
Each $V_\alpha$ then naturally yields a set-theoretic structure, $\model{V}_\alpha$, whose domain is $V_\alpha$, and which interprets  `$\in$' as membership-restricted-to-$V_\alpha$, i.e.\ $\Setabs{\langle x, y\rangle \in V_\alpha \times V_\alpha}{x \in y}$. We then have the following result, using full second-order logic: $\model{M} \models \ZF$ iff $\model{M} \isomorphic \model{V}_\alpha$ for some strongly inaccessible $\alpha$.\footnote{\textcite{Zermelo:GM}. For an accessible proof, see \textcite[\S8.A]{ButtonWalsh:PMT}.} There is an analogous quasi-categoricity result for \LT:\footnote{\textcite[\S8.C]{ButtonWalsh:PMT} prove this for Potter's theory (see \S\ref{s:1:history:dp}); the same proof works for \LT. The same remark applies to the other results mentioned in this section.  We could obtain external categoricity using only first-order logic, if we augmented \LT with some axiom of the form `there are exactly $n$ levels'.}
\begin{thm}[{in full second-order logic}]\label{thm:LTexternalcat}$\model{M} \models \LT$ iff $\model{M} \isomorphic \model{V}_\alpha$ for some $\alpha > 0$.
\end{thm}\noindent
This shows that any two hierarchies satisfying \LT (read that phrase as `any models of \LT') are structurally identical (read that phrase as `are isomorphic') for so far as they both run (read that phrase in the light of the well-ordering of the $V_\alpha$s, established in the model theory). In short, \LT is quasi-categorical, on a model-theoretic (`external') way of understanding quasi-categoricity.

There is also, though, an object-language (`internal') way to understand quasi-categoricity.\footnote{This has been brought out by  \textcites{Parsons:UNN}{Parsons:MTO}{McGee:HLML}{VaananenWang:ICAST}. The remainder of this section presents specific elements of \textcite[ch.11]{ButtonWalsh:PMT}.} Since this idea is less familiar, I will spend some time unpacking it.

In embracing \ref{ext}, \LT assumes that everything is a pure set. There is a quick-and-dirty way to avoid this assumption. First, introduce a new predicate, $\purepred$; intuitively, this should apply to the pure sets. Next, relativise \LT to $\purepred$, via the following formula:\footnote{Here, `$\subseteq$' and `$\levpred$' should be defined in terms of `$\varin$' rather than `$\in$'; similarly for `$\levof$' in Theorem \ref{thm:LTInternal}. For now, we can treat `$\purepred$' as a primitive; but see Definition \ref{def:purepred}.}
\begin{align*}
	\LT(\purepred, \varin) \coloneq{} 
	& (\forall a : \purepred)(\forall b : \purepred)(\forall x (x \varin a \liff x \varin b)\lonlyif a = b) \land {} \\
	& \forall F(\forall a: \purepred)(\exists b : \purepred)\forall x(x \varin b \liff (F(x) \land x \varin a)) \land {}\\
	& (\forall a : \purepred)(\exists s : \levpred)a \subseteq s \land{}\\
	& \forall x \forall y(x \varin y \lonlyif (\purepred(x) \land \purepred(y)))
\end{align*}\noindent
The first three conjuncts tell us that the pure sets satisfy \LT;\footnote{With one insignificant caveat (see footnotes \ref{fn:lt:cheap} and \ref{fn:lt:nonnull}): whereas classical logic guarantees that any model of \LT contains an empty set,  $\LT(\purepred, \varin)$ allows that there may be no pure sets.} the last says that, when we use `$\varin$', we restrict our attention to membership facts between pure sets. Using this formula, I can now state the internal quasi-categoricity result (I have labelled the lines to facilitate its explanation):
\begin{thm}\label{thm:LTInternal} This is a deductive theorem of impredicative second-order logic:
	\begin{align}
		(\LT(\purepred_1, &\varin_1) \land \LT(\purepred_2, \varin_2))\lonlyif \nonumber\\
		\exists R(
		& \label{lti:VtoV} 
		\forall v \forall y(R(v,y) \lonlyif (\purepred_1(v) \land \purepred_2(y))) \land {}\\
		& \label{lti:exhausts} 
		((\forall v : \purepred_1)\exists y R(v,y) \lor (\forall y : \purepred_2)\exists v R(v,y)) \land {}\\
		& \label{lti:preserves} 
		\forall v\forall y\forall x\forall z ((R(v, y) \land R(x, z)) \lonlyif (v \varin_1 x \liff y \varin_2 z))\land {}\\
		& \label{lti:functional} 
		\forall  v \forall y\forall z ((R(v,y) \land R(v,z)) \lonlyif y = z)\land {}\\
		& \label{lti:injective} 
		\forall  y \forall v\forall x ((R(v,y) \land R(x, y)) \lonlyif v= x)\land {}\\
		& \label{lti:segmentsa} 
		\forall v (\exists y R(v, y) \lonlyif (\forall x  \subseteq_1 \levof{}_1{v}) \exists z R(x, z))\land {}\\
		& \label{lti:segmentsb} 
		\forall y (\exists v R(v, y) \lonlyif (\forall z  \subseteq_2 \levof{}_2{y}) \exists x R(x, z)))
	\end{align}
\end{thm}\noindent
Intuitively, the point is this. Suppose two people are using their versions of \LT, subscripted with `$1$' and `$2$' respectively. Then there is some second-order entity, a relation $R$, which takes us between their sets \eqref{lti:VtoV}, exhausting the sets of one or the other person \eqref{lti:exhausts}; which preserves membership \eqref{lti:preserves};  which is functional \eqref{lti:functional} and injective \eqref{lti:injective}; and whose domain is an initial segment of one \eqref{lti:segmentsa} or the other's  \eqref{lti:segmentsb} hierarchy. Otherwise put: \LT is (internally) quasi-categorical.

As a bonus, this internal \emph{quasi}-categoricity result can be lifted into an internal \emph{total}-categoricity result. To explain how, consider this abbreviation (where `$P$' is a second-order function-variable):
$$\somany x \Phi(x) \coloneq \exists P(\forall x \Phi(P(x)) \land (\forall y : \Phi)\exists ! x\ P(x) = y)$$
This formalizes the idea that there are as many $\Phi$s as there are objects \emph{simpliciter}, i.e., that there is a bijection between the $\Phi$s and the universe. We can use this notation to state our internal total-categoricity result:
\begin{thm}\label{thm:LTInternalFull} This is a deductive theorem of impredicative second-order logic:
	\begin{align*}
		(\LT(\purepred_1,& \varin_1) \land \somany x\,  \purepred_1(x) \land \LT(\purepred_2, \varin_2) \land \somany x\, \purepred_2(x))\lonlyif \\
		\exists R(&\forall v \forall y(R(v,y) \lonlyif {}(\purepred_1(v) \land \purepred_2(y))) \land {} \\
		&(\forall  v : \purepred_1)\exists !y R(v,y) \land (\forall y : \purepred_2) \exists! v R(v,y)\land {} \\
		&\forall v\forall y\forall x\forall z ((R(v, y) \land R(x, z)) \lonlyif (v \varin_1 x \liff y \varin_2 z)))
	\end{align*}
\end{thm}\noindent
Intuitively, if both \LT-like hierarchies are as large as the universe, then there is a structure-preserving \emph{bijection} between them. To see the significance of this result, note that it is common to claim that there are \emph{absolutely infinitely many} pure sets. Whatever exactly this is meant to mean, it must surely entail that $\somany x\, \purepred(x)$. So Theorem \ref{thm:LTInternalFull} tells us that absolutely infinite \LT-like hierarchies are (internally) isomorphic.

\section{\LT as a subtheory of \ZF}\label{s:1:ltsubzf}
I have shown that any cumulative hierarchy satisfies \LT, so that, in setting up a cumulative hierarchy, our only freedom of choice concerns its height. To make all of this more familiar, though, it is worth commenting on \LT's relationship to \ZF, the `industry standard' set theory.

Unsurprisingly, \ZF proves \LT. In more detail: working in \ZF, define the $V_\alpha$s as usual; we can then show that the $V_\alpha$s are the levels;\footnote{\emph{Proof sketch.} Working in \ZF, fix $\alpha$, and suppose for induction that $(\forall \beta < \alpha)(\forall x \subseteq V_\beta)(\levpred(x) \liff \exists \delta\, \ x = V_\delta)$. {Fix $V_\gamma \subseteq V_\alpha$; then $V_\gamma = \pot{\Setabs{V_\delta}{V_\delta\in V_\gamma}} = \pot{\Setabs{s \in V_\gamma}{\levpred(s)}}$ by the induction hypothesis, which is a level by Lemma \ref{lem:es:acc}. Similarly,} if $s \subseteq V_\alpha$ is a level, then $\exists \delta\, \ s = V_\delta$.} so \ref{lt:strat} holds as every set is a subset of some $V_\alpha$. 

Of course, \ZF is much stronger than \LT, since \LT deliberately says nothing about the height of the cumulative hierarchy. If we want to set up a tall hierarchy, then three axioms naturally suggest themselves (where `$P$' is a second-order function-variable in the statement of \ref{lt:rep}):\footnote{For \ref{lt:cre}, cf.\ \textcites[142]{Montague:STHOL}[212]{Scott:AST}[20--1]{Potter:S}[61--2]{Potter:STP}. For \ref{lt:inf}, see \textcite[68--70]{Potter:STP} and \citepossess[8]{Boolos:IA} axiom \ref{boolos:inf}, which I discuss in \S\ref{s:1:history:boolos}.}
\begin{listaxiom}
	\labitem{Endless}{lt:cre} $(\forall s : \levpred)(\exists t : \levpred)s \in t$
	\labitem{Infinity}{lt:inf} $(\exists s : \levpred)((\exists q : \levpred)q \in s \land (\forall q : \levpred)(q \in s \lonlyif (\exists r : \levpred)q \in r \in s))$
	\labitem{Unbounded}{lt:rep} $\forall P\forall a (\exists s : \levpred)(\forall x \in a) P(x) \in s$
\end{listaxiom}\noindent
\ref{lt:cre} says there is no last level. \ref{lt:inf} says that there is an infinite level, i.e.\ a level with no immediate predecessor. \ref{lt:rep} states that the hierarchy of levels is so tall that no set can be mapped unboundedly into it. We now have some nice facts, whose proofs I leave to the reader:\footnote{Cf.\ \textcites[212]{Scott:AST}[20--4]{Potter:S}[47--9, 61--2]{Potter:STP}.}
\begin{prop}\label{prop:lt:zfconnections}\textcolor{white}{.}
	\begin{listn-0}
		\item\label{gets:suf} \LT proves Separation, Union, and Foundation.
		\item $\LT + \text{\ref{lt:cre}}$ proves Pairing and Powersets.
		\item\label{gets:inf} $\LT + \text{\ref{lt:cre}} + \textnormal{\ref{lt:inf}}$ proves Zermelo's axiom of infinity.\footnote{i.e.\ $(\exists w \ni \emptyset)(\forall x \in w)x \cup \{x\} \in w$.}
		\item\label{gets:zffin} $\LT + \text{\ref{lt:cre}} +  \lnot\text{\ref{lt:inf}}$ is equivalent to \ZFfin.\footnote{The latter is the theory with all of \ZF's axioms except that: (i) Zermelo's axiom of infinity is replaced with its negation; and (ii) it has a new axiom, $\forall a(\exists t \supseteq a)(t\text{ is transitive})$.}
		\item $\LT + \text{\ref{lt:inf}} + \textnormal{\ref{lt:rep}}$ proves \textnormal{\ref{lt:cre}}.
		\item\label{gets:zf} $\LT + \text{\ref{lt:inf}} + \text{\ref{lt:rep}}$ is equivalent to \ZF.
	\end{listn-0}
\end{prop}\noindent
Facts \eqref{gets:suf}--\eqref{gets:inf} show that $\LT + \textnormal{\ref{lt:cre}} + \textnormal{\ref{lt:inf}}$ extends Zermelo's \Zermelo. This extension is {strict}, since \ref{lt:strat} is independent from \Zermelo.\footnote{\label{fn:mathias}\textcite[293ff]{Potter:STP} makes a similar point. The independence is immediate from the fact that there are models of (even second-order) \Zermelo which fail to satisfy $\forall a (\exists c \supseteq a)(c\text{ is transitive})$; see \textcite[111]{Drake:ST}. For detailed discussions of \Zermelo's weaknesses, as either a first- or second-order theory, see \textcites{Mathias:SMZST}{Uzquiano:MSOZST}. (As mentioned in the introduction, although I have formulated \LT as a second-order theory, it has a natural first-orderization. Read uniformly as either first-order or second-order theories, and closing under provability, the point is: $\Zermelo \subsetneq \LT + \textnormal{\ref{lt:cre}} + \textnormal{\ref{lt:inf}} \subsetneq \ZF$.)} Fact \eqref{gets:zf} then offers a neat way to conceive of \ZF, as extending the theory which holds of any cumulative hierarchy, i.e.\ \LT, with specific claims about the hierarchy's height. 

\section{Conclusion, and \LT's predecessors}\label{s:1:history}
The theory \LT holds of every cumulative hierarchy. Since \LT is also quasi-categorical, the only choice we have, in setting up a cumulative hierarchy, is over the hierarchy's height. 

I will close this paper by discussing \LT's predecessors, in roughly chronological order.

\subsection{Scott}\label{s:1:history:scott}
At a talk in 1957, Scott presented what seems to have been the first theory of stages. This was an axiomatic theory of \emph{ranks}, in the sense of the $V_\alpha$s. Writing `$a < b$' for `$a$ has lesser rank than $b$', Scott's suggested axioms were \ref{ext} and:\footnote{\textcite{Scott:NRST}; I have tweaked the presentation slightly.} 
\begin{listaxiom}
	\labitem{}{scott1} $\forall a \forall b(a < b \liff (\exists x < b)x\nless a)$
	\labitem{}{scott2} $\forall F(\forall a((\forall x < a)F(x) \lonlyif F(a)) \lonlyif \forall a F(a))$
	\labitem{}{scott3} $\forall F\forall a \exists b \forall x(x \in b \liff (F(x) \land x < a))$
\end{listaxiom}
This 1957 theory is clearly satisfied in any $\model{V}_\alpha$ with $\alpha > 0$, when $\in$ and $<$ are given the obvious interpretations. However, it has some unintended models. 
\begin{example}
	Let the domain have two sets: $\emptyset$ and a Quine atom $a = \{a\}$. Let $a < \emptyset$. This is a model of the 1957 theory, since $<$ is trivially a well-order, and since the only sets given by the third axiom are $\emptyset$ and $\{a\} = a$. 
\end{example}
\begin{example}
	Let the domain have four sets: $\emptyset$, $\{\emptyset\}$, $\{\{\emptyset\}\}$ and $\{\emptyset, \{\emptyset\}\}$. Permute the usual rank relation, so that $\{\emptyset\} < \emptyset < \{\{\emptyset\}\}, \{\emptyset, \{\emptyset\}\}$, with $\{\{\emptyset\}\}$ and $\{\emptyset, \{\emptyset\}\}$ incomparable.
\end{example}\noindent
At a talk in 1967, Scott provided a vastly improved theory of stages. I will present the 1967 theory in a slightly simplified form, starting with a definition given later by Potter (see \S\ref{s:1:history:dp}):
\begin{define}\label{def:acc}
	For each set $a$, let $\accumulation{a} = \Setabs{x}{(\exists c \in a)(x \in c \lor x \subseteq c)}$, if it exists.
\end{define}\noindent
Scott's 1967 theory treats the notion of \emph{level} as a primitive, which applies to certain sets. Temporarily using bold-face letters to range over these \emph{levels}, the 1967 theory comprises just \ref{ext}, \ref{sep}, and two further axioms:\footnote{\textcite[208--9]{Scott:AST}. Scott allowed urelements, which I am ignoring for ease of presentation (though see \S\ref{s:1:app:ur:conventional}).}
\begin{listaxiom}
	\labitem{{Accumulation}}{scott:acc} 
	$\forall \stage{s}\ \, \stage{s} = \accumulation{\Setabs{\stage{r}}{\stage{r} \in \stage{s}}}$
	\labitem{{Restriction}}{scott:res} 
	$\forall a \exists \stage{s}\ \, a \subseteq \stage{s}$
\end{listaxiom}
Scott's 1967 theory (unlike his 1957 theory) does not explicitly state that the levels are well-ordered; instead, the 1967 theory \emph{proves} the well-ordering of the levels (cf.\ \S\ref{s:1:wodiscussion}).\footnote{\citepossess[211--2]{Scott:AST} proof uses the idea of a \emph{grounded} set, introduced by \textcite{Montague:PGC}.}  We have Scott to thank for a truly remarkable bit of mathematics-cum-conceptual-analysis.

Scott's 1967 theory obviously inspires \ST: compare his \ref{scott:res} axiom with my \ref{st:stage} (and \ref{lt:strat}), and his \ref{scott:acc} axiom with my Lemma \ref{lem:st:slice}.\ref{slice:acc} (and Lemma \ref{lem:es:acc}). Moreover, Scott's 1967 theory and \ST make exactly the same claims about sets (cf.\ Theorem \ref{thm:LTST}). But I used \ST in \S\ref{s:1:st}, rather than Scott's 1967 theory, since \ST is easier to motivate. In particular, Scott simply instructs us to write `$\stage{s} \in \stage{t}$' for `$\stage{s}$ is before $\stage{t}$', and his justification of \ref{scott:acc} amounts to stipulating that `a given level is \emph{nothing more than} the accumulation of all the members and subsets of all the \emph{earlier} levels'.\footnote{Both quotes from \textcite[209]{Scott:AST}; his emphasis; variables adjusted to match surrounding text.} Both claims are very natural, and they are true in \LT; but it is not immediately obvious that they are true of the Basic Story I told in the introduction. (In fairness to Scott, he does not start with that story, but with a related justificatory tale.)

\subsection{Boolos and Shoenfield}\label{s:1:history:boolos}
The second source of inspiration for \ST is Boolos. He first presented a theory of stages in \cite*{Boolos:ICS}, which included explicit axioms stating that the stages are well-ordered;\footnote{See \citepossess[223--4]{Boolos:ICS} I--IV and Induction Axioms.} this theory has several similarities with \cite{Shoenfield:ML}.\footnote{\textcite[238--40]{Shoenfield:ML}.} Boolos then presented a better theory of stages in \cite*{Boolos:IA}, explicitly drawing from Scott's 1967 theory to prove (rather than assume) a principle of induction for stages.\footnote{\textcites[211--2]{Scott:AST}[11--12]{Boolos:IA}; Boolos cites \citepossess[327]{Shoenfield:AST} presentation of Scott.} My theory \ST tweaks Boolos's \cite*{Boolos:IA} theory in three ways. 

\emph{First.} Boolos has qualms about how to justify \ref{ext};\footnote{\textcite[10--11]{Boolos:IA}.} I have no such qualms. 

\emph{Second.} Boolos aims to vindicate the traditional Zermelian axioms of Foundation, Union, Pairing, Separation, Powersets, and Infinity. To secure these last two axioms, his \cite*{Boolos:IA} theory contains:
\begin{listaxiom}
	\labitem{Inf}{boolos:inf} $\exists\stage{t}(\exists \stage{r} \ \stage{r} < \stage{t} \land (\forall \stage{r} < \stage{t})\exists \stage{s}(\stage{r} < \stage{s} < \stage{t}))$
	\labitem{Net}{boolos:net} $\forall \stage{r} \forall \stage{s} \exists\stage{t}(\stage{r} < \stage{t} \land \stage{s} < \stage{t})$
\end{listaxiom}
Boolos's \ref{boolos:inf} guarantees there is a stage with infinitely many predecessors, and his \ref{boolos:net} guarantees that there is no last stage. Since \ST is deliberately silent on the height of any cumulative hierarchy, it has no similar axioms. However, if I had wanted to augment \ST with the claim that there is no last stage, I would have offered $\forall \stage{s} \exists \stage{t}\, \ \stage{s} < \stage{t}$ (cf.\ \ref{lt:cre}, from \S\ref{s:1:ltsubzf}). Boolos's \ref{boolos:net} says more than this; it guarantees that stages are directed. Boolos's proof of Pairing relies upon this directedness,\footnote{\citepossess[19]{Boolos:IA} proof is as follows. Fix $a$ and $b$; by \ref{st:stage}, there are $\stage{r}$ and $\stage{s}$ with $a \foundat \stage{r}$ and $b \foundat \stage{s}$. By \ref{boolos:net}, there is some $\stage{t}$ after both $\stage{r}$ and $\stage{s}$. So by \ref{boolos:spec} there is a set whose members are exactly $a$ and $b$.} but I cannot see why Boolos felt independently entitled to adopt \ref{boolos:net} rather than the weaker principle. 

\emph{Third.} The remainder of Boolos's \cite*{Boolos:IA} theory comprises \ref{st:ord}, \ref{st:stage}, and these two axioms:\footnote{\textcite[8]{Boolos:IA} formulates \ref{boolos:spec} as a first-order scheme, but considers the second-order axiom on the next page.}
\begin{listaxiom}
	\labitem{{When}}{boolos:when} 
	$\forall \stage{s} \forall a(a \foundat \stage{s} \liff (\forall x \in a)x \foundby \stage{s})$
	\labitem{Spec}{boolos:spec} $\forall F \forall \stage{s} ((\forall x : F)x \foundby \stage{s} \lonlyif \exists a\forall x(F(x) \liff x \in a))$
\end{listaxiom}
In the presence of \ref{ext}, the axioms \ref{boolos:when}+\ref{boolos:spec} are equivalent to \ST's \ref{st:pri}+\ref{st:spec}; but we need \ref{ext} to prove the right-to-left direction of \ref{boolos:when} from \ref{st:pri}+\ref{st:spec} (see Lemma \ref{lem:st:conversepri}). Moreover, given Boolos's qualms about \ref{ext}, he cannot provide an intuitive justification for the right-to-left direction of \ref{boolos:when}. If $(\forall x \in a)x \foundby \stage{s}$, then there should certainly be some $b \foundat \stage{s}$ such that $\forall x(x \in b \liff x \in a)$; but only \ref{ext} can justify the assertion that $b = a$. Crucially for Boolos's aims, though, Powersets can fail if we replace \ref{boolos:when}+\ref{boolos:spec} with \ref{st:pri}+\ref{st:spec} in Boolos's theory: without \ref{ext} or the right-to-left direction of \ref{boolos:when}, we might keep finding new empty sets at every stage in the hierarchy; there will then be no stage by which every subset of a set has been found, and hence no stage at which any powerset can be found.

\subsection{Scott and Montague}\label{s:1:history:sm}
I now want to return to Scott's 1967 theory. As mentioned in \S\ref{s:1:history:scott}, this theory initially takes the notion of \emph{level} as primitive. However, Scott notes that the primitive can be eliminated, by proving within the 1967 theory that $s$ is a level iff $\pot{s} \subseteq s \land (\forall a \in s)(\exists h \in s)(\forall k \subseteq h)(\pot{k} \in s \land (\pot{k} \in h \lor a \subseteq \pot{k}))$. Scott developed this ideologically-spartan theory in joint work with Montague; they described their theory as `rank free', so I will call it \rankfreetheory.\footnote{\textcites[139]{Montague:STHOL}[161--2]{MontagueScottTarski}[214]{Scott:AST}.} It has just three axioms: \ref{ext}, \ref{sep}, and 
\begin{listaxiom}
	\labitem{Hierarchy}{hier} $\forall a \exists h (\forall k\subseteq h){(\exists s = \pot{k})(s \in h \lor a \subseteq s)}$
\end{listaxiom}
The point of calling it `rank free' was to highlight that \rankfreetheory takes no stance on the number of ranks in the hierarchy. More precisely, we have the external quasi-categoricity result that $\model{M} \models \rankfreetheory$ iff $\model{M} \isomorphic \model{V}_\alpha$ for some $\alpha > 0$ (assuming full second-order logic; cf.\ Theorem \ref{thm:LTexternalcat}). To establish this, Montague and Scott first say that $h$ is a \emph{hierarchy} iff $(\forall k \subseteq h)(h \subseteq \pot{k} \lor \pot{k} = \bigcap(h \setminus \pot{k}))$. They then let $\text{R}a \coloneq \bigcap \{\pot{h} : h\text{ is a hierarchy} \land a \subseteq \pot{h}\}$ for each $a$, and show that $\text{R}a$ serves the role of $a$'s `rank' (cf.\ \LT's notion of $\levof{a}$, as laid down in Definition \ref{def:levof}).

Unfortunately, as Scott himself put it, the deductions from these axioms and definitions `are quite lengthy'.\footnote{\textcite[214]{Scott:AST}. Indeed, it occupies 13 dense sides of \textcite[161--74]{MontagueScottTarski}. The two key definitions are 22.7 and 22.21.} This led Scott to dismiss the significance of \rankfreetheory, writing: `there seems to be no technical or conceptual advantage in reducing the number of primitive notions to the minimum.'\footnote{\textcite[214]{Scott:AST}.}

Still, these lengthy deductions were intended to form a section of a monograph on axiomatic set theory. A complete manuscript of this monograph exists, \textcite{MontagueScottTarski}, containing very minor markups, handwritten notes to the printers, and an accompanying list of `Things to be Done' which amounts to nothing more than writing an Introduction and dealing with the mundane logistics of publication. Everything, in short, was almost ready to print.

Sadly, it was never printed. This was a serious loss. As explained in \S\ref{s:1:st}, there are good {philosophical} reasons for `reducing the number of primitive notions to the minimum.' Moreover, whilst Montague's and Scott's \emph{deductions} were `quite lengthy', the \emph{axioms} of \rankfreetheory are quite elegant. The lengthiness of the deductions from \rankfreetheory is down to the  awkwardness of the definitions of \emph{hierarchy} and $\text{R}a$. If Montague and Scott had been aware of the definition of \emph{history} and \emph{level}, as given in Definition \ref{def:history}, they could have given some much briefer deductions. Indeed, these definitions make it easy to prove that  \rankfreetheory and \LT are equivalent. One direction of this equivalence is easy:
\begin{prop}[\LT]\label{prop:SPtoMS} \rankfreetheory holds.
\end{prop}
\begin{proof}
	It suffices to prove \ref{hier}. Fix $a$, let $h = \Setabs{s \in \levof{a}}{\levpred(s)}$ and fix $k \subseteq h$. Now $\pot{k} = \levof{k}$ by Lemma \ref{lem:lt:levof}.\ref{levoffatten}; so if $\pot{k} = \levof{k} \notin h$, then $\levof{k} \notin \levof{a}$, so $a \subseteq \levof{a} \subseteq \levof{k} = \pot{k}$ by Lemma \ref{lem:lt:levof}.\ref{levofquick}.
\end{proof}\noindent
For the other direction of the quivalence, I must first prove some quick facts in \rankfreetheory:
\begin{lem}[\rankfreetheory]\label{lem:pot} For all $a$:
	\setcounter{ncounts}{0}
	\begin{listn}
		\item \label{ms:naughty} if $\pot{a}$ exists, then $\pot{a} \notin a$
		\item\label{ms:exists} $\pot{a}$ exists 		\item\label{ms:nextlevel} if every member of $a$ is a level, then $\pot{a}$ is a level
	\end{listn}
\end{lem}	
\begin{proof}
	\emphref{ms:naughty} {If $\pot{a} \in a$, then $(\forall c \subseteq \pot{a})c \in \pot{a}$. But this is impossible: by} \ref{sep}, let $d = \Setabs{x \in \pot{a}}{x \notin x}$; then $d \notin \pot{a}$
	
	\emphref{ms:exists} Fix $a$, and let $h$ witness \ref{hier}. Let $k = h$, so that $\pot{h}$ exists and $\pot{h} \in h \lor a \subseteq \pot{h}$, i.e.\ $a \subseteq \pot{h}$ by \eqref{ms:naughty}. Since $\pot{h}$ is potent by Lemma \ref{lem:null:potent}, $\pot{a} \subseteq \pot{h}$ exists by \ref{sep} on $\pot{h}$. 
	
	\emphref{ms:nextlevel} Using \ref{sep} and \eqref{ms:exists}, let $h = \Setabs{s \in \pot{a}}{\levpred(s)}$. I will first prove that $\pot{h} = \pot{a}$, and then that $h$ is a history, so that $\pot{h} = \pot{a}$ is a level. 
	
	To see that $\pot{a} = \pot{h}$: since $h \subseteq \pot{a}$, we have $\pot{h} \subseteq \pot{\pot{a}} =  \pot{a}$ by Lemmas \ref{lem:null:potent}--\ref{lem:ext:idempotent}; and if $x \in \pot{a}$ then $x \subseteq r \in a$ for some level $r$, so $r \in h$, and hence $x \in \pot{h}$.
	
	To see that $h$ is a history, fix $s \in h$; it suffices to show that $s = \pot{(s \cap h)}$. Since $s$ is a level,  $\pot{(s \cap h)} \subseteq \pot{s} = s$ by Lemmas \ref{lem:ext:idempotent}--\ref{lem:es:order}. To see that $s \subseteq \pot{(s \cap h)}$, fix $x \in s$; now $x \subseteq r \in s$ for some level $r$ by Lemma \ref{lem:es:acc}; and {$r \subseteq s \in \pot{a}$ by Lemma \ref{lem:es:order}, so $r \in \pot{a}$ by Lemma \ref{lem:null:potent} and hence $r \in h$}; so $x \subseteq r \in (s \cap h)$, i.e.\ $x \in \pot{(s \cap h)}$. 
\end{proof}	
\begin{prop}[\rankfreetheory] \LT holds.
\end{prop}
\begin{proof}
	It suffices to prove \ref{lt:strat}. Fix $a$, and let $h$ witness \ref{hier}, i.e., $(\forall k \subseteq h)(\pot{k} \in h \lor a \subseteq \pot{k})$. Let $k = \Setabs{s \in h}{\levpred(s)}$. By Lemma \ref{lem:pot}.\ref{ms:nextlevel}, $\pot{k}$ is a level. Now if $\pot{k} \in h$, then $\pot{k} \in k$, contradicting Lemma \ref{lem:pot}.\ref{ms:naughty}; so $a \subseteq \pot{k}$. 
\end{proof}\noindent
This last proof helps to explain the intuitive idea behind \rankfreetheory's axiom \ref{hier}.\footnote{Cf.\ \textcite[162]{MontagueScottTarski}.} Roughly, the $h$ guaranteed to exist by \ref{hier} has this property: for any initial sequence of levels $k \subseteq h$, the next level after all of them is $\pot{k}$; and if $a$ is not a subset of $\pot{k}$, then $\pot{k}$ is in $h$; and \emph{hence} (but here I invoke a transfinite induction) the members of $h$ are all the levels up to the first level including $a$. In short, the fundamental idea behind \rankfreetheory is quite elegant.

\subsection{Derrick and Potter}\label{s:1:history:dp}
As mentioned in \S\ref{s:1:lt}, my definition of \emph{level} is inspired by Derrick and Potter,\footnote{See especially \textcites[16--20]{Potter:S}[41--7]{Potter:STP}.} but I have simplified it. Here is a little more detail about that simplification. In his \cite*{Potter:S} book, Potter explicitly built on Scott's 1967 theory and also on Derrick's unpublished lecture notes.\footnote{\textcites[183--4]{Potter:IST}[22]{Potter:S}[vii, 54]{Potter:STP}.} Now, Scott's \ref{scott:acc} axiom (see \S\ref{s:1:history:scott}) formalizes the claim that `a given level is \emph{nothing more than} the accumulation of all the members and subsets of all the \emph{earlier} levels'.\footnote{\textcite[209]{Scott:AST}.} This suggests the use of the $\accumulation$-operator, and so Potter offers Definition \ref{def:acc}.\footnote{\textcites[16]{Potter:S}[41, 50]{Potter:STP}.} Potter then supplies the definition of \emph{history} and \emph{level} given in Definition \ref{def:history}, but using $\accumulation$ rather than $\pot$. So, Potter stipulates that $h$ is a history iff $(\forall x \in h)x = \accumulation{(x \cap h)}$, and that $s$ is a level iff $s = \accumulation{h}$ for some history $h$. Potter then proves that, so defined, the levels are well-ordered. And his own theory of levels is, in effect, just \LT, with this slightly different explicit definition of `$\levpred$'.\footnote{There are three other small differences: (1) Potter allows urelements; (2) he provides a first-order theory; (3) he offers a slightly more restricted version of \ref{sep}, whose second-order formulation is $\forall F (\forall s : \levpred)\exists b\forall x(x \in b \liff (F(x) \land x \in s))$, but this trivially entails the unrestricted version of \ref{sep} given (Potter's version of) \ref{lt:strat}.} But the use of $\pot$, rather than $\accumulation$, simplifies things significantly, as illustrated by the brevity of \S\ref{s:1:LTwo}.

\startappendix
\section{Adding urelements}\label{s:1:app:ur:conventional}
In this paper, I restricted my attention to pure sets.\footnote{\textcites[139]{Montague:STHOL}[214]{Scott:AST}{Potter:S}{Potter:STP} accommodate urelements from the outset.} This was only for ease of exposition; in this appendix and the next, I will remove this simplifying assumption. 

To accommodate urelements, we must tweak the Basic Story. The easiest way to do this (which I revisit in \S\ref{s:urelements:abs}) is to assume that the urelements are `always' available to be collected into sets:
\begin{storytime}\textbf{The Urelemental Story.}
	Sets are arranged in stages. Every set is found at some stage. At any stage $\stage{s}$: for any things, each of which is either a set found before $\stage{s}$ or an urelement, we find a set whose members are exactly those things. We find nothing else at $\stage{s}$.
\end{storytime}\noindent 
To formalize this Story, we need a new primitive predicate, enabling us to distinguish sets from urelements: we take $\setpred$ as primitive, and define $\urpred(x) \coloneq \lnot \setpred(x)$. Stage Theory with Urelements, \STU, now has six axioms:\footnote{As in footnote \ref{fn:lt:cheap}: \STU gives us a stage $\stage{s}$ `for free', so that $\Setabs{x}{\urpred(x)}$ exists by  \ref{stu:spec}.}
\begin{listaxiom}
	\labitem{Empty-U}{u:empty} 
	$(\forall u : \urpred)\forall x\ \, x \notin u$
	\labitem{Ext-U}{u:ext} $(\forall a : \setpred)(\forall b : \setpred)( \forall x(x \in a \liff x \in b) \lonlyif a = b)$
	\labitem{\ref{st:ord}}{dummy1} $\forall \stage{r} \forall \stage{s}\forall \stage{t}(\stage{r} < \stage{s} < \stage{t} \lonlyif \stage{r} < \stage{t})$
	\labitem{Staging-U}{stu:stage} $(\forall a : \setpred)\exists \stage{s}\, \ a \foundat \stage{s}$
	\labitem{Priority-U}{stu:pri} 
	$\forall \stage{s}(\forall a : \setpred)(a \foundat \stage{s} \lonlyif (\forall x \in a)(\urpred(x) \lor x \foundby \stage{s}))$
	\labitem{Spec-U}{stu:spec} 
	$\forall F \forall \stage{s}((\forall x : F)(\urpred(x) \lor x \foundby \stage{s}) \lonlyif (\exists a : \setpred)(a \foundat \stage{s} \land \forall x(F(x) \liff x \in a)))$
\end{listaxiom}
\ref{u:empty} says that no urelement has any members; the other axioms relativise \ST to sets. As in \S\ref{s:1:st}, any cumulative hierarchy obviously satisfies \STU, on the assumption that the urelements are all `always' available to be arbitrarily collected into sets.

We obtain Level Theory with Urelements, \LTU, by tweaking \LT's key definitions. Specifically, I offer the following re-definition:\footnote{The first level is therefore $\Setabs{x}{\urpred(x)}$. This follows \textcites[139]{Montague:STHOL}[16]{Potter:S}[41]{Potter:STP}. By contrast, \citepossess[214]{Scott:AST} first level is $\emptyset$, and the urelements are members of every subsequent level.}
\begin{define}[for \S\ref{s:1:app:ur:conventional} only] 	\label{def:ltu:lev}
	Say that $a$ is potent iff $\forall x((\urpred(x) \lor (\exists c: \setpred)x \subseteq c \in a) \lonlyif x \in a)$. 
	Let $\pot{a} \coloneq \Setabs{x}{\urpred(x) \lor (\exists c: \setpred)x \subseteq c \in a}$, if it exists. Say that $\histpred(h)$ iff $(\forall x \in h)x = \pot(x \cap h)$. Say that $\levpred(s)$ iff $(\exists h : \histpred)s = \pot{h}$. 
\end{define}\noindent
The axioms of \LTU are then \ref{u:empty}, \ref{u:ext}, \ref{lt:strat} (with `$\levpred$' as redefined) and:\footnote{Together, \ref{lt:strat} and \ref{u:sep} deliver the existence of $\Setabs{x}{\urpred(x)}$; see the previous two footnotes.}
\begin{listaxiom}
	\labitem{Sep-U}{u:sep} 
	$\forall F\forall a(\exists b : \setpred)\forall x(x \in b \liff (F(x) \land x \in a))$
\end{listaxiom}
The proofs of \S\S\ref{s:1:LTwo}--\ref{s:1:LTST} go now through with trivial changes. Specifically, the (redefined) levels are well-ordered, and \STU and \LTU make exactly the same demands on sets and urelements. 

The (quasi-)categoricity results of \S\ref{s:1:quasicat} also carry over to \LTU. Let $\model{A}$ and $\model{B}$ be models of \LTU in full second-order logic, and suppose there is a bijection between their respective collections of urelements, $\urpred^\model{A}$ and $\urpred^\model{B}$. This bijection can be lifted to a quasi-isomorphism: $\model{A}$ and $\model{B}$ are isomorphic `so far as they go', but the levels of one may outrun the other. This external result can also be `internalised', yielding results analogous to Theorems \ref{thm:LTInternal} and \ref{thm:LTInternalFull}. 

Note that \LTU, like \LT before it, takes no stance on the height of the hierarchy. In particular, it has no version of Replacement. In this regard, \LTU differs sharply from ZF(C)U, which is something like the `industry standard' for iterative set theory with urelements. It is particularly noteworthy that \LTU allows that the set of urelements may be larger than any pure set.\footnote{\LTU could therefore be used in place of e.g.\ Menzel's ZFCU$'$ \parencite*[67--71]{Menzel:WS}, which is designed to accommodate the claim that the set of urelements is not equinumerous with any pure set.} (For a trivial example, suppose there are exactly 3 urelements and exactly 2 levels; for a less trivial example, suppose there are exactly $\beth_{\omega+1}$ urelements but only an $\omega+\omega$ sequence of levels.) 

\section{Adding absolutely infinitely many urelements}\label{s:urelements:abs}
The Urelemental Story accommodates urelements in a humdrum way. However, there has been recent interest in a less humdrum approach, according to which there are \emph{absolutely infinitely many} urelements. Here is a brisk, three-premise argument in favour of that approach, inspired by Christopher Menzel:\footnote{\textcite[41ff]{Menzel:IEP}; cf.\ \textcite[271--5]{Rumfitt:BST}. \textcite[57]{Menzel:WS} also offers a second (very different) argument to the same conclusion.}
\begin{listl-0}
		\item\label{ord:inf} There are absolutely infinitely many levels in the cumulative hierarchy.
		\item\label{ord:corresp} There are at least as many ordinals as there are levels in the cumulative hierarchy.	
		\item\label{ord:sui} Ordinals are not really sets; they are urelements.
	\end{listl-0}
Each premise is not implausible,\footnote{Claim \eqref{ord:inf} can be motivated by a principle of plenitude concerning sets. Claim \eqref{ord:corresp} can be motivated by combining the fact that the levels of any (pure) cumulative hierarchy are well-ordered (see \S\ref{s:1:wodiscussion}) with the idea that any system of well-ordered objects exemplifies an ordinal (provided that the objects are all members of some set). Claim \eqref{ord:sui} can be motivated by a kind of platonistic structuralism, according to which ordinals are indeed \emph{objects}, but not \emph{sets}, since sets have structure which is not purely order-theoretic. For the record, I do not subscribe to this kind of platonistic structuralism.} and they jointly  entail that there are absolutely infinitely many urelements. In this appendix, I will explore that idea (without endorsing it).

\subsection{Preliminary motivations and observations}\label{s:urelementsinf:prelim}
There is an immediate technical issue: in this kind of cumulative setting, no set has absolutely infinitely many members.\footnote{\emph{Pace} \textcites[44--51]{Menzel:IEP}[71--9]{Menzel:WS}. Note that my argument does not involve Powersets (which Menzel ultimately rejects). Menzel escapes formal inconsistency, whilst retaining (a first-order version of) \ref{sep}, only because his set-theoretic object language has no way to pick out a suitable map, $P$, which witnesses the absolute infinity of his set $\Setabs{x}{\urpred(x)}$.} This follows from a simple version of Cantor's Theorem. For reductio, suppose that some set, $a$, has absolutely infinitely many members. As discussed in \S\ref{s:1:quasicat}, this entails that $\somany x\ x \in a$, i.e.\ there is a map, $P$, such that $\forall x\ P(x) \in a$ and $(\forall y \in a)\exists ! x\ P(x) = y$. By $P$'s injectivity and \ref{sep},\footnote{I take it that rejecting \ref{sep} is not an option in this setting; though see Pt.\ref{pt:blt} for an approach which rejects \ref{sep}.} there is some $d = \Setabs{x \in a}{x \notin P^{-1}(x)}$. Since $P(d) \in a$, contradiction follows familiarly.

So: if there are absolutely infinitely many urelements, then there is no set of all urelements.\footnote{\textcite[330--1]{Uzquiano:NRRPP} also suggests the use of a set theory with urelements but no set of urelements, though for somewhat different reasons.} But the existence of such a set is a trivial consequence of \ref{stu:spec}, as laid down in \S\ref{s:1:app:ur:conventional}. So, those who think that there are absolutely infinitely many urelements must reject \ref{stu:spec}. Furthermore, since \ref{stu:spec} follows from the third sentence of the Urelemental Story of \S\ref{s:1:app:ur:conventional}, they must change their story. 

Many alternative stories are possible, but the simplest approach is simply to bolt a Limitation of Size principle onto the Urelemental Story, insisting that the Basic Story remains correct of the pure sets, whilst denying that any set is absolutely infinite. This leads to the following:\footnote{Cf.\ \textcite[331]{Uzquiano:NRRPP}.}
\begin{storytime}\textbf{The \Uinf{}relemental Story.}
	Sets are arranged in stages. Every set is found at some stage. At any stage $\stage{s}$: for any things---provided both that (i) there are not absolutely infinitely many of them, and that (ii) each of them is either a set found before $\stage{s}$ or an urelement---we find a set whose members are exactly those things. We find nothing else at $\stage{s}$.  (NB: since the Basic Story is correct of the pure sets, we do not find {absolutely} infinitely many pure sets before $\stage{s}$.)
\end{storytime}\noindent
In the remainder of this appendix, I will briefly sketch (equivalent) stage-theoretic and level-theoretic formalizations of this Story. For readability, I leave all proofs to the reader, with hints in footnotes.

\subsection{Stage-theoretic approach: \STUfew}\label{s:stufew}
To axiomatize the \Urelementalinf Story, we need a predicate, `$\purepred$', to pick out the pure sets (cf.\ \S\ref{s:1:quasicat}). Since we have assumed that the Basic Story holds of the pure sets, we can define `$\purepred$' explicitly:
\begin{define}\label{def:purepred}
	Say that $a$ is pure, $\purepred(a)$, iff both $\setpred(a)$ and there is some transitive $c \supseteq a$ whose members are all sets.
\end{define}\noindent
To axiomatize the \Urelementalinf Story, we also need a way to formalize `there are absolutely infinitely many $\Phi$s'. There are familiar concerns about the possibility of formalizing this idea.\footnote{See e.g.\ \textcite[279]{McGee:TPTTC}.} Nonetheless, if there are absolutely infinitely many $\Phi$s, then certainly $\somany x \Phi(x)$ (cf.\ \S\ref{s:1:quasicat}). Conversely, if $\somany x \Phi(x)$, then no property can have \emph{more} instances than $\Phi$. So, `$\somany x \Phi(x)$' will serve as our proxy for `there are absolutely infinitely many $\Phi$s'.\footnote{\label{fn:somanyduck}Very little of what I say depends upon this particular choice of proxy. In particular, I rely upon its logical properties only when claiming that both \STUfew and \LTUfew prove \ref{u:sep}, and in my remarks on the quasi-categoricity of \LTUfew.}

I can now lay down the theory \STUfew. Its axioms are \ref{u:empty}, \ref{u:ext}, \ref{st:ord}, \ref{stu:stage}, \ref{stu:pri}, and the following: 
\begin{listaxiom}
	\labitem{Spec-\Uinf}{stfew:spec} 
	$\forall F \forall \stage{s}((\lnot\somany x F(x) \land (\forall x : F)(\urpred(x) \lor x \foundby \stage{s}))\lonlyif {}$\\
	$\phantom{\forall F\forall\stage{s}(}(\exists a : \setpred)(a \foundat \stage{s} \land \forall x(F(x) \liff x \in a)))$
	\labitem{LoS-\Uinf}{set:few} $(\forall a : \setpred)\lnot\somany x\ x \in a$
	\labitem{Pure-\Uinf}{pure:few} $\forall F \forall \stage{s}((\forall x : F)(\purepred(x) \land x \foundby \stage{s}) \lonlyif \lnot\somany x F(x))$
	\labitem{Many-\Uinf}{u:many} $\somany x \urpred(x)$
\end{listaxiom}
In brief: \ref{stfew:spec} restricts \ref{stu:spec} to capture conditions (i)--(ii) of the \Urelementalinf Story; \ref{set:few} enshrines Limitation of Size, which follows from condition (i) plus the fact that `we find nothing else' at any stage; \ref{pure:few} formalizes the parenthetical `NB' of the Story; and \ref{u:many} formalizes the claim that there are absolutely infinitely many urelements.

\subsection{Level-theoretic approach: \LTUfew}\label{s:ltufew}
\STUfew is a multi-sorted, stage-theoretic, formalization of the \Urelementalinf Story. That Story can instead be given a single-sorted formalization, \LTUfew. To do this, I start by tweaking \LT's key definitions:
\begin{define}[for \S\ref{s:urelements:abs} only] 	\label{def:levelininfcontext} Say that $a$ is potent iff $(\forall x : \setpred)(\exists c(x \subseteq c \in a) \lonlyif x \in a)$. Let $\pot{a} = \Setabs{x}{\setpred(x) \land \exists c(x \subseteq c \in a)}$, if it exists. Say that $\histpred(h)$ iff $(\forall x \in h)x = \pot{(x \cap h)}$. Say that $\levpred(s)$ iff $(\exists h : \histpred)s = \pot{h}$.
\end{define}\noindent
Using these redefinitions, we can prove analogues of Lemmas \ref{lem:es:order}--\ref{lem:es:comparability} from \S\ref{s:1:LTwo}. Specifically, given \ref{u:ext} and \ref{u:sep}, we can prove that the levels (so defined) are potent, transitive, pure,\footnote{Since they are transitive, they witness their own purity.} and well-ordered by $\in$.

I can now lay down \LTUfew. It uses a  primitive one-place function symbol, $\levpof$, where `$\levpof{a}$' should be read as \emph{$a$'s level-index}. (I discuss the use of this primitive in \S\ref{s:urelements:ltuinfc}.) Then \LTUfew has six axioms: \ref{u:empty}, \ref{u:ext},  \ref{set:few}, \ref{u:many}, and two axioms governing $\levpof$:
\begin{listaxiom}
	\labitem{Leveller}{leveller} 
	$(\forall a : \setpred)((\exists s : \levpred)\levpof{a} = s \land {}$\\
	$\phantom{(\forall a : \setpred)(}(\forall x : \setpred)(x \in a \lonlyif \levpof{x} \in \levpof{a}) \land {}$\\
	$\phantom{(\forall a : \setpred)(}(\forall s : \levpred)(s \in \levpof{a} \lonlyif (\exists x : \setpred)(x \in a \land s \subseteq \levpof{x})))$
	\labitem{Consolidation}{consolidation}
	$\forall F ((\lnot\somany x F(x) \land \exists a(\forall x : F)(\urpred(x) \lor \levpof{x} \in a)) \lonlyif {}$\\
	$\phantom{\forall F (}(\exists b : \setpred)\forall x(F(x) \liff x \in b))$		
\end{listaxiom}\noindent
To understand these axioms, note that \LTUfew guarantees that the (pure) levels are well-ordered by membership.\footnote{To see this, note \LTUfew proves \ref{u:sep}, and combine this with the remarks after Definition \ref{def:levelininfcontext}.} Now, \ref{leveller} ensures that the levels index the sets; intuitively, $a$'s level-index is the least level greater than the level-index of every set in $a$. \ref{consolidation} then allows us to find all the impure sets we would want to find `at' any given level. Finally, note that \LTUfew proves a pure-analogue of \ref{lt:strat}:\footnote{Use induction on levels, together with the second conjunct of \ref{leveller}.}
\begin{lem}[\LTUfew]\label{lem:ltphi:purelev}
	If $a$ is pure, then $a \subseteq \levpof{a}$
\end{lem}\noindent
Consequently, \LTUfew's pure sets can be thought of as satisfying \LT. Indeed, if we define `$x \varin y$ as `$\purepred(x) \land \purepred(y) \land x \in y$', then $\LTUfew \proves \LT(\purepred, \varin)$, as defined in \S\ref{s:1:quasicat}. It follows that \LTUfew is externally and internally (quasi-)categorical: any two hierarchies satisfying \LTUfew have quasi-categorical pure sets; moreover, if there is a bijection between the hierarchies' urelemental bases, their impure sets are quasi-categorical. (However, \LTUfew's analogue of Theorem \ref{thm:LTexternalcat} is more restricted: if $\model{M}$ is a standard, set-sized, model of \LTUfew, then $|\urpred^\model{M}|$ is regular.)\footnote{Assuming Choice. \emph{Proof.} Let $\kappa = |\urpred^\model{M}|$. By \ref{consolidation}, every smaller-than-$\kappa$ subset of $\urpred^\model{M}$ is in $\setpred^\model{M}$. So $\kappa$ is infinite, by \ref{u:many}. For each $\lambda < \kappa$, there are $\kappa^\lambda$ subsets of $\urpred^\model{M}$ with cardinality $\lambda$, so that $\kappa^\lambda \leq |\setpred^\model{M}|$. So if $\cofinality{\kappa} < \kappa$, then by K\"onig's Theorem $\kappa < \kappa^{\cofinality{\kappa}} \leq |\setpred^\model{M}|$, contradicting \ref{u:many}; hence $\cofinality{\kappa} = \kappa$. (Thanks to Gabriel Uzquiano for suggesting I consider how \LTUfew interacts with regular cardinals.)}

In fact, \LTUfew and \STUfew are provably equivalent, concerning sets and urelements. To prove that \LTUfew interprets \STUfew, tweak the $*$-translation of \S\ref{s:1:LTST}, so that $(x \foundat \stage{s})^* \coloneq \levpof{x} \subseteq \stage{s}$.\footnote{Stipulate that $(\setpred(x))^* \coloneq \setpred(x)$.} It is then easy to show that $\LTUfew \proves \STUfew^*$ (cf.\ Lemma \ref{lem:lt:sttrans}).

To show that \STUfew interprets \LTUfew, first note that  \STUfew proves \ref{u:sep} and the converse of \ref{stu:pri} (cf.\ Lemmas \ref{lem:st:sep}--\ref{lem:st:conversepri}). Then tweak Definition \ref{def:slice} (cf.\ Definition \ref{def:levelininfcontext}):
\begin{define}[for \S\ref{s:urelements:abs} only] \label{def:sliceininfcontext}
	Let $\slice{s} \coloneq \Setabs{x \foundby \stage{s}}{\purepred(x)}$. Say that $a$ is a \emph{slice} iff $a = \slice{s}$ for some $\stage{s}$.
\end{define}\noindent
It follows that the slices are the levels (in the senses of Definitions \ref{def:levelininfcontext} and \ref{def:sliceininfcontext}; cf.\ Lemma \ref{lem:st:levelsslices}).\footnote{For the analogue of Lemma \ref{lem:st:slice}: \ref{u:ext}, \ref{pure:few}, and \ref{stfew:spec} guarantee that $\slice{s}$ exists for each stage $\stage{s}$; in clauses \eqref{slice:foundup}--\eqref{slice:foundat}, the quantifier `$\forall a$' becomes `$(\forall a : \purepred)$'; and note that each slice witnesses its own purity.} 
We can then interpret \LTUfew's unique primitive, $\levpof$, via $\rankof{}$, defined as follows:\footnote{Since \STUfew does not prove that all stages are comparable (cf.\ the discussion of Boolos's \ref{boolos:net} from \S\ref{s:1:history:boolos}), it takes several steps to vindicate Definition \ref{def:stuinf:rank}. First: show that stages obey $<$-induction. Second: show that if $a \foundat \stage{s}$ and $\lnot\exists\stage{r}(a \foundat{r} < \stage{s})$ and $a \foundat \stage{t}$ and $\lnot\exists\stage{r}(a \foundat{r} < \stage{t})$ then $\slice{s} = \slice{t}$; it follows that $\rankof{a}$ is a slice. Third: combine this with the fact that the slices are levels, to show that $\rankof$ behaves like $\levpof$.}
\begin{define}\label{def:stuinf:rank}
	For each set $a$, let $\rankof{a} \coloneq \bigcap\Setabs{\slice{s}}{a \foundat \stage{s} \land \lnot \exists \stage{r}(a \foundat \stage{r} < \stage{s})}$.
\end{define}
\begin{thm}\label{thm:LTUinfpSTUinf}
	$\STUfew \proves \phi^\rankof{}$ iff $\LTUfew \proves \phi$, for any \LTUfew-sentence $\phi$, where $\phi^\rankof{}$ is the formula which results from $\phi$ by replacing each instance of $\levpof$ with $\rankof{}$.
\end{thm}\noindent
The upshot is that no information about sets or urelements is lost or gained in moving from \STUfew to \LTUfew. Since any hierarchy which is described by the \Urelementalinf Story satisfies \STUfew, it also satisfies \LTUfew. And \LTUfew is quasi-categorical. Our work on the \Urelementalinf Story is complete. 

\subsection{Eliminating primitives and first-orderisation}\label{s:urelements:ltuinfc}
Or rather: almost complete. Given the discussion of \S\ref{s:1:st}, we may want to eliminate \LTUfew's primitive, $\levpof$. This is easily done within second-order logic: just conjoin \ref{leveller} and \ref{consolidation}, and bind $\levpof$ with a (second-order) existential quantifier. But if we are willing to make some further assumptions, then we can eliminate $\levpof$ using certain \emph{first-order} functions.\footnote{\textcites[1047]{LevyVaught:PPR}[299]{Uzquiano:MSOZST} present a somewhat similar method for defining the rank of a set (via functions on ordinals). Here I treat first-order functions as sets of ordered pairs in the normal way, and $x \in \domain{f}$ abbreviates $\exists y\, \tuple{x, y} \in f$. Of course, a fully first-order version of \LTUfew would need to define `$\somany x F(x)$' differently (cf.\ footnote \ref{fn:somanyduck}).} 

Roughly, a \emph{ranking-function}: \eqref{frank:trans} has a transitive domain (setting aside urelements); and \eqref{frank:llike} behaves like $\levpof$ where defined. More formally:
\begin{define}\label{def:frank}
	Say that a function $f$ is a \emph{ranking-function} iff, for all $a \in \domain{f}$, both: 
	\begin{listn-0}
		\item\label{frank:trans} $\setpred(a)$ and $(\forall x : \setpred)(x \in a \lonlyif x \in \domain{f})$; and 
		\item\label{frank:llike} $\levpred(f(a))$ and $(\forall x : \setpred)(x \in a \lonlyif f(x) \in f(a))$ and $(\forall s : \levpred)(s \in f(a) \lonlyif (\exists x \in a)s \subseteq f(x))$. 
	\end{listn-0}
	Say that $\rankspred{f}{a}$ iff $f$ is a ranking-function with $a \in \domain{f}$.
\end{define}\noindent
It is easy to show that ranking-functions agree wherever they are defined, i.e.:
\begin{lem}[\ref{u:ext}, \ref{u:sep}]\label{lem:ltuinfc:rankersagree}
	If $\rankspred{f}{a}$ and $\rankspred{g}{a}$, then $f(a) = g(a)$.
\end{lem}\noindent
We can now replace \ref{leveller}, in \LTUfew, with $(\forall a : \setpred)\exists f\, \rankspred{f}{a}$. Note that this claim is \emph{independent} of \LTUfew: it guarantees that every set is a member of some set, and so guarantees that the hierarchy has no final stage  (cf.\ \ref{lt:cre} from \S\ref{s:1:ltsubzf}). Still, this allows us to define $\levpof{a} \coloneq \bigcap \Setabs{f(a)}{\rankspred{f}{a}}$. We can use this definition in \ref{consolidation}, and prove \ref{leveller} via Definition \ref{def:frank}.

\section*{Acknowledgements}
Thanks to Neil Barton, Sharon Berry, Luca Incurvati, Juliette Kennedy, \O{}ystein Linnebo, Michael Potter, Chris Scambler, Will Stafford, James Studd, Rob Trueman, Gabriel Uzquiano, Sean Walsh, and anonymous referees for \emph{Bulletin of Symbolic Logic}.

\stopappendix

\chapter[Part 2]{Level Theory, Part \thechapter
	\chapsubhead{Axiomatizing the bare idea of a potential hierarchy}}\label{pt:pst}

\noindent\textcolor{blue}{This document contains preprints of Level Theory, Parts 1--3. All three papers are forthcoming at \emph{Bulletin of Symbolic Logic}.}

\begin{quote}
	\textbf{Abstract.} Potentialists think that the concept of set is importantly modal. Using tensed language as a heuristic, the following bare-bones story introduces the idea of a potential hierarchy of sets: `Always: for any sets that existed, there is a set whose members are exactly those sets; there are no other sets'. Surprisingly, this story already guarantees well-foundedness
	and persistence. Moreover, if we assume that time is linear, the ensuing modal set theory is almost definitionally equivalent with non-modal set theories; specifically, with Level Theory, as developed in Part \ref{pt:lt}.
\end{quote}

\begin{epigraph}
	{What we need to do is to replace the language
			of time and activity by the more bloodless language of potentiality and
			actuality.}
	{\textcite[293]{Parsons:ICS}}
\end{epigraph}
\noindent 
Potentialists, such as Charles Parsons, \O{}ystein Linnebo, and James Studd, think that the concept of \emph{set} is importantly modal. Put thus, potentialism is a broad church; different potentialists will disagree on the precise details of the relevant modality.\footnote{See e.g.\ \textcites{Fine:RUQ}[209]{Linnebo:PHS}[264--5]{Linnebo:PMML}[61--5]{Linnebo:TO}[706--7]{Studd:ICS}[144--53]{Studd:EML}.} My aim is shed light on potentialism, in general, using Level Theory, \LT, as introduced in Part \ref{pt:lt}.

I start by formulating Potentialist Set Theory, \PST. This uses a tensed logic to formalize the bare idea of a `potential hierarchy of sets'.\footnote{This is \citepossess{Linnebo:PHS} phrase.} Though \PST is extremely minimal, it packs a surprising punch (see \S\S\ref{s:2:story}--\ref{s:2:wodiscussion}).

In the vanilla version of \PST, we need not assume that time is linear. However, if we make that assumption, then the resulting theory is almost definitionally equivalent to \LT, its non-modal counterpart (see \S\S\ref{s:2:LPST}--\ref{s:2:ns:2}). This equivalence allows me to clarify Hilary Putnam's famous claim, that modal and non-modal set theories express the same facts (see \S\ref{s:2:equivalencethesis}). Putting my cards on the table: I am not a potentialist, in part because I am so sympathetic with Putnam's claim.

This paper presupposes familiarity with Part \ref{pt:lt}. My notation conventions are as in Pt.\ref{pt:lt} \S\ref{s:1:prelim},  with the addition that I use $\vec{x}$ for an arbitrary sequence, writing things like $F(\vec{x})$ rather than $F(x_1, \ldots, x_n)$. For readability, all proofs are relegated to the appendices.

\section{Tense and possibility}\label{s:2:story}
Many potentialists hold that temporal language serves as a useful heuristic for their favoured mathematical modality. To illustrate the idea, consider what Studd calls the \emph{Maximality Thesis}: `any sets \emph{can} form a set.'\footnote{\textcite[699]{Studd:ICS}. \textcite[206--8]{Linnebo:PHS} formulates a similar thesis.} This Thesis is given a modal formulation. But, as Studd notes, it can be glossed temporally: `any sets \emph{will} form a set'. Of course, no potentialist will take this temporal gloss literally. Nobody, after all, wants to countenance absurd questions like `which pure sets existed at noon today?', or `which pure sets will exist by teatime?'\footnote{For further issues, see e.g.\ \textcites[\S{}II]{Parsons:ICS}[706]{Studd:ICS}[49]{Studd:EML}.} The idea, to repeat, is just that temporal language is a useful \emph{heuristic} for the potentialist's preferred modality.

To elaborate on this heuristic, consider the bare-bones story of (pure) sets, which I told and explored in Part \ref{pt:lt}, and which I will repeat here:
\begin{storytime}\textbf{The Basic Story.}
	Sets are arranged in stages. Every set is found at some stage. At any stage $\stage{s}$: for any sets found before $\stage{s}$, we find a set whose members are exactly those sets. We find nothing else at $\stage{s}$.
\end{storytime}\noindent
We can regard the stages of this Story as moments of time. Regarded thus, the Basic Story adopts the \emph{tenseless} view of time, according to which moments are just a special kind of object. But this tenseless approach serves potentialists poorly. At no stage is there a set of all the sets which are found at any stage, so this tenseless Story falsifies the claim `any sets \emph{will} form a set'. 

Familiarly, though, time can also be thought of in a 
\emph{tensed} fashion. On the tensed approach, we do not quantify over moments or stages; rather, we use primitive temporal operators, like `it \emph{was} the case that\ldots' or `previously: \ldots'. And we can retell the Basic Story in tensed terms:
\begin{storytime}\textbf{The Tensed Story.}
	Always: for any sets that {existed}, there {is} a set whose members are exactly those sets; there {are} no other sets. 
\end{storytime}\noindent
Unlike the Basic Story, this Tensed Story is compatible with the claim `(always:) any sets \emph{will} form a set'.

Note, though, that I say `is compatible with', rather than `entails'. If time abruptly ends, then some things will never form a set. And, by design, the Tensed Story is compatible both with the claim that time abruptly ends, and that time is endless. Otherwise put: it says nothing at all about the `height' of any potential hierarchy. This silence is deliberate, for potentialists might disagree about questions of `height'. 

Still, once potentialists have agreed to use tense as an heuristic for their preferred modality, I do not see how they could doubt that the Tensed Story holds of every potential hierarchy of sets. In what follows, then, I take it for granted that the Tensed Story presents us with the \emph{bare idea} of a potential hierarchy

\section{Temporal logic for past-directedness}\label{s:2:pdlogic}
My first goal is to axiomatize the Tensed Story. For this, I will employ a temporal logic. In particular, I use a negative free second-order logic which assumes that time is past-directed. Here is a brief sketch of this \emph{past-directed-logic}, with fuller explanations in footnotes. (Let me take this opportunity to flag that I am wholly indebted to Studd for the idea of investigating potentialism via temporal logic; see \S\ref{s:2:studd}.) 

We use `$\existspred{x}$' as an existence predicate; it abbreviates `$x = x$'. We prohibit consideration of never-existent entities, and we insist that quantification and atomic truth require existence.\footnote{\label{fn:comprehensionisvanilla}\label{fn:negativefreelogic}So, we adopt the axiom scheme $\posssub\existspred{x}$, and inference rules so that these schemes hold: 
	(1) $\exists x \phi \lonlyif \exists x(\existspred{x} \land \phi)$ and $\exists F \phi \lonlyif \exists F(\existspred{F} \land \phi)$, for any formula $\phi$; 
	(2) $\alpha(\vec{x}) \lonlyif (\existspred{x_1} \land \ldots \land \existspred{x_n})$, for any atomic $\alpha(\vec{x})$ with all free variables displayed; 
	(3) $F(\vec{x}) \lonlyif \existspred{F}$ for any `$F$'. We also have second-order Comprehension, i.e.\ the scheme $\exists F\forall \vec{x}(F(\vec{x}) \liff \phi)$, for any $\phi$ not containing `$F$'.}
We have three temporal operators (with their obvious duals):\footnote{i.e.\ $\allearliersub \phi  \coloneq \lnot \someearliersub \lnot \phi$ and $\alllatersub \phi \coloneq \lnot \somelatersub \lnot \phi$ and $\necsub \phi \coloneq \lnot \posssub \lnot \phi$.}
\begin{listbullet}
	\item[$\someearlier$:] A past-tense operator; gloss `$\someearlier\phi$' as `previously: $\phi$' or `it \emph{was} the case that $\phi$' . 
	\item[$\somelater$:] A future-tense operator; gloss `$\somelater\phi$' as `eventually: $\phi$'  or `it \emph{will be} the case that $\phi$'.
	\item[$\poss$:] An unlimited temporal operator; gloss `$\poss\phi$' as `sometimes: $\phi$'.
\end{listbullet}
We have Necessitation rules: if $\phi$ is a theorem, then so are both $\allearlier\phi$ and $\alllater\phi$. We then lay down schemes as follows:\footnote{See e.g.\ \textcite[41]{Goldblatt:LTC} for all but the last scheme.}
\begin{align*}
	&\allearlier(\phi \lonlyif \psi) \lonlyif (\allearlier \phi \lonlyif \allearlier \psi) & 
	&\alllater(\phi \lonlyif \psi) \lonlyif (\alllater \phi \lonlyif \alllater \psi)\\
	&\someearlier \alllater \phi \lonlyif \phi & 
	& \somelater\allearlier \phi \lonlyif \phi\\
	&\allearlier \phi \lonlyif \allearlier\allearlier \phi &
	&\someearlier(\phi \land \allearlier\phi) \lonlyif \allearlier(\phi \lor \someearlier\phi)
\end{align*}
The first two schemes are familiar distribution principles. The second two schemes ensure appropriate past/future interaction. The fifth scheme ensures that \emph{before} is transitive. The last scheme characterizes \emph{past-directedness}.\footnote{\label{fn:def:pathconnected}i.e.\ this frame-condition: $(\forall \worldvar{v} \leq \worldvar{w})(\forall \worldvar{u} \leq \worldvar{w})(\exists \worldvar{t} \leq \worldvar{v})\worldvar{t} \leq \worldvar{u}$. Equivalently: if $\worldvar{v}$ and $\worldvar{u}$ are path-connected, then $(\exists \worldvar{t} \leq \worldvar{v})\worldvar{t} \leq \worldvar{u}$. We say that worlds are \emph{path-connected} iff they are related by the reflexive, symmetric, transitive closure of the accessibility relation.} Given past-directedness, `sometimes: $\phi$' amounts to `it was, is, will be, or was going to be the case that $\phi$'. So we adopt this scheme:
	$$\poss \phi \liff (\someearlier \phi \lor \phi \lor \somelater \phi \lor \someearlier\somelater\phi)$$
It follows that $\poss$ obeys S5.  This completes my sketch of past-directed-logic. 

In what follows, I will assume that potentialists are happy to use this temporal logic.\footnote{Though note that not all potentialists have used temporal logics; see \S\ref{s:parsonslinnebo}.} However, it is worth repeating that our potentialist only regards time as an \emph{heuristic}. Ultimately, they want $\poss$ to express their favoured mathematical modality. So they will need to explain how (and why) their favoured modality decomposes into other operators, $\someearlier$ and $\somelater$, which obey past-directed-logic. This is a non-trivial demand; but, for the purposes of this paper, I assume it can be met.

\section{Potentialist Stage Theory}\label{s:2:PST}
Armed with past-directed-logic, the Tensed Story is easy to axiomatize. Let \PST, for Potentialist Set Theory, be the result of adding these four axioms to past-directed-logic:
\begin{listaxiom}
 	\labitem{Mem{$_{\posssub}$}}{m:mem} 
	$\forall a \nec \forall x(\poss x \in a \lonlyif \nec(\existspred{a} \lonlyif x \in a))$
	\labitem{Ext{$_{\posssub}$}}{m:ext}
	$\forall a \nec \forall b(\nec \forall x(\poss x \in a \liff \poss x \in b) \lonlyif \poss a = b)$
	\labitem{Priority{$_{\someearliersub}$}}{pst:pri} 
	$\forall a (\forall x \in a) \someearlier \existspred{x}$
	\labitem{Spec{$_{\someearliersub}$}}{pst:spec}
	$\forall F((\forall x : F)\someearlier\existspred{x} \lonlyif \exists a  \forall x(F(x) \liff x \in a))$
\end{listaxiom}\noindent
The first two axioms are not explicit in the Tensed Story, but I take it they are supposed to be something like analytic: roughly, \ref{m:mem} says that each set $a$ has its members essentially, and \ref{m:ext} says that if everything which could (ever) be in $a$ could be in $b$, and vice versa, then $a = b$ (when they exist).\footnote{See \textcites[286(3)]{Parsons:ICS}[711--12]{Studd:ICS}[]{Studd:EML}[215]{Linnebo:PHS}[211--2]{Linnebo:TO}.} The next two axioms are explicit in the Story: \ref{pst:pri} says that a set's members existed before the set itself, and \ref{pst:spec} says that, if every $F$ existed earlier, then the set of $F$s exists.  So all of \PST's axioms are obviously true of the Tensed Story. 

It is worth comparing \PST with Stage Theory, \ST (see Pt.\ref{pt:lt} \S\ref{s:1:st}). Indeed, we could equally think of \PST as Potentialis\emph{ed Stage} Theory, since it is little more than the most obvious reworking of \ST using tensed operators.\footnote{But \PST is indeed \emph{more}: \PST assumes past-directedness, and \ST has no comparable assumption about stages. (Cf.\ the discussion of \citepossess{Boolos:IA} Net in Pt.\ref{pt:lt} \S\ref{s:1:history:boolos}.) For the technical role of past-directedness, see the end of \S\ref{s:2:app:pst}.}

\section{The inevitability of well-foundedness and persistence}\label{s:2:wodiscussion}
I have just shown that \PST is a good formalization of the Tensed Story. As explained in \S\ref{s:2:story}, though, this Story articulates the \emph{bare idea} of a potential hierarchy of sets. It follows that any potential hierarchy satisfies \PST. This is significant, since \PST is surprisingly rich.

To gauge \PST's depths, I will explain how it relates to Level Theory, \LT, the non-modal theory which axiomatizes the (tenseless) Basic Story (see Pt.\ref{pt:lt} \S\S\ref{s:1:st}--\ref{s:1:wodiscussion}). According to \LT, the sets are arranged into well-ordered \emph{levels}, where levels are sets which goes proxy for the stages of the Basic Story. Now, \PST proves the following result (see \S\ref{s:2:app:pst}):
\begin{thm}[\PST]\label{thm:pst:key}
	Where $\maxlev(s)$ abbreviates $(\existspred{s} \land \forall x\ x \subseteq s)$:
	\begin{listn-0}
		\item\label{key:lt} $\LT$ holds
		\item\label{key:persistence} $\forall x \alllater\existspred{x}$
		\item\label{key:maxlev} $(\exists s : \levpred)\maxlev(s)$
		\item\label{key:prune} $(\forall s : \levpred)\poss\maxlev(s)$
	\end{listn-0}
\end{thm}\noindent
If we consider a Kripke model of \PST: \eqref{key:lt} says that every possible world comprises a hierarchy of sets, arranged into well-ordered levels. Among other things, this yields \emph{well-foundedness}, i.e.\ $\forall F(\exists x F(x) \lonlyif (\exists x : F) (\forall z : F)z \notin x)$.\footnote{Indeed, \PST proves a modal version of well-foundedness; see Lemma \ref{cor:pst:lob}.} Then \eqref{key:persistence} is a statement of \emph{persistence}; it says that, once a set exists, it exists forever after. Last, \eqref{key:maxlev} says that every world has a maximal level, and \eqref{key:prune} says that every level is some world's maximal level. So, the worlds in a Kripke model of \PST are, in effect, just arbitrary, persistent, initial segments of an (actualist) \LT-hierarchy of pure sets. 

I will develop the link between \PST and \LT over the next few sections. First, I want to highlight the significance of Theorem \ref{thm:pst:key}. The Tensed Story does not involve an explicit statement of well-foundedness or persistence. So one \emph{might} try to entertain versions of the Tensed Story wherein well-foundedness or persistence fail: that is, one might try to entertain a potential hierachy wherein time had no beginning, or wherein sets fade in and out of existence. But the foregoing remarks show that all such speculation is incoherent: every potentialist hierarchy \emph{must} obey well-foundedness and persistence, since every potentialist hierarchy obeys \PST, and \PST proves Theorem \ref{thm:pst:key}. Echoing Scott, then, we see `how little choice there is in setting up' a potential hierarchy of sets.\footnote{\textcite[210]{Scott:AST}. That quote is discussed in Pt.\ref{pt:lt} \S\ref{s:1:wodiscussion}; this section `modalizes' that discussion.}

\section{Linear Potentialist Stage Theory}\label{s:2:LPST}
So far, our potentialist has assumed that time is past-directed (to use the tensed-heuristic). If we also assume that time is \emph{linear}, then we can obtain even deeper connections between \PST and \LT. I will spell out these connections in \S\S\ref{s:2:ns:1:ded}--\ref{s:2:ns:2}; first, I must say a bit about linearity. 

Formally, we can insist on linearity by adding these schemes to past-directed-logic:\footnote{See e.g.\ \textcite[78]{Goldblatt:LTC}. These allow us to prove the schemes $\posssub \phi \liff(\someearliersub \phi \lor \phi \lor \somelatersub \phi \lor \someearliersub \somelatersub \phi)$ and $\someearliersub (\phi \land \allearliersub \phi) \lonlyif \allearliersub (\phi \lor \someearliersub \phi)$ of \S\ref{s:2:pdlogic}.}
\begin{align*}
	\poss \phi &\liff (\someearlier\phi \lor \phi \lor \somelater\phi) & 
	\someearlier\somelater\phi &\lonlyif \poss\phi & 
	\somelater\someearlier\phi&\lonlyif \poss\phi
\end{align*}
As in \S\ref{s:2:pdlogic}, potentialists who want to use this linear-logic must explain why their favoured notion of mathematical possibility vindicates such linearity; this is a non-trivial challenge, but again I will not push it.\footnote{Though cf.\ the discussion of Boolos's \cite*{Boolos:ICS} theory in Pt.\ref{pt:lt} \S\ref{s:1:history:boolos}, and footnote \ref{fn:Studdtroubles}.} When using \PST with this linear logic, I write \LPST, for linear-\PST.
	
By combining Theorem \ref{thm:pst:key} with the assumption of linearity, we can simplify our ideology considerably. Intuitively, linearity allows us to gloss `previously' as `when there are fewer things', and to gloss `eventually' as `when there are more things'. More precisely, we recursively define a translation, $\bullet$, whose only non-trivial clauses are as follows:\footnote{So $(\exists x \phi)\mltint \coloneq \exists x \phi\mltint$, $(\exists X\phi)\mltint \coloneq \exists X \phi\mltint$, $(\poss \phi)\mltint \coloneq \poss \phi\mltint$, $(\lnot\phi)\mltint \coloneq \lnot \phi\mltint$, $(\phi \land \psi)\mltint\coloneq( \phi\mltint\land\psi\mltint)$, and $\alpha\mltint \coloneq \alpha$ for atomic $\alpha$; we choose variables to avoid clashes.}
\begin{align*}
	(\someearlier \phi)\mltint &\coloneq \exists x \poss (\lnot\existspred{x} \land \phi\mltint) \\
	(\somelater\phi)\mltint & \coloneq (\exists x : \maxlev)\poss(\exists v\ x \in v \land \phi\mltint)	
\end{align*}
It is then easy to prove:
\begin{prop}[\LPST]\label{prop:pstl:mltint}
	$\phi \liff \phi\mltint$ for any \LPST-formula $\phi$
\end{prop}\noindent
We can therefore rewrite \LPST, without loss, as a modal theory which uses a \emph{single} primitive modal operator, $\poss$, which obeys S5 (for more, see \S\ref{s:2:app:lpst}).
	
We can go even further, though, and eliminate \emph{all} modal notions from \LPST. The rough idea is straightforward. Theorem \ref{thm:pst:key} says that \emph{levels simulate possible worlds, and vice versa}. By assuming linearity, we can obtain results which say: \emph{actual hierarchies simulate potential hierarchies, and vice versa}. 

That way of putting things is, however, rather rough. The details of the simulation are in fact quite fiddly. I will therefore divide my discussion into three sections. In \S\ref{s:2:ns:1:ded}, I consider a \emph{deductive} version of this simulation. This is suitable for first-order versions of \LT and \LPST, which I call \LTfo and \LPSTfo.\footnote{These arise just by replacing the single second-order axiom, \ref{sep} or \ref{pst:spec}, with its obvious first-order schematisation, and abandoning Comprehension.} In \S\ref{s:2:ns:1:sem}, I consider a \emph{semantic} version of this first-order simulation. Finally, in \S\ref{s:2:ns:2}, I consider deductive and semantic versions of this simulation for (various) \emph{second-order} versions of \LT and \LPST. 

\section{Near-synonymy: first order, deductive}\label{s:2:ns:1:ded}
To interpret \LTfo in \LPSTfo, we will simply replace what \emph{happens} with what \emph{could happen}. More precisely, we consider the following translation; following Studd, I call $\phi\modalize$ the \emph{modalization} of $\phi$:\footnote{\textcites[708]{Studd:ICS}[154]{Studd:EML}; cf.\ also \textcites[115--6]{Linnebo:PS}[213]{Linnebo:PHS}.} 
\begin{align*}
	\alpha\modalize & \coloneq \poss \alpha, \text{for atomic }\alpha & 
	(\phi \land \psi)\modalize &\coloneq (\phi\modalize \land \psi\modalize) \\
	(\lnot \phi)\modalize &\coloneq \lnot\phi\modalize & 
	(\exists x \phi)\modalize &\coloneq \poss \exists x \phi\modalize 
\end{align*}
Conversely, to interpret \LPSTfo in \LTfo, we take the hint suggested by Theorem \ref{thm:pst:key}, and simply regard possible worlds as levels. More precisely, we consider the following translation; I call $\phi\levelling$ the \emph{levelling} of $\phi$:\footnote{\textcites[224--5]{Linnebo:PHS}[719]{Studd:ICS}[173]{Studd:EML} consider similar maps. We choose new variables (to avoid clashes) in the clauses for $(\someearliersub \phi)^s$, $(\somelatersub \phi)^s$ and $(\posssub \phi)^s$.} 
\begin{align*}
	(x = y)\levelling &\coloneq (x = y \subseteq s)
	&
	(x \in y)\levelling &\coloneq (x \in y \subseteq s) \\
	(\phi \land \psi)\levelling &\coloneq (\phi\levelling \land \psi\levelling)
	& 
	(\lnot\phi)\levelling &\coloneq \lnot \phi\levelling \\
	(\exists x \phi)^{s} &\coloneq (\exists x \subseteq s)\phi^{s} & 
	(\poss \phi)\levelling &\coloneq (\exists t : \levpred)\phi^t
	\\
	(\someearlier \phi)^{s} &\coloneq (\exists t : \levpred)(t \in s \land \phi^t) & 
	(\somelater \phi)^{s} &\coloneq (\exists t : \levpred)(s \in t \land \phi^t) 
\end{align*}\noindent
Note that levelling is defined using variables; to illustrate: $(x \in y)\levelling$ is $(x \in y \subseteq s)$, but $(\poss x \in y)\levelling$ is $(\exists t : \levpred)x \in y \subseteq t$. We now have a deep result about modalization and levelling (see \S\ref{s:2:app:ns:deductive}):\footnote{Studd proves similar results. Compare: 
	\eqref{ns:ded:m} with
	\textcites[Theorem 23 p.719]{Studd:ICS}[Proposition 18 p.263]{Studd:EML}; 
	\eqref{ns:ded:ml} with 
	\textcites[Lemma 24 p.719]{Studd:ICS}[Lemma 20 p.263]{Studd:EML}; 
	\eqref{ns:ded:l} with 
	\textcites[Lemma 25 p.719]{Studd:ICS}[Lemma 19 p.263]{Studd:EML}; 	
	\eqref{ns:ded:lm} with \textcites[720]{Studd:ICS}[Proposition 22 p.263]{Studd:EML}.
	
	Clauses \eqref{ns:ded:m}--\eqref{ns:ded:l} do not require temporal-linearity. Clause \eqref{ns:ded:lm} does. To see this, consider a model of \PST with four worlds and accessibility relations exhaustively specified by: $\worldvar{w} < \worldvar{v} < \worldvar{u}$ and $\worldvar{w} < \worldvar{t}$ and $\worldvar{w} < \worldvar{u}$. Where $D(\worldvar{x})$ is $\worldvar{x}$'s first-order domain, let $D(\worldvar{w}) = \{\emptyset\}$; $D(\worldvar{v}) = D(\worldvar{t}) = \powerset\{\emptyset\}$ and $D(\worldvar{u}) = \powerset\powerset\{\emptyset\}$.}
\begin{thm}\label{thm:ns:deductive}
	For any \LTfo-formula $\phi$:
	\begin{listn-0}
		\item\label{ns:ded:m} If $\LTfo \proves \phi$, then $\LPSTfo \proves \phi\modalize$
		\item\label{ns:ded:ml} $\LTfo \proves \phi \liff (\phi\modalize)\levelling$
	\end{listn-0}
	For any \LPSTfo-formula $\phi$:
	\begin{listn}
		\item\label{ns:ded:l} If $\LPSTfo \proves \phi$, then $\LTfo \proves \levpred(s) \lonlyif \phi\levelling$
		\item\label{ns:ded:lm} $\LPSTfo \proves \maxlev(s) \lonlyif (\phi \liff (\phi\levelling)\modalize)$
	\end{listn}
\end{thm}\noindent
This result entails that modalization and levelling are faithful (see Corollary \ref{cor:mutuallyfaithfully}). But Theorem \ref{thm:ns:deductive} is much stronger than a statement of mutual faithful interpretability; it is almost a \emph{definitional equivalence} between \LTfo and \LPSTfo. This claim, though, requires some explanation.\footnote{I know of no existing analogue of definitional equivalence between non-modal and modal theories (such as \LTfo and \LPSTfo); this is my best attempt to provide such an analogue. For a general overview to definitional equivalence in non-modal settings, see e.g.\  \textcite[ch.5]{ButtonWalsh:PMT}.} 

Roughly speaking, to say that two theories are definitionally equivalent is to say that each interprets the other, and that combining the interpretations gets us back exactly where we began. To make this rough idea precise for the case of first-order theories, we say that \textspaced{S} and \textspaced{T} are definitionally equivalent iff there are interpretations $I$ and $J$ such that for any $\textspaced{S}$-formula $\phi$: 
(1) if $\textspaced{S} \proves \phi$ then $\textspaced{T} \proves \phi^I$; and 
(2) $\textspaced{S} \proves \phi \liff (\phi^I)^J$; and for any $\textspaced{T}$-formula $\phi$:
(3) if $\textspaced{T} \proves \phi$ then $\textspaced{S} \proves \phi^J$; and 
(4) $\textspaced{T} \proves \phi \liff (\phi^J)^I$. 
Clauses (1) and (3) tell us we have interpretations; clauses (2) and (4) make precise the idea that `combining the interpretations gets us back exactly where we began'. 

The clauses of Theorem \ref{thm:ns:deductive} are extremely similar to those of a paradigm definitional equivalence. So, Theorem \ref{thm:ns:deductive} is almost a statement of definitional equivalence. Almost; but not quite. We must say something about $s$ in clauses \eqref{ns:ded:l} and \eqref{ns:ded:lm} of Theorem \ref{thm:ns:deductive}, thereby disrupting the similarity. So: we do not have a definitional equivalence; but we \emph{almost} do. 

Since `almost-definitional-equivalence' is quite long-winded, and definitional equivalence is sometimes known as `synonymy', I call this a (deductive) \emph{near-synonymy} between \LTfo and \LPSTfo. 

\section{Near-synonymy: first-order, semantic}\label{s:2:ns:1:sem}
Theorem \ref{thm:ns:deductive} is deductive, but we can extract semantic content from it. (In what follows, my discussion of modal semantics should be understood in terms of \emph{connected} Kripke structures, i.e.\  variable domain Kripke structures where all worlds are path-connected.)\footnote{See footnote  \ref{fn:def:pathconnected} for the definition of \emph{path-connected}.} 

Modalization is defined syntactically, but it has obvious semantic import: as noted, it tells us to replace what \emph{happens} with what \emph{could happen}. This motivates a definition:\footnote{See \textcite[154--5]{Studd:EML}.}
\begin{define}
	Let $\model{P}$ be any connected Kripe structure. Its \emph{flattening}, $\actualize{\model{P}}$, is the following non-modal structure:  
	$\actualize{\model{P}}$'s domain is $\model{P}$'s global domain; and $\actualize{\model{P}} \models a \in b$ iff $\model{P} \mmodels \poss a \in b$.
\end{define}\noindent
Levelling has similar semantic import: it tells us to regard possible worlds as levels. So:
\begin{define}
	Let $\model{A}$ be any non-modal structure. Its \emph{potentialization}, $\potentialize{\model{A}}$, is the following connected Kripke structure: $\potentialize{\model{A}}$'s worlds are those $s$ such that $\model{A}\models \levpred(s)$; accessibility is given by $r < s$ iff $\model{A} \models r \in s$; $\potentialize{\model{A}}$'s global domain is just $\model{A}$'s domain; $\potentialize{\model{A}} \mmodels[s] {a \in b}$ iff $\model{A} \models a \in b \subseteq s$; and $\potentialize{\model{A}} \mmodels[s] a = b$ iff $\model{A} \models a = b \subseteq s$. 
\end{define}\noindent
By considering flattening and potentialization, we can move between models of \LTfo and connected Kripke models of \LPSTfo. To make this movement almost seamless (but only almost; see below), we need one last general construction; intuitively, this construction will allow us to take a Kripke structure, $\model{P}$, and create a new structure, $\model{P}_{f}$, by disrupting the `identities' of $\model{P}$'s worlds (and perhaps duplicating some worlds):
\begin{define}
	Let $\model{P}$ be any connected Kripke structure. Let $f$ be any surjection whose range is $\model{P}$'s set of worlds. Then $\model{P}_{f}$ is the following connected Kripke structure: 
	$\model{P}_{f}$'s set of worlds is $\domain{f}$;
	accessibility is given by $\worldvar{v} < \worldvar{w}$ in $\model{P}_{f}$ iff $f(\worldvar{v}) < f(\worldvar{w})$ in $\model{P}$; 
	$\model{P}_{f}$ has the same global domain as $\model{P}$;
	and $\model{P}_{f} \mmodels[\worldvar{w}] R(\vec{a})$ iff $\model{P}  \mmodels[f(\worldvar{w})] R(\vec{a})$ for all $R$ (including identity).
\end{define}\noindent 
We now have the following result (see \S\ref{s:2:app:semantic}):\footnote{Clause \eqref{ns:sem:ap} requires linearity, since $\actualize{\model{P}}$ has well-ordered levels.}
\begin{thm}\label{thm:ns:semantic}
	\textcolor{white}{.}
	\begin{listn-0}
		\item\label{ns:sem:a} If $\model{P} \mmodels \LPSTfo$, then  $\actualize{\model{P}} \models \LTfo$
		\item\label{ns:sem:ap} If $\model{P} \mmodels \LPSTfo$, then there is a surjection $f$ such that $\model{P} = ({\potactualize{\model{P}}})_{f}$
		\item\label{ns:sem:p} If $\model{A} \models \LTfo$, then $\potentialize{\model{A}} \mmodels \LPSTfo$
		\item\label{ns:sem:pa} If $\model{A} \models \LTfo$, then $\model{A} = {\actpotentialize{\model{A}}}$
	\end{listn-0}
\end{thm}\noindent
This is a semantic reworking of Theorem \ref{thm:ns:deductive}. Consequently, it is \emph{almost} a statement of (semantic) definitional equivalence. Recall that, roughly speaking, two theories are definitionally equivalent iff each interprets the other, and that combining the interpretations gets us back exactly where we began. In \S\ref{s:2:ns:1:ded}, I precisely defined this idea for (non-modal) first-order theories in deductive terms. The same idea can be defined in semantic terms. To say that $\textspaced{S}$ and $\textspaced{T}$ are definitionally equivalent is to say that they (respectively, and uniformly from interpretations) define operations, $g$ and $h$, such that: if $\model{B} \models \textspaced{T}$, then both (1) $g\model{B} \models \textspaced{S}$ and (2) $\model{B} = hg\model{B}$; and if $\model{A} \models \textspaced{S}$, then both (3) $h\model{A} \models \textspaced{T}$ and (4) $\model{A} = gh\model{A}$. Clauses (1) and (3) tell us that we have interpretations; clauses (2) and (4) make precise the idea that `combining the interpretations gets us back exactly where we began'. 

Theorem \ref{thm:ns:semantic} has a very similar shape. So it is almost a (semantic) statement of definitional equivalence between \LTfo and \LPSTfo. Again, though: almost, but not quite. Clause \eqref{ns:sem:ap} of Theorem \ref{thm:ns:semantic} does not tell us that $\model{P} = {\potactualize{\model{P}}}$, as a definitional equivalence would require, but introduces a slight wrinkle. So I will say that we have a semantic \emph{near-synonymy}.

The wrinkle I just mentioned is unavoidable. Fix some $\model{O} \mmodels \LPSTfo$ and $f$ so that $\model{O} \neq \model{O}_{f}$. Clearly $\actualize{\model{O}}  = \actualize{(\model{O}_{f})}$, so that ${\potactualize{\model{O}}} = {\potactualize{(\model{O}_{f})}}$; so we cannot in general have that $\model{P} = {\potactualize{\model{P}}}$. Moreover, this scarcely depends upon the specific definitions of flattening and potentialization; it is an inevitable consequence of the fact that modal semantics has an extra degree of freedom compared with non-modal semantics (the `identities' of worlds, which $f$ can disrupt). 

\section{Near-synonymy: second-order}\label{s:2:ns:2}
I have outlined near-synonymies for the first-order theories \LTfo and \LPSTfo. I now want to consider near-synonymies for the second-order theories. 

In what follows, I assume that \LT's (non-modal) background logic treats second-order identity as co-extensionality, i.e.\ $\forall F \forall G(\forall \vec{x}(F(\vec{x}) \liff G(\vec{x})) \lonlyif F = G)$. Similarly, I assume that all potentialists treat second-order identity as co-intensionality, i.e.:
\begin{listsoaxiom}
	\labitem{Coint}{m:coint} 
	$\forall F \nec \forall G(\nec \forall x_1 \ldots \nec \forall x_n(\poss F(\vec{x}) \liff \poss G(\vec{x})) \lonlyif \poss F = G)$
\end{listsoaxiom}
To take things further, though, I must separately consider two different approaches to second-order entities:  \emph{necessitism} and \emph{contingentism}.\footnote{I use `necessitism' and `contingentism' in roughly \citepossess{Williamson:MLM} sense, though note that the relevant modality here is \emph{potentialist}.}

\subsection{Second-order necessitism}\label{s:2:ns:2:nec}
Second-order necessitism treats second-order entities as necessary existents. We can implement this formally via these axioms:
\begin{listsoaxiom}
	\labitem{Ex$_\necessitistmarker$}{nec:nec} $\existspred{F}$, for any second-order variable  `$F$'
	\labitem{Comp$_\necessitistmarker$}{nec:comp} $\exists F \nec \forall x_1 \ldots \nec \forall x_n(\poss F(\vec{x}) \liff \poss\phi)$, for any formula $\phi$ not containing `$F$'
	\labitem{Inst$_\necessitistmarker$}{m:inst}
	$\forall F \nec \forall x_1 \ldots \nec \forall x_n(\poss F(\vec{x}) \lonlyif \nec ((\existspred{x_1} \land \ldots \land \existspred{x_n}) \lonlyif F(\vec{x})))$
\end{listsoaxiom}
The scheme \ref{nec:nec} guarantees that every second-order entity is a necessary existent. \ref{nec:comp} is a kind of potentialized Comprehension principle. Then \ref{m:inst} guarantees that second-order entities have their instances essentially (cf.\ \ref{m:mem}). 

Let \LPSTnec, for necessitist-\LPST, add these axioms and \ref{m:coint} to \LPST.\footnote{Note that we retain plain vanilla Comprehension; see footnote \ref{fn:comprehensionisvanilla}.} Unsurprisingly, our earlier results are easily extended, to show that \LT and \LPSTnec are deductively and semantically near-synonymous (see Theorems \ref{thm:ns:deductive:nec} and \ref{thm:ns:semantic:nec}).

\subsection{Second-order contingentism}\label{s:2:ns:2:con}
In contrast with necessitism, second-order contingentism holds that a second-order entity exists iff all its (possible) instances do. Contingentists will therefore spurn \ref{nec:nec}, \ref{m:inst}, and \ref{nec:comp}, and instead adopt:
\begin{listsoaxiom}
	\labitem{Ex$_\contingentmarker$}{m:poss} $\poss \existspred{F}$, for any second-order variable `$F$'
	\labitem{Inst$_\contingentmarker$}{con:inst} 
	$\forall F \nec \forall x_1 \ldots \nec \forall x_n(\poss F(\vec{x}) \lonlyif \nec(\existspred{F} \lonlyif F(\vec{x})))$
\end{listsoaxiom}
retaining plain-vanilla Comprehension. Call the result \LPSTcon, for contingentist-\LPST.

Potentialists who treat (monadic) second-order quantification as plural quantification are likely to be contingentists.\footnote{This is \citepossess{Boolos:TB} suggested interpretation of monadic second-order logic. For the link to contingentism, see \textcites[249]{Williamson:MLM}[157--62]{Studd:EML}. The discussion in this paragraph closely follows Studd.} After all, necessitism proves $\exists F \lnot \poss \exists a \nec \forall x(\poss F(a) \liff \poss x \in a )$; read plurally, this contradicts the Maximality Thesis, that any sets can form a set (see \S\ref{s:2:story}). Moreover, the same example establishes that \LT and \LPSTcon are \emph{not} deductively near-synonymous. After all, \LT proves $\exists  F \lnot \exists a \forall x(F(a) \liff x \in a)$, whose modalization will contradict the Maximality Thesis. 

Instead, \LPSTcon is deductively and semantically near-synonymous with a weakened version of \LT. To obtain this weakening, note that contingentists, in effect, restrict second-order entities to the worlds at which their instances occur. Since worlds go proxy for levels, the non-modal equivalent should restrict second-order entities to those which are bounded by levels. Specifically, let $\vec{x}\subseteq s$ abbreviate $(x_1 \subseteq s \land \ldots \land x_n \subseteq s)$, and let $F \sofoundat s$ abbreviate $\forall\vec{x}(F(\vec{x}) \lonlyif \vec{x} \subseteq s)$. Then bounded Level Theory, \LTb, is the theory whose axioms are \ref{ext}, \ref{sep}, \ref{lt:strat}, and:
\begin{listaxiom}
	\labitem{Strat$_\boundedmarker$}{ltb:strat} $\forall F(\exists s : \levpred)F \sofoundat s$
	\labitem{Comp$_\boundedmarker$}{ltb:comp} $(\forall s : \levpred)(\exists F \sofoundat s)(\forall \vec{x} \subseteq s)(F(\vec{x}) \liff \phi)$, for any $\phi$ not containing `$F$'	
\end{listaxiom}
with \ref{ltb:comp} \emph{replacing} the usual Comprehension scheme. Our earlier results can then be extended, to show that \LPSTcon and \LTb are near-synonymous, both deductively and for a Henkin semantics (see Theorems \ref{thm:ns:deductive:con} and \ref{thm:ns:semantic:con}).

So far, deductive and semantic results have gone hand-in-hand. However, they can be prised apart, by considering \emph{full} semantics for second-order logic. For non-modal structures, full (actualist) semantics treats the (monadic) second-order domain as the powerset of the first-order domain. For connected Kripke structures, full contingentist semantics treats a world's (monadic) second-order domain as the powerset of that world's first-order domain. This full semantics is sufficiently rich, that \LPSTcon is not merely near-synonymous with \LTb, but with \LT \emph{itself} (see Theorem \ref{thm:ns:semantic:full}).

\section{The significance of the near-synomies}\label{s:2:equivalencethesis}
The following table summarises the near-synonymies of \S\S\ref{s:2:ns:1:ded}--\ref{s:2:ns:2}: 
\begin{center}
	\setlength\tabcolsep{0.025\textwidth}
	\noindent\begin{tabular}{@{}p{0.05\textwidth}p{0.075\textwidth}p{0.10\textwidth}p{0.15\textwidth}@{}}
		& & deductive & semantic \tabularnewline\addlinespace
		\hline \addlinespace
		\LTfo & \LPSTfo & $\checkmark$& $\checkmark$\\
		\LT & \LPSTnec & $\checkmark$& $\checkmark$\\
		\LTb & \LPSTcon & $\checkmark$ & $\checkmark$\\
		\LT & \LPSTcon & $\times$& full only
	\end{tabular}
\end{center}
To appreciate the significance of these results, consider Paula, a potentialist who uses linear time as an heuristic for her favourite mathematical modality. Paula admires the mathematical work undertaken within $\ZF_1$. However, she regards $\ZF_1$ as lamentably \emph{actualist}, since it lacks modal operators. Fortunately, there is an extension of \LPSTfo---call it $\clearme{LPZF}_1$--- which is near-synonymous with $\ZF_1$.\footnote{Let $\clearme{LPZF}_1  = \LPSTfo \cup \Setabs{\phi^{\posssub}}{\phi \in \ZF_1}$; the near-synomyny holds as $\LTfo \subset \ZF_1$ (see Pt.\ref{pt:lt} \S\ref{s:1:ltsubzf}).} Leaning on this near-synonymy, Paula can regard (worryingly actualist) $\ZF_1$ as a notational-variant of (reassuringly potentialist) $\clearme{LPZF}_1$. Indeed, by modalization and levelling, Paula can move fluidly between $\ZF_1$ and $\clearme{LPZF}_1$. 

The same idea cuts the other way. Actualist Alan may initially be somewhat perplexed by the boxes and diamonds which pepper Paula's work. But Alan need not remain confused for long: modalization and levelling allow him to make perfect sense of Paula, as using a notational-variant of $\ZF_1$. 

\subsection{Outlining an Equivalence Thesis}\label{s:2:argumentforET}
The ease with which Paula and Alan can communicate with each other, despite their philosophical differences, suggests a further thought: 
\begin{storytime}\textbf{The Potentialist/Actualist Equivalence Thesis.}
	Actualism and potentialism do not disagree; they are different but equivalent ways to express the same facts.
\end{storytime}\noindent
Putnam was the foremost proponent of such a Thesis.\footnote{\textcite[8--9]{Putnam:MWF} specifically uses the phrase `the same facts'.} I will say more about Putnam in \S\ref{s:2:putnamperse}; first, I want to assess the Equivalence Thesis directly. Specifically, I want to consider the following, concrete argument for the Equivalence Thesis:
\begin{listl-0}
	\item\label{pet:lt} \LT correctly axiomatizes the idea of an actual hierarchy of sets.
	\item\label{pet:pst} \LPST correctly axiomatizes the idea of a (linear) potential hierarchy of sets.
	\item\label{pet:de} Theories like \LT and \LPST are near-synonymous.
	\item[So:] the Equivalence Thesis obtains.
\end{listl-0}
I am very sympathetic to this argument. However, I am not yet certain of its soundness. In the remainder of this section, I will explain how the argument is best resisted, but also suggest that the Equivalence Thesis remains plausible in the face of such resistance.

The first two premises of the argument are perfectly secure: I established \eqref{pet:lt} in Pt.\ref{pt:lt} \S\S\ref{s:1:st}--\ref{s:1:wodiscussion}, and \eqref{pet:pst} in \S\S\ref{s:2:story}--\ref{s:2:wodiscussion} of this paper. But I should emphasise the caveat in \eqref{pet:pst}. Whilst every potentialist should accept \PST, embracing linearity requires a further step. So: this argument for the Equivalence Thesis can be resisted, straightforwardly, by denying that potentialists can/should assume linearity. 

Premise \eqref{pet:de}, however, contains a sneaky weasel-clause, `theories \emph{like}\ldots'. I will criticise this weaseling in \S\ref{s:2:ET:conpot}. My more pressing concern, though, is whether we could even \emph{hope} to infer the Equivalence Thesis from \eqref{pet:lt}--\eqref{pet:de}.\footnote{\textcite[\S\S5.6, 5.8, 14.7]{ButtonWalsh:PMT} offer some complementary thoughts, about the difficulties of drawing phliosophical conclusions from formal equivalences.}

\subsection{On drawing philosophical conclusions from formal equivalences}\label{s:2:ET:philfromform}
Near-synomy is an extremely tight, formal, equivalence between modal and non-modal theories. Still, theories can be a equivalent in some purely \emph{formal} sense, whilst being non-equivalent in \emph{other} important senses.

To illustrate, suppose Noddy systematically calls red things `green' and green things `red'. Defining interpretations by swapping colour-predicates, Noddy's theory of the empirical world may be definitionally equivalent with my own. Still, if we hold fixed the interpretation of colour-predicates, then we will say that Noddy's theory is simply mistaken; Noddy says `grass is red', but grass is green. 

This noddy example illustrates a simple moral: whether formally equivalent theories `express the same facts' depends upon how firmly we have pinned down the interpretation of the theories' expressions. In the case of Noddy, the relevant expressions colour-predicates. In discussing the Potentialist/Actualist Equivalence Thesis, the relevant expressions are quantifiers and modal operators. And this indicates how discussions of the Equivalence Thesis are likely to play out. 

Suppose you think that we have a firm grasp on the concepts used within the metaphysics of mathematics. In particular, suppose you are convinced that there is a clear difference in meaning between `there is' and `there could be' (as used by potentialists), which does not depend upon their use in any particular \emph{formal} theories. The near-synonymies essentially ask you to move between what `there is' and what `there could be'. Given your prior conviction, you will regard this as a change in subject matter. So you will insist that actualism and potentialism make different claims, and reject the Equivalence Thesis.

Suppose instead, though, that you embrace a rather different attitude. You think that, in advance of any particular formal theorising, it is not entirely clear how one might go about distinguishing between the meanings of `there is' and `there could be' (in mathematical contexts). Indeed, you think that any differences in their meaning would have to be revealed by differences in their use. In that case, you will likely find the argument of \S\ref{s:2:argumentforET} extremely compelling. After all, the near-synonymies establish that there is no significant difference between `$\exists$' in \LTfo and `$\poss\exists$' in \LPSTfo.\footnote{\textcite[588]{Soysal:WUSNS} makes a similar point against any potentialists who treat mathematical possibility as a primitive notion. However, Soysal states that `the potential and [actual] iterative hierarchies are isomorphic, and modal and non-modal set theories are mutually interpretable'. Mutual interpretability is insufficient to support this point (see the \emph{Second} point of \S\ref{s:2:putnamperse}); and it is imprecise to describe potentialist and actualist hierarchies as isomorphic. Soysal's point is better made by appealing to near-synonymy.} 

\subsection{Equivalence and contingentist-potentialism}\label{s:2:ET:conpot}
The case of \LTfo and \LPSTfo is, though, the very simplest case. The situation concerning second-order theories is more complicated, and this merits scrutiny. 
	
Consider Edna, a potentialist who (i) embraces contingentism and (ii) thinks that time is endless, who also (iii) uses second-order logic, whilst (iv) eschewing the full semantics. So Edna embraces an extension of \LPSTcon.\footnote{See \S\ref{s:app:ltb} for details of Edna's theory. By the results of \S\ref{s:2:ns:2} and \S\ref{s:app:ltb}, if Edna drops any of (i)--(iv), then her favourite theory will be near-synonymous (in some salient sense) with \LT itself, rather than \LTb.} As we saw in \S\ref{s:2:ns:2:con}, though, this theory is \emph{not} near-synonymous (whether deductively or using Henkin semantics) with an extension of \LT; we must retreat to \LTb . Edna therefore takes issue with the weasel-clause in premise \eqref{pet:de} of the argument for the Equivalence Thesis. Indeed, she goes further, rebutting the argument as follows: actualists will insist that $\exists F \lnot \exists a \forall x (F(x) \liff x \in a)$; the modalization of this claim is $\poss \exists F \lnot \poss \exists a \nec \forall x(\poss F(x) \liff \poss x \in a)$; this is inconsistent with her favourite potentialist set theory; so potentialism and actualism genuinely \emph{disagree}.\footnote{Thanks to Geoffrey Hellman and \O{}ystein Linnebo for raising concerns along these lines.}
	
This rebuttal of the Equivalence Thesis is exactly as strong as our grasp on the relevant ideology. \emph{If} we have a firm grasp of Edna's intended potentialist modality (independently of the formalism), and how that modality contrasts with actuality, and of the sense of (higher-order) quantification, and why contingentism (but not the use of full second-order semantics) is suitable, \emph{then} Edna's rebuttal will succeed. For, in that case, attempts to move between discussing what `there is' and what `there could be' will amount to a change in truth-value, and therefore also a change in subject matter. But if our grasp of the relevant ideology is insufficiently firm, then Edna's worry will melt away. Edna, then, presents us with an interesting way to resist the Equivalence Thesis, which dovetails with the line of resistance offered in \S\ref{s:2:ET:philfromform}. 

The upshot is that the failure or success of the Equivalence Thesis turns on whether potentialists can supply us with a sufficiently firm grasp of their favoured metaphysical-mathematical-modal concepts. I am genuinely unsure whether they can, but I cheerfully present this as a challenge. 
	
\subsection{Putnam on the equivalence of modal and non-modal theories}\label{s:2:putnamperse}
To conclude my discussion of the Equivalence Thesis, I want to revisit Putnam. As mentioned in \S\ref{s:2:argumentforET}, the Thesis is hugely indebted to Putnam, who claimed in \cite*{Putnam:MWF} that modal and non-modal theories are `equivalent'. However, it is worth emphasizing a few of the differences between Putnam's \cite*{Putnam:MWF} claim and my Equivalence Thesis.
	
\emph{First.} Putnam did not say much about the modality he had in mind, except to connect `$\poss$' with possible `standard concrete models for Zermelo set theory'.\footnote{\textcite[20--1]{Putnam:MWF}.} My discussion is restricted to a potentialist modality, though I have deliberately left room for various different versions of potentialism.\footnote{\textcite[262--6]{Linnebo:PMML} offers good reasons to suggest that Putnam \emph{should} have considered a potentialist modality.}
	
\emph{Second.} Putnam did not precisely define the formal notion of `equivalence' he had in mind. He sometimes considers the \emph{mutual interpretability} of modal and non-modal theories;\footnote{E.g.\ \textcite[8]{Putnam:MWF} `the primitive terms of each admit of definition by means of the primitive terms of the other theory, and then each theory is a deductive consequence of the other.'} but mutual interpretability is far too weak to sustain anything like the Equivalence Thesis.\footnote{\textcite[260--2]{Linnebo:PMML} makes this point. To bring it out in another way, note that $\textspaced{PA}$ and $\textspaced{PA} + \lnot\textspaced{Con}(\textspaced{PA})$ are mutually interpretable, but are surely \emph{not} equivalent ways to express the same facts.} By contrast, my formal notion of `equivalence' is \emph{near-synonymy}. 
	
\emph{Third.} Putnam ultimately retracted his version of the Equivalence Thesis.\footnote{\textcite[11.Dec.2014]{Putnam:blog}.} He claimed that mathematics is `about proofs, ways of conceiving of mathematical problems, mathematical approaches, and much more', and worried that his interpretation would not preserve such things. Now, these considerations might tell against Putnam's \cite*{Putnam:MWF} claim; but they only highlight the plausibility of my Equivalence Thesis. My near-synonymies simply formalize the intuitive and obvious point that \LT's levels simulate \LPST's possible worlds, and vice versa (see \S\ref{s:2:LPST}); this simulation straightforwardly preserves proofs; and this is precisely why it is so plausible that \LT and \LPST do not really differ over `ways of conceiving mathematical problems, mathematical approaches', or anything else that matters.
	
\emph{Fourth.} Having decided that modal and non-modal formulations of set theory genuinely disagree, Putnam came to favour the former, on the grounds that non-modal set theories face `a generalization of a problem first pointed out by Paul Benacerraf\ldots e.g.\ are sets a kind of function or are functions a sort of set?'\footnote{\textcite[13.Dec.2014]{Putnam:blog}.} Again, this might detract from Putnam's \cite*{Putnam:MWF} claim, but it has no force against my Equivalence Thesis.  If \LT and \LPST are equally good in all other regards---as I think they might be---then choosing potentialism (with its distinctive modality) \emph{over} actualism (with its distinctive ontology) is exactly as arbitrary as saying that functions are a kind of set (rather than vice versa). 

\section{Conclusion, and predecessors}\label{s:2:alternates}
The Tensed Story articulates the bare idea of a potential hierarchy of sets. \PST axiomatizes that bare idea. Whilst it takes no stance on the height of any potential hierarchy, it ensures persistence and well-foundedness. Moreover, versions of \PST are near-synonymous with versions of the non-modal theory \LT. And these near-synonymies both sharpen and leave plausible the idea that there is no deep difference between actualism and potentialism. 

I will close this paper by comparing \PST with some alternative potentialist set theories. 

\subsection{Parsons and Linnebo}\label{s:parsonslinnebo}
In formulating their modal set theories, Parsons and Linnebo do not use a temporal logic.\footnote{\textcites{Parsons:ICS}{Parsons:SM}{Linnebo:PHS}[ch.12]{Linnebo:TO}.} Instead, they use a single modal operator, $\poss$, whose background logic is S4.2, and which can be glossed as `now and henceforth'.

The asymmetry of this operator generates a deep expressive problem.\footnote{For related problems, see \textcites[723--4]{Studd:ICS}[169--71]{Studd:EML}.} Stated non-modally: there is a stage (the initial stage) at which nothing has any members. Potentialists should therefore want to be able to prove: \emph{possibly, nothing has any members}, i.e.\ $\poss \forall x \forall y\, \ x \notin y$. But this cannot be a theorem for Parsons or Linnebo. To see why, suppose otherwise; then $\nec \poss \forall x \forall y\, \ x \notin y$ is also a theorem, by Necessitation; but this is catastrophic, for it catastrophically entails that there is always a later moment at which nothing has any members. 

This problem does not arise in \PST. There, $\poss$ obeys S5, and \PST proves $\poss \forall x \forall y\, x \notin y$.

\subsection{Studd}\label{s:2:studd}
In using a tensed logic to formulate \PST, I am entirely indebted to Studd. Moreover, Studd proves a result like Theorem \ref{thm:ns:deductive} for his modal set theory, MST. So my \PST is similar to Studd's MST, and owes a great deal to it. However, it is worth noting two differences.

The minor difference concerns our versions of \ref{pst:pri}. Studd's MST has: all of $a$'s members are found \emph{together} before $a$ is found.\footnote{\textcites*[712]{Studd:ICS}[164--5]{Studd:EML}.} My \PST has: \emph{each} of $a$'s members is found before $a$ is found. The slight difference emerges only at limit worlds:\footnote{Where $\worldvar{w}$ is a limit world iff $(\forall \worldvar{u} < \worldvar{w})\exists \worldvar{v}(\worldvar{u} < \worldvar{v} < \worldvar{w})$.} in Studd's MST, $a$ exists at a limit world iff $a$ existed earlier; in my \PST, $a$ exists at a limit world iff all of $a$'s members existed earlier.

The major difference concerns the richness of Studd's modal schemes. Studd's MST explicitly adopts modal axioms which guarantee \emph{linearity}, \emph{persistence}, \emph{well-ordering}, and that time is \emph{endless}.\footnote{\textcites[704]{Studd:ICS}[152, 252]{Studd:EML} guarantees persistence via Barcan-formulas; see also \textcites[210]{Linnebo:PHS}[207]{Linnebo:TO}. \textcites[702--4]{Studd:ICS}[152, 251--2]{Studd:EML} guarantees well-ordering via a L\"ob-scheme; \textcites[296]{Parsons:ICS}[318]{Parsons:SM}[216]{Linnebo:PHS}[206]{Linnebo:TO} guarantee well-ordering via non-modal means.} My \PST only assumes past-directedness, and instead \emph{proves} persistence and well-foundedness (see \S\ref{s:2:wodiscussion}). Proof has three virtues over explicit assumption. First: my \PST is considerably leaner than Studd's MST. Second: it will be strictly easier for potentialists to try to explain why they are entitled to assume past-directedness, than to try to justify Studd's richer assumptions.\footnote{\label{fn:Studdtroubles}To illustrate: \textcite[144--53]{Studd:ICS} glosses $\allearliersub$ as `however the lexicon is interpreted by preceding interpretations' and $\alllatersub$ as `however the lexicon is interpreted by succeeding interpretations'. I worry that Studd does not manage to show that, so glossed, these operators should obey the schemes for linearity, persistence, or well-ordering. However, past-directedness might well be justifiable; and from there we can prove persistence and well-foundedness, via Theorem \ref{thm:pst:key}.} Third: as in \S\ref{s:2:wodiscussion}, the proofs of persistence and well-foundedness show `how little choice there is in setting up' a potential hierarchy.

\startappendix
\section{Elementary results concerning \PST}\label{s:2:app:pst}
The time has come to prove the results stated in Part \ref{pt:pst}. I will start with some elementary results within \PST, building up to Theorem \ref{thm:pst:key} of \S\ref{s:2:wodiscussion}. My proofs are semantic, relying on standard soundness and completeness results for (connected) Kripke frames. I use bold letters, $\worldvar{w}, \worldvar{v}, \worldvar{u}, \ldots$, for arbitrary worlds (note that this differs from my use of bold letters in Pt.\ref{pt:lt} and Pt.\ref{pt:blt}).

In what follows, we must not assume that expressions like `$\Setabs{x}{\phi(x)}$' are rigid designators; we should read `$a = \Setabs{x}{\phi(x)}$' as abbreviating `$\forall x(x \in a \liff \phi(x))$', which may be true in one world and false in another. Similarly, recall that `$a = \pot{b}$' abbreviates `$\forall x(x \in a \liff \exists c(x \subseteq c \in b))$'. 

I start with two very elementary results:
\begin{lem}[\PST]\label{lem:modalcore:ext}
	\ref{ext} holds.
\end{lem}
\begin{proof}
	Fix $a$ and $b$ at $\worldvar{w}$, and assume  $\mmodels[\worldvar{w}] \forall x(x \in a \liff x \in b)$. Fix $x$ at world $\worldvar{u}$, now $\mmodels[\worldvar{u}] \poss x \in a$ iff $\mmodels[\worldvar{w}] x \in a$ (by \ref{m:mem}) iff $\mmodels[\worldvar{w}] x \in b$ iff $\mmodels[\worldvar{u}] \poss x \in b$; so $\mmodels[\worldvar{u}] \poss a = b$ by \ref{m:ext}. Hence $\mmodels[\worldvar{w}] a = b$. 
\end{proof}
\begin{lem}[\PST]\label{lem:pst:sep} 
	\ref{sep} holds.
\end{lem}
\begin{proof}
	Using Comprehension, let $G$ be given by $\forall x(G(x) \liff (F(x) \land x \in a))$. If $G(x)$, then $\someearlier \existspred{x}$ by \ref{pst:pri}; so some $b = \Setabs{x}{G(x)} = \Setabs{x \in a}{F(x)}$ exists by \ref{pst:spec}.
\end{proof}\noindent 
Since \PST proves \ref{ext} and \ref{sep}, it proves the key results of Pt.\ref{pt:lt} \S\ref{s:1:LTwo}, concerning the well-ordering of \emph{levels}, in the sense of Pt.\ref{pt:lt} Definition \ref{def:level}. This next result establishes that all of the key notions of that Definition are (weakly) rigid:
\begin{lem}[\PST]\label{lem:pst:robust}\textcolor{white}{.}
	\begin{listn-0}
		\item\label{subset:robust} 	$\forall a(\forall b \subseteq  a)\nec(\existspred{a} \lonlyif (\existspred{b} \land b \subseteq a))$
		\item\label{pot:exists} $\forall a \exists b\ \pot{a} = b$
		\item\label{pot:robust} $\forall a (\forall b = \pot{a}) \nec (\existspred{a} \lonlyif (\existspred{b} \land b = \pot{a}))$
		\item\label{histpred:robust} $(\forall h : \histpred)\nec(\existspred{h} \lonlyif \histpred(h))$
		\item\label{levpred:robust} $(\forall s : \levpred)\nec(\existspred{s} \lonlyif \levpred(s))$
	\end{listn-0}
\end{lem}
\begin{proof}
	\emphref{subset:robust} 	
	Fix $a$ and $b$ at $\worldvar{w}$ such that $\mmodels[\worldvar{w}] b \subseteq a$. Let $a$ exist at $\worldvar{v}$; by \ref{sep} at $\worldvar{v}$ there is $c$ at $\worldvar{v}$ such that $\mmodels[\worldvar{v}] c = \Setabs{x \in a}{\poss x \in b}$; I claim that $\mmodels[\worldvar{v}] c = b$. Fix $x$ at $\worldvar{u}$: if $\mmodels[\worldvar{u}] \poss x \in c$, then $\mmodels[\worldvar{v}] x \in c$ by \ref{m:mem}, so $\mmodels[\worldvar{v}] \poss x \in b$ and $\mmodels[\worldvar{u}] \poss x \in b$; if $\mmodels[\worldvar{u}] \poss x \in b$, then $\mmodels[\worldvar{w}] x \in b \subseteq a$ by \ref{m:mem}, so that $\mmodels[\worldvar{v}] x \in a$  and $\mmodels[\worldvar{v}] \poss x \in b$, i.e.\ $\mmodels[\worldvar{v}] x \in c$, so that $\mmodels[\worldvar{u}] \poss x \in c$. Hence $\mmodels[\worldvar{v}] c = b$ by \ref{m:ext}.
	
	\emphref{pot:exists} Fix $a$. If $\exists z(x \subseteq z \in a)$, then $\someearlier \existspred{x}$ by \ref{pst:pri} and \eqref{subset:robust}. So using \ref{pst:spec} we have some $b$ such that $b = \pot{a} = \Setabs{x}{\exists z(x \subseteq z \in a)}$. 
	
	\emphref{pot:robust}
	Fix $a$ and $b$ at $\worldvar{w}$ with $\mmodels[\worldvar{w}] b = \pot{a}$. Let $a$ exist at $\worldvar{v}$, and using \eqref{pot:exists} fix $c$ such that $\mmodels[\worldvar{v}] c = \pot{a}$. Now $\mmodels[\worldvar{v}] b = c$, by \ref{m:ext} and \eqref{subset:robust}.
	
	\emphreffrom{histpred:robust}\emphref{levpred:robust}
	By \eqref{subset:robust} and \eqref{pot:robust}. 
\end{proof}\noindent
We can now show that levels persist, and also that every world has a maximal level:
\begin{lem}[\PST]\label{lem:pst:levelspersist}	$(\forall s : \levpred)\alllater(\existspred{s} \land \levpred(s))$
\end{lem}
\begin{proof}
	Let $s$ be a level in $\worldvar{w}$. For induction on levels (i.e.\ Pt.\ref{pt:lt} Theorem \ref{thm:es:wo}), suppose that $\mmodels[\worldvar{w}] (\forall r : \levpred)(r \in s \lonlyif \alllater(\existspred{r} \land \levpred(r)))$. Fix $\worldvar{v} > \worldvar{w}$; using \ref{pst:spec} fix $t$ such that $\mmodels[\worldvar{v}] t = 	\pot{\Setabs{x}{(\exists r : \levpred)(x \subseteq r \land \poss r \in s)}}$. I claim that $\mmodels[\worldvar{v}] s = t$; the result will then follow by induction on levels in $\worldvar{w}$ and Lemma \ref{lem:pst:robust}.\ref{levpred:robust}. 	
	
	If $\mmodels[\worldvar{u}] \poss x \in s$, then $\mmodels[\worldvar{w}] x \in s$; so by Pt.\ref{pt:lt} Lemma \ref{lem:es:acc} there is some $r$ such that $\mmodels[\worldvar{w}] x \subseteq r \in s \land \levpred(r)$ ; now $\mmodels[\worldvar{v}] \existspred{r} \land \levpred(r)$ by the induction hypothesis and Lemma \ref{lem:pst:robust}.\ref{levpred:robust}, and $\mmodels[\worldvar{v}] x \subseteq r$ by Lemma \ref{lem:pst:robust}.\ref{subset:robust}; so $\mmodels[\worldvar{v}] x \in t$ and hence $\mmodels[\worldvar{u}] \poss x \in t$. The converse is similar. So $\mmodels[\worldvar{u}] \poss s = t$, and $\mmodels[\worldvar{v}] s = t$ by \ref{m:ext}.
\end{proof}
\begin{lem}[\PST]\label{lem:pst:maxlevel}
	$(\exists s : \levpred)(\forall r : \levpred)(r \subseteq s \land (r \neq s \liff \someearlier\existspred{r}))$
\end{lem}
\begin{proof}
	Using \ref{pst:spec}, let $h = \Setabs{r}{\levpred(r) \land \someearlier\existspred{r}}$. I claim that $h$ is a history. Fix $r \in h$. Clearly $\pot{(r\cap h)} \subseteq \pot{r} = r$ as levels are potent. Conversely, if $a \in r$ then there is some level $q$ such that $a \subseteq  q \in r$ by Pt.\ref{pt:lt} Lemma \ref{lem:es:acc}, and since $\someearlier\existspred{r}$ we have $\someearlier\existspred{q}$; so $q \in r \cap h$ and hence $a \in \pot{(r\cap h)}$. Generalising, $r \subseteq \pot{(r\cap h)}$. So $h$ is a history. Using Lemma \ref{lem:pst:robust}.\ref{pot:exists}, let $s = \pot{h}$. By construction, $s$ is a level. I claim that $s$ has the required properties.
	
	For reductio, suppose that $\someearlier\existspred{s}$; then $s \in h \subseteq s$, contradicting the well-ordering of levels; so $\lnot\someearlier\existspred{s}$. 
	
	Suppose $r \neq s$. Then either $r \in s$ or $s \in r$ by the well-ordering of levels; but if $s \in r$ then $\someearlier\existspred{s}$ by \ref{pst:pri}, a reductio. So $r \in s$. Hence $\someearlier\existspred{r}$ by \ref{pst:pri}, and also $r \subseteq s$ as $s$ is transitive. 
\end{proof}\noindent
From here, we can prove a L\"ob-like scheme for \PST:
\begin{lem}[\PST]\label{cor:pst:lob}
	$\nec(\allearlier\phi \lonlyif \phi) \lonlyif \nec \phi$, for all $\phi$
\end{lem}
\begin{proof}
	For reductio, suppose this is false at $\worldvar{w}$, i.e.\ $\mmodels[\worldvar{w}] \nec(\allearlier\phi \lonlyif \phi)$ but $\mmodels[\worldvar{w}] \poss \lnot \phi$. So $\mmodels[\worldvar{v}] \lnot \phi$ for some $\worldvar{v}$. Since $\mmodels[\worldvar{v}] \allearlier\phi \lonlyif \phi$, there is $\worldvar{u} < \worldvar{v}$ such that $\mmodels[\worldvar{u}] \lnot \phi$. For brevity, let:
		$$\psi(x) \text{ abbreviate }(\lnot \phi \land \levpred(x) \land \lnot \someearlier \existspred{x} \land (\forall q : \levpred)q \subseteq x)$$	
	Now $\mmodels[\worldvar{v}] (\exists s : \levpred)\someearlier\psi(s)$, by Lemmas \ref{lem:pst:levelspersist}--\ref{lem:pst:maxlevel}. Using induction on levels, let $s$ be the $\in$-minimal level in $\worldvar{v}$ such that $\mmodels[\worldvar{v}] \someearlier \psi(s)$. So there is $\worldvar{t} < \worldvar{v}$ with $\mmodels[\worldvar{t}] \psi(s)$. Since $\mmodels[\worldvar{t}] \lnot \phi$ and $\mmodels[\worldvar{t}] \allearlier\phi \lonlyif \phi$ by assumption, there is $\worldvar{t}_0 < \worldvar{t}$ with $\mmodels[\worldvar{t}_0] \lnot \phi$. Using Lemma \ref{lem:pst:maxlevel}, fix $r$ such that $\mmodels[\worldvar{t}_0] \psi(r)$. Now $\mmodels[\worldvar{t}] \levpred(r) \land \someearlier\existspred{r}$ by Lemma \ref{lem:pst:levelspersist}, so $\mmodels[\worldvar{t}] r \in s$ by Lemma \ref{lem:pst:maxlevel} and choice of $s$. So $\mmodels[\worldvar{v}] \levpred(r) \land r \in s \land \someearlier \psi(r)$ by Lemma \ref{lem:pst:levelspersist}, contradicting the choice of $s$.
\end{proof}\noindent
This effectively licenses schematic induction on worlds, enabling us to prove the main result of \S\ref{s:2:wodiscussion}:
\begin{namedthm}[\textbf{\ref{thm:pst:key}} \textnormal{(\PST)}]
	Where $\maxlev(s)$ abbreviates $(\existspred{s} \land \forall x\ x  \subseteq s)$:
	\begin{listn-0}
		\item[\eqref{key:lt}] $\LT$ holds
		\item[\eqref{key:persistence}] $\forall x \alllater\existspred{x}$
		\item[\eqref{key:maxlev}] $(\exists s : \levpred)\maxlev(s)$
		\item[\eqref{key:prune}] $(\forall s : \levpred)\poss\maxlev(s)$
	\end{listn-0}
\end{namedthm}
\begin{proof}
	\emphref{key:lt} It suffices to prove \ref{lt:strat}, i.e.\ that $\forall a (\exists s : \levpred)a \subseteq s$. Fix $\worldvar{w}$, and suppose for induction on worlds that $\mmodels[\worldvar{v}] \forall a(\exists s : \levpred)a \subseteq s$ for all $\worldvar{v} < \worldvar{w}$. Using Lemma \ref{lem:pst:maxlevel}, fix $s$ such that $\mmodels[\worldvar{w}] \levpred(s) \land  \lnot\someearlier\existspred{s} \land (\forall r : \levpred)r \subseteq s$. Suppose $\mmodels[\worldvar{w}] x \in a$; by \ref{pst:pri} there is some $\worldvar{u} < \worldvar{w}$ such that $\mmodels[\worldvar{u}] \existspred{x}$; by assumption there is $r$ such that $\mmodels[\worldvar{u}] \levpred(r) \land x \subseteq r$; now $\mmodels[\worldvar{w}] x \subseteq r \in s$ by Lemmas \ref{lem:pst:robust}--\ref{lem:pst:levelspersist}, so that $x \in s$ as $s$ is potent. Hence $\mmodels[\worldvar{w}] a \subseteq s$. The result follows by Lemma \ref{cor:pst:lob}.
	
	\emphreffrom{key:persistence}\emphref{key:maxlev} Combine \ref{lt:strat} with Lemmas \ref{lem:pst:robust}--\ref{cor:pst:lob}. 
	
	\emphref{key:prune} Fix $\worldvar{w}$, and suppose for induction on worlds that $\mmodels[\worldvar{v}] (\forall s : \levpred)\poss\maxlev(s)$ for all $\worldvar{v} < \worldvar{w}$. Let $s$ be such that $\mmodels[\worldvar{w}] \levpred(s)$. If $\mmodels[\worldvar{w}] \someearlier \existspred{s}$ then $\mmodels[\worldvar{w}] \poss \maxlev(s)$ by our supposition and Lemma \ref{lem:pst:robust}. Otherwise, $\mmodels[\worldvar{w}] (\forall r : \levpred)r \subseteq s$ by the well-ordering of levels and Lemma \ref{lem:pst:maxlevel}, so that $\mmodels[\worldvar{w}] \maxlev(s)$ by \ref{lt:strat}. The result follows by Lemma \ref{cor:pst:lob}. 
\end{proof}\noindent
To round things off, note that \LT's key notions are robust under modalization:
\begin{lem}[\PST]\label{lem:pst:inext}
	\textcolor{white}{.}
	\begin{listn-0}
		\item\label{inext:scheme} $\phi\modalize(\vec{x})$ iff $\nec\phi\modalize(\vec{x})$, for any \LT-formula $\phi(\vec{x})$
		\item\label{inext:sub1} if $\existspred{b} \land b \subseteq a$, then $(b \subseteq a)\modalize$
		\item\label{inext:sub2} if $(b \subseteq a)\modalize$, then $\nec(\existspred{a} \lonlyif b \subseteq a)$
		\item\label{inext:pot1} if $\existspred{b}$ and $(b = \pot{a})\modalize$, then $\existspred{a}$ and $b = \pot{a}$
		\item\label{inext:pot2} if $\existspred{b}$ and $b = \pot{a}$, then $\existspred{a}$ and $(b = \pot{a})\modalize$
		\item\label{inext:hist} if $\existspred{h}$, then $\histpred(h) \liff \histpred\modalize(h)$
		\item\label{inext:lev} if $\existspred{s}$, then $\levpred(s) \liff \levpred\modalize(s)$
	\end{listn-0}
\end{lem}
\begin{proof}
	\emphref{inext:scheme} A routine induction on complexity, using the fact that $\poss$ obeys S5.
	
	\emphreffrom{inext:sub1}\emphref{inext:sub2} Straightforward.
	
	\emphref{inext:pot1} Suppose $\mmodels[\worldvar{w}] \existspred{b}$ and $\mmodels[\worldvar{w}] (b = \pot{a})\modalize$, i.e.\ $\mmodels[\worldvar{w}] \nec \forall x(\poss x \in b \liff (\exists z(x \subseteq z \in a))\modalize)$. 
	
	I first show that $\mmodels[\worldvar{w}] \existspred{a}$. By \ref{sep} there is $c$ at $\worldvar{w}$ such that $\mmodels[\worldvar{w}] c = \Setabs{x \in b}{\poss x \in a}$; I claim $a = c$ using \ref{m:ext}. 
	Fix $x$ at $\worldvar{u}$. If $\mmodels[\worldvar{u}] \poss x \in c$ then clearly $\mmodels[\worldvar{u}] \poss x \in a$. Conversely, if $\mmodels[\worldvar{u}] \poss x \in a$, then letting $x= z$ we have $\mmodels[\worldvar{u}] (\exists z(x \subseteq z \in a))\modalize$ by \eqref{inext:sub1}, hence $\mmodels[\worldvar{w}] \poss x \in b$ so that $\mmodels[\worldvar{w}] x \in b$ and hence $\mmodels[\worldvar{w}] x \in c$ i.e.\ $\mmodels[\worldvar{u}] \poss x \in c$. 
	
	I now show that $\mmodels[\worldvar{w}] b = \pot{a}$. If $\mmodels[\worldvar{w}] x \in b$, then $\mmodels[\worldvar{w}] (\exists z(x \subseteq z \in a))\modalize$, i.e.\ there is $\worldvar{u}$ and $z$ such that $\mmodels[\worldvar{u}] (x \subseteq z \in a)\modalize$; now $\mmodels[\worldvar{w}] x \subseteq z \in a$ by \eqref{inext:sub2} and as $\mmodels[\worldvar{w}] \existspred{a}$. Conversely, if $\mmodels[\worldvar{w}] x \subseteq z \in a$ for some $z$, then $\mmodels[\worldvar{w}] (\exists z(x \subseteq z \in a))\modalize$ by \eqref{inext:sub1}, so $\mmodels[\worldvar{w}] \poss x \in b$ and so $\mmodels[\worldvar{w}] x \in b$. 
	
	\emphref{inext:pot2} Similar to \eqref{inext:pot1}.
	
	\emphreffrom{inext:hist}\emphref{inext:lev} By \eqref{inext:scheme} and \eqref{inext:pot1}--\eqref{inext:pot2}. 
\end{proof}\noindent
All the results of this appendix can be first-orderized straightforwardly. Keen readers will also notice that the proofs of this appendix have made no apparent use of the assumption of past-directedness. Indeed: the only role for past-directedness is to supply us with a possibility operator, $\poss$, which is unrestricted and obeys S5. 
	
\section{Results concerning \LPST}\label{s:2:app:lpst}
I will now turn from \PST to \LPST. As mentioned in \S\ref{s:2:LPST}, linearity allows us to define away $\someearlier$ and $\somelater$ via the map $\phi \mapsto \phi\mltint$. To guarantee that this is so, we use the results of \S\ref{s:2:app:pst} to prove Lemma \ref{prop:pstl:mltint} by a simple induction on complexity; I leave this to the reader.
	
Evidently, $\LPST\mltint$ is a unimodal S5 theory. However, it may be worth noting that it can be given a simpler presentation. Let \MLT be a unimodal S5 theory whose set-theoretic axioms are \ref{m:mem}, \ref{m:ext}, \ref{sep}, and clauses \eqref{key:maxlev}--\eqref{key:prune} of Theorem \ref{thm:pst:key}. The proofs of Lemmas \ref{lem:modalcore:ext}--\ref{lem:pst:robust} go through in \MLT with only tiny adjustments; and it is easy to show that $\MLT \proves \poss\phi\mltint \liff (\someearlier\phi \lor \phi \lor \somelater\phi)\mltint$ for each \LPST-formula $\phi$. It follows that $\LPST\mltint \dashv\vdash \MLT$. By Lemma \ref{prop:pstl:mltint}, then, \LPST and \MLT are (strictly) definitionally equivalent.
	
\subsection{Deductive near-synonymy}\label{s:2:app:ns:deductive}
The key results concerning \LPST, though, are the near-synonymies. I will start with the first-order deductive near-synonymy:
\begin{namedthm}[\ref{thm:ns:deductive}]
	For any \LTfo-formula $\phi$ not containing $s$:
	\begin{listn-0}
		\item[\eqref{ns:ded:m}] If $\LTfo \proves \phi$, then $\LPSTfo \proves \phi\modalize$
		\item[\eqref{ns:ded:ml}] $\LTfo \proves \phi \liff (\phi\modalize)\levelling$
	\end{listn-0}
	For any \LPSTfo-formula $\phi$ not containing $s$:
	\begin{listn}
		\item[\eqref{ns:ded:l}] If $\LPSTfo \proves \phi$, then $\LTfo \proves \levpred(s) \lonlyif \phi\levelling$
		\item[\eqref{ns:ded:lm}] $\LPSTfo \proves \maxlev(s) \lonlyif (\phi \liff (\phi\levelling)\modalize)$
	\end{listn}
\end{namedthm}
\begin{proof}
	\emphref{ns:ded:m} 
	\ref{ext}$\modalize$ is \ref{m:ext}. For \ref{lt:strat}$\modalize$, use Theorem \ref{thm:pst:key}.\ref{key:maxlev}  and Lemma \ref{lem:pst:inext}. For the \ref{sep}$\modalize_1$-instances, fix suitable $\phi$; fix $a$ at $\worldvar{w}$; by \ref{sep} we have some $b$ in $\worldvar{w}$ such that $\mmodels[\worldvar{w}] b = \Setabs{x \in a}{\phi\modalize}$. Fix $x$ at $\worldvar{u}$; now $\mmodels[\worldvar{u}] \poss x \in b$ iff $\mmodels[\worldvar{w}] x \in b$ iff $\mmodels[\worldvar{w}] \phi\modalize \land x \in a$ iff $\mmodels[\worldvar{u}] \phi\modalize \land \poss x \in a$ by Lemma \ref{lem:pst:inext}.
	
	\emphref{ns:ded:ml} 
	A routine induction on complexity.
	
	\emphref{ns:ded:l} 
	The well-ordering and potency of levels yields the levelling of each underlying logical principle. It is then straightforward to obtain the levelling of each \LPSTfo axiom is then straightforward.
	
	\emphref{ns:ded:lm}
	An induction on complexity. The cases of atomic formulas, conjunctions and quantifiers are easy, relying on \ref{m:mem} and Lemma \ref{lem:pst:inext}.\ref{inext:sub1}--\ref{inext:sub2}. 
	
	For quantifiers: using the induction hypothesis, \LPST proves that, if $\maxlev(s)$ then: 
	$(\exists x \phi)$ 
	iff $(\exists x \subseteq s)(\phi^{s})$
	iff $\poss(\exists x \subseteq {s})(\phi^{s})\modalize$
	iff $((\exists x\phi)^{s})\modalize$. 
	
	For modal operators, I will prove the case for $\someearlier$ (the others are similar). Fix $\worldvar{w}$ and, using Theorem \ref{thm:pst:key}.\ref{key:maxlev}, let $\mmodels[\worldvar{w}] \maxlev(s)$; I claim that $\mmodels[\worldvar{w}] \someearlier \phi \liff ((\someearlier\phi)\levelling)\modalize$. 
	
	Suppose $\mmodels[\worldvar{w}] \someearlier \phi$, i.e.\ there is $\worldvar{v} < \worldvar{w}$ such that $\mmodels[\worldvar{v}] \phi$. Using Theorem \ref{thm:pst:key}.\ref{key:maxlev}, let $\mmodels[\worldvar{v}] \maxlev(r)$. By the induction hypothesis, $\mmodels[\worldvar{v}] \phi \liff (\phi^r)\modalize$; so $\mmodels[\worldvar{w}] (\phi^r)\modalize$. Hence $\mmodels[\worldvar{w}] \levpred(r) \land \poss r \in s \land (\phi^r)\modalize$; now by Lemma \ref{lem:pst:inext} we have $\mmodels[\worldvar{w}] \poss \exists r(\levpred\modalize(r) \land \poss  r \in s \land (\phi^r)\modalize)$, i.e.\ $\mmodels[\worldvar{w}] ((\someearlier \phi)^r)\modalize$
	
	Suppose $\mmodels[\worldvar{w}] ((\someearlier \phi)^r)\modalize$, i.e.\ for some $\worldvar{v}$ and some $r$ at $\worldvar{v}$ we have $\mmodels[\worldvar{v}] \levpred\modalize(r) \land \poss  r \in s \land (\phi^r)\modalize$. Using Theorem \ref{thm:pst:key}.\ref{key:prune} and Lemma \ref{lem:pst:inext}, fix $\worldvar{u}$ such that $\mmodels[\worldvar{u}] \maxlev(r)$; note that $\mmodels[\worldvar{u}] (\phi^r)
	\modalize$, so that $\mmodels[\worldvar{u}] \phi$ by the induction hypothesis. Moreover, $\worldvar{u} < \worldvar{w}$, as $\mmodels[\worldvar{v}] \poss r \in s$ and we have assumed linearity. So $\mmodels[\worldvar{w}] \someearlier \phi$.
\end{proof}\noindent
Theorem \ref{thm:ns:deductive} straightforwardly entails that modalization and levelling are \emph{faithful}:
\begin{cor}\label{cor:mutuallyfaithfully}
	\textcolor{white}{.}
	\begin{listn-0}
		\item\label{syntactic:possintfaithful} $\LPSTfo \proves \phi\modalize$ iff $\LTfo  \proves \phi$, for any \LTfo-formula $\phi$ 
		\item\label{syntactic:sintfaithful} $\LTfo  \proves \levpred(s) \lonlyif \phi\levelling$ iff $\LPSTfo \proves \phi$, for any \LPSTfo-formula $\phi$ not containing $s$
	\end{listn-0}
\end{cor}\noindent
I leave the proof to the reader. The reader can also prove these two second-order versions of Theorem \ref{thm:ns:deductive}, mentioned in \S\ref{s:2:ns:2}:
\begin{thm}\label{thm:ns:deductive:nec}
	Theorem \ref{thm:ns:deductive} holds for \LT and \LPSTnec, where we enrich modalization and levelling with these clauses:
	\begin{align*}
		\alpha \modalize &\coloneq \poss\alpha\text{, for atomic }\alpha &		(\exists F \phi)\modalize &\coloneq \poss \exists F \phi\modalize\\
		(F = G)\levelling &\coloneq F = G& 
		F(\vec{x})\levelling &\coloneq (F(\vec{x}) \land \vec{x} \subseteq s)&
		(\exists F \phi)\levelling &\coloneq \exists F \phi\levelling
	\end{align*}
\end{thm}
\begin{thm}\label{thm:ns:deductive:con}
	Theorem \ref{thm:ns:deductive} holds for \LTb and \LPSTcon, where we enrich modalization as above, but instead enrich levelling as follows:
	\begin{align*}
		(F = G)\levelling &\coloneq (F = G \sofoundat s) &
		F(\vec{x})\levelling &\coloneq (F(\vec{x}) \land F \sofoundat s) &
		(\exists F \phi)\levelling& \coloneq (\exists F \sofoundat s)\phi\levelling
	\end{align*}
\end{thm}

\subsection{Semantic near-synonymy}\label{s:2:app:semantic}
I now consider the semantic near-synonymies. The first-order result follows from two lemmas, which are proved by a routine induction on complexity:
\begin{lem}\label{lem:modalize-actualize}
	If $\model{P} \mmodels \LPSTfo$, then $\model{P} \mmodels \phi\modalize(\vec{a})$ iff $\actualize{\model{P}} \models \phi(\vec{a})$, for any $\vec{a}$ in $\actualize{\model{P}}$'s domain and any \LTfo-formula $\phi(\vec{x})$ with free variables displayed.
\end{lem}
\begin{lem}\label{lem:levelling-potentialize}
	If $\model{A} \models \LTfo$, then $\model{A} \models \phi^r(\vec{a})$ iff $\potentialize{\model{A}} \mmodels[r] \phi(\vec{a})$, for any $\vec{a}$ from $\model{A}$'s domain, any $r$ such that $\model{A} \models \levpred(r)$, and any  \LPSTfo-formula $\phi(\vec{x})$ with free variables displayed.
\end{lem}
\begin{namedthm}[\ref{thm:ns:semantic}]
	\textcolor{white}{.}
	\begin{listn-0}
		\item[\eqref{ns:sem:a}] If $\model{P} \mmodels \LPSTfo$, then  $\actualize{\model{P}} \models \LTfo$
		\item[\eqref{ns:sem:ap}] If $\model{P} \mmodels \LPSTfo$, then there is a surjection $f$ such that $\model{P} = ({\potactualize{\model{P}}})_{f}$
		\item[\eqref{ns:sem:p}] If $\model{A} \models \LTfo$, then $\potentialize{\model{A}} \mmodels \LPSTfo$
		\item[\eqref{ns:sem:pa}] If $\model{A} \models \LTfo$, then $\model{A} = {\actpotentialize{\model{A}}}$
	\end{listn-0}
\end{namedthm}
\begin{proof}
	\emphref{ns:sem:a}
	By Theorem \ref{thm:ns:deductive}.\ref{ns:ded:m} and Lemma \ref{lem:modalize-actualize}.
	
	\emphref{ns:sem:ap}
	Let $W$ be the set of $\model{P}$'s worlds; let $L = \Setabs{s}{\actualize{\model{P}} \mmodels \levpred(s)}$ be the set of ${\potactualize{\model{P}}}$'s worlds. Using Theorem \ref{thm:pst:key}.\ref{key:maxlev}, for each $\worldvar{w} \in W$, let $f(\worldvar{w})$ be the maximal level in $\worldvar{w}$. 
	
	I claim that $f : W \functionto L$ is a surjection. To show that $L \subseteq \range{f}$, fix $s \in L$, i.e.\ $\actualize{\model{P}} \models \levpred(s)$. Let $\worldvar{w}$ be such that $\model{P} \mmodels[\worldvar{w}] \existspred{s}$; now $\model{P} \mmodels[\worldvar{w}] \levpred(s)$ by Lemmas \ref{lem:modalize-actualize} and \ref{lem:pst:inext}, and there is $\worldvar{v}$ such that $\model{P} \mmodels[\worldvar{v}] \maxlev(s)$ by Theorem \ref{thm:pst:key}.\ref{key:prune}; so $f(\worldvar{v}) = s $. The proof that $\range{f} \subseteq L$ is similar but simpler. 
	
	Now $\model{P}$ and ${\potactualize{\model{P}}}$ share a global domain, since $\poss\existspred{x}$ is a schema of our logic (see footnote \ref{fn:negativefreelogic}). They agree on membership and identity by construction. So $\model{P} = ({\potactualize{\model{P}}})_{f}$.
	
	\emphref{ns:sem:p} 
	By Theorem \ref{thm:ns:deductive}.\ref{ns:ded:l} and Lemma \ref{lem:levelling-potentialize}.

	\emphref{ns:sem:pa} 	
	By \ref{lt:strat}, $\model{A}$ and ${\actpotentialize{\model{A}}}$ have the same domain, and they agree on membership by construction.	
\end{proof}\noindent
As discussed in \S\ref{s:2:ns:2}, we also have two second-order versions of Theorem \ref{thm:ns:semantic} which hold for full or Henkin semantics. 
\begin{thm}\label{thm:ns:semantic:nec}
	Theorem \ref{thm:ns:semantic} holds for \LT and \LPSTnec, where we extend flattening and potentialization with these clauses:
	
	\emph{Flattening}:  $\actualize{\model{P}}$'s second-order domain is $\model{P}$'s global second-order domain; and $\actualize{\model{P}} \models F(\vec{a})$ iff $\model{P} \mmodels \poss F(\vec{a})$. 
	
	\emph{Potentialization}: $\potentialize{\model{A}}$'s global second-order domain is $\model{A}$'s second-order domain; and $\potentialize{\model{A}} \mmodels[s] F(\vec{a})$ iff $\model{A} \models F(\vec{a})  \land \vec{a} \subseteq s$.
\end{thm}
\begin{thm}\label{thm:ns:semantic:con}
	Theorem \ref{thm:ns:semantic} holds for \LTb and \LPSTcon, where we extend flattening and potentialization as above, \emph{and} add a further clause for potentialization, to allow variable second-order domains: $\potentialize{\model{A}} \mmodels[s] \existspred{F}$ iff $\model{A} \models F \sofoundat s$.
\end{thm}\noindent
As mentioned in \S\ref{s:2:ns:2:con}, if we invoke \emph{full} semantics, we can obtain a final semantic result. Recall that, with full semantics, first-order domains determine second-order domains. (In the modal setting: full contingentist semantics specifies that a world's monadic second-order domain is the powerset of that world's first-order domain.) So, when we are using full semantics, we can forget about second-order entities, allowing them to `take care of themselves', and simply use the definitions of flattening and potentialization that were given for \emph{first-order} theories. We then have a near-synonymy as follows:
\begin{thm}\label{thm:ns:semantic:full}
	Using full semantics, Theorem \ref{thm:ns:semantic} holds for \LT and \LPSTcon, with flattening and potentialization exactly as defined in \S\ref{s:2:ns:1:sem}.
\end{thm}
\begin{proof}
	Clauses \eqref{ns:sem:ap}--\eqref{ns:sem:pa} are left to the reader. To establish \eqref{ns:sem:a}, suppose $\model{P} \mmodels \LPSTcon$. So $\model{P} \mmodels[\worldvar{w}] \LT$ for each world $\worldvar{w}$, by Theorem \ref{thm:pst:key}. Now \LT is externally quasi-categorical by Pt.\ref{pt:lt}  Theorem \ref{thm:LTexternalcat}, and membership is modally robust by \ref{m:mem} and \ref{m:ext}. So, given any two worlds of $\model{P}$, one is an initial segment of the other. Hence $\actualize{\model{P}} \models \LT$. 
\end{proof}

\section{Equivalences concerning \LTb}\label{s:app:ltb}
In \S\ref{s:2:ET:conpot}, I considered Edna, a contingentist who holds that time is endless. To formalise the claim `time is endless', we have the modal axiom $\somelater\top$. Let \LPSTconplus be the result of adding this axiom to \LPSTcon. So, Edna's theory is \LPSTconplus. 

By contrast, consider the principle $\alllater \bot \lor \somelater\alllater\bot$. Over \LPSTcon, this amounts to the statement `time has an end'. Call this theory \LPSTconminus. 

Actualists can mirror such talk about the `end of time'. The sentence \ref{lt:cre}, from Pt.\ref{pt:lt} \S\ref{s:1:ltsubzf}, states that the (actualist) hierarchy has no last level. For brevity, let \LTbplus be \LTb + \ref{lt:cre}, and let \LTbminus be \LTb + $\lnot$\ref{lt:cre}. It is easy to confirm that \LPSTconplus is near-synonymous with \LTbplus, and that \LPSTconminus is near-synonymous with \LTbminus. 

However, \LTbplus and \LTbminus merit discussion in their own right. Fairly trivially, \LTbminus is identical to \LT + $\lnot$\ref{lt:cre}. More interestingly, \LTbplus can be regarded as a notational variant of the \emph{first-order} theory \LTfo + \ref{lt:cre}, i.e.\ \LTfoplus. Specifically: there is an interpretation which is identity over the first-order entities and bi-interpretability over the second-order entities.\footnote{Thanks to James Studd, Albert Visser, and Sean Walsh for discussion of this case.} Here is the point in detail. We interpret \LTbplus in \LTfoplus using a translation, $\sotofomark$, which tells us to regard $n$-place second-order variables as an odd way to talk about sets of $n$-tuples. Formally, its only non-trivial clauses are:
\begin{align*}
	(Y^n(x_1, \ldots, x_n))\sotofo &\coloneq \tuple{x_1, \ldots, x_n} \in Y^n\\
	(\forall Y^n\phi)\sotofo &\coloneq \forall Y^n((\forall z \in Y^n)(z\text{ is an }n\text{-tuple}) \lonlyif \phi\sotofo)
\end{align*}
where we treat $n$-tuples via Wiener--Kuratowski,\footnote{So e.g.\ $(Y^2(x_1, x_2))\sotofo$ is $(\exists z \in Y^2)\forall y(y \in z \liff (\forall w(w \in y \liff w = x_1) \lor \forall w(w \in y \liff (w = x_1 \lor w = x_2))))$.} and regard each capital, superscripted, variable as just a new first-order variable. This yields a very tight connection between \LTfoplus and \LTbplus:
\begin{thm}\label{thm:engulf}\textcolor{white}{.}
	\begin{listn-0}
		\item\label{engulf:back} $\LTfoplus \proves \phi$ iff $\LTbplus \proves \phi$,  for first-order $\phi$
		\item\label{engulf:forth} $\LTbplus \proves \phi$ iff $\LTfoplus \proves \phi\sotofo$, for second-order $\phi$
	\end{listn-0}
	Moreover, \LTbplus proves that $\sotofomark$ is identity over the first-order entities and an isomorphism over the second-order entities. 
\end{thm}
\begin{proof} 
	\emphref{engulf:back} It suffices to show that \LTbplus proves the \ref{sep} scheme. Fix a formula $\phi$. Fix $a$. By \ref{lt:strat}, there is some level $s \supseteq a$. Using \ref{ltb:comp}, there is $F \sofoundat s$ such that $(\forall x \subseteq s)(F(x) \liff \phi)$. By \ref{ext} and the \ref{sep} axiom, we have $b = \Setabs{x \in a}{F(x)}$; now $b = \Setabs{x \in a}{\phi}$, as required, since levels are transitive. 
	
	\emphref{engulf:forth} To establish \ref{ltb:comp}$\sotofo$: fix a level $s$; using \ref{lt:cre}, let $t$ be the  $2n\mathord{+}1^{\text{th}}$ level after $s$; then use the Separation scheme to obtain $F^n = \Setabs{\tuple{x_1, \ldots, x_n} \in t}{\phi\sotofo}$, noting that $\tuple{x_1, \ldots, x_n} \in t$ iff $x_{i} \subseteq s$ for all $1 \leq i \leq n$. \ref{ltb:strat}$\sotofo$ follows from \ref{lt:strat}.  
	
	To establish the `moreover' clause: in \LTbplus, define an isomorphism (which blurs types) via $\tau(F^n) = \Setabs{\tuple{x_1, \ldots, x_n}}{F^n(\vec{x})}$.
\end{proof}\noindent
Note that  Theorem \ref{thm:engulf} is \emph{not} a definitional equivalence: definitional equivalence is unavailable, since \LTfoplus and \LTbplus have different grammars. This difference aside, \LTfoplus and \LTbplus are as tightly linked as we could want. Moreover, since \LPSTconplus and \LTbplus are near-synonymous, Theorem \ref{thm:engulf} allows us to regard \LPSTconplus, which is a modal second-order theory, as a notational variant of \LTfoplus, which is a non-modal first-order theory.\footnote{Cf.\ \textcite[179--80]{Studd:EML}.}

\section*{Acknowledgements}
Special thanks to \O{}ystein Linnebo and James Studd, for extensive discussion. Thanks also to Sharon Berry, Geoffrey Hellman, Luca Incurvati, Will Stafford, Rob Trueman, Albert Visser,  Sean Walsh, and anonymous referees for \emph{Bulletin of Symbolic Logic}.

\stopappendix

\chapter[Part 3]{Level Theory, Part \thechapter
	\chapsubhead{A boolean algebra of sets arranged in well-ordered levels}}\label{pt:blt}

\noindent\textcolor{blue}{This document contains preprints of Level Theory, Parts 1--3. All three papers are forthcoming at \emph{Bulletin of Symbolic Logic}.}

\begin{quote}
	\textbf{Abstract.} On a very natural conception of sets, every set has an absolute complement. The ordinary cumulative hierarchy dismisses this idea outright. But we can rectify this, whilst retaining classical logic. Indeed, we can develop a boolean algebra of sets arranged in well-ordered levels. I show this by presenting Boolean Level Theory, which fuses ordinary Level Theory (from Part \ref{pt:lt}) with ideas due to Thomas Forster, Alonzo Church, and Urs Oswald. BLT neatly implement Conway's games and surreal numbers; and a natural extension of BLT is definitionally equivalent with ZF.
\end{quote}

\begin{epigraph}
	{Like all walls it was ambiguous, two-faced. What was inside it and what was outside it depended upon which side you were on.}
	{\textcite[1]{LeGuin:D}}
\end{epigraph}
\noindent 
Building on work by Alonzo Church and Urs Oswald, Thomas Forster has provided a pleasingly different way to think about sets. As in the ordinary cumulative hierarchy, the sets are stratified into well-ordered levels. But, unlike the ordinary cumulative picture, the sets form a boolean algebra. In particular, every set has an absolute complement, in the sense that $\forall a \exists c \forall x(x \in a \liff x \notin c)$. In this paper, I develop an axiomatic theory for this conception of set: Boolean Level Theory, or \BLT.

I start by outlining the bare-bones idea of a complemented hierarchy of sets, according to which sets are arranged in stages, but where each set is found alongside its complement. I axiomatize this bare-bones story in the most obvious way possible, obtaining Boolean Stage Theory, \BST. It is clear that any complemented hierarchy satisfies \BST (see \S\S\ref{s:3:story}--\ref{s:3:bst}). Unfortunately, \BST has multiple primitives. To overcome this, I develop Boolean Level Theory, \BLT. The only primitive of \BLT is $\in$, but \BLT and \BST say exactly the same things about sets. As such, any complemented hierarchy satisfies \BLT. Moreover, \BLT is quasi-categorical (see \S\S\ref{s:3:blt}--\ref{s:3:quasicat}). I then provide two interpretations using \BLTzf (an obvious extension of \BLT): we can regard \ZF as a proper part of \BLTzf; but \ZF is definitionally equivalent to \BLTzf (see \S\S\ref{s:3:helow}--\ref{s:3:de}). I close by explaining how to implement Conway's games and surreal numbers in \BLT (see \S\ref{s:3:conway}).

This paper is the third in a triptych. It closely mirrors Part \ref{pt:lt}, but can be read in isolation. Let me repeat, though, that Part \ref{pt:lt} is hugely indebted to the work of Dana Scott, Richard Montague, George Boolos, John Derrick, and Michael Potter; this paper inherits those debts.\footnote{See in particular \textcites[139]{Montague:STHOL}[\S22]{MontagueScottTarski}{Scott:NRST}{Scott:AST}[8--11]{Boolos:ICS}[]{Boolos:IA}[16--22]{Potter:S}[ch.3]{Potter:STP}.}

Some remarks on notation (which is exactly as in Pt.\ref{pt:lt} \S\ref{s:1:prelim}). I use second-order logic throughout. Mostly, though, this is just for convenience. Except when discussing quasi-categoricity (see \S\ref{s:3:quasicat}), any second-order claim can be replaced with a first-order schema in the obvious way. I use some simple abbreviations (where $\Psi$ can be any predicate whose only free variable is $x$, and $\lhd$ can be any infix predicate):
\begin{align*}
	(\forall x : \Psi)\phi&\coloneq \forall x(\Psi(x) \lonlyif \phi) & 	(\forall x \lhd y)\phi&\coloneq\forall x(x \lhd y \lonlyif \phi)\\
	(\exists x : \Psi)\phi&\coloneq\exists x(\Psi(x) \land \phi)& 
	(\exists x \lhd y)\phi&\coloneq\exists x(x \lhd y \land \phi)
\end{align*}
I also concatenate infix conjunctions, writing things like $a \subseteq r \in s \in t$ for $a \subseteq r \land r \in s \land s \in t$. And I run these devices together; so $(\forall x \notin x \in a)x \subseteq a$ abbreviates $\forall x((x \notin x \land x \in a) \lonlyif x \subseteq a)$. When I announce a result or definition, I list in brackets the axioms I am assuming. For readability, all proofs are relegated to the appendices.

\section{The Complemented Story}\label{s:3:story}
Here is a very natural image of sets: \emph{sets are not just collections of objects; sets partition the universe, and both sides of the partition yield a set.} There is the set of sheep; and there is the set of non-sheep. There is the set of natural numbers; and there is the set of everything else. There is the empty set; and there is the universal set. 

Many will reject this image out of hand. Supposedly, the paradoxes of na\"ive set theory have taught us that there is no universal set; for if there were a universal set $V = \Setabs{x}{x=x}$, then Separation would entail the existence of the Russell set $\Setabs{x}{x \notin x}$, which is a contradiction.

That reasoning, though, is too quick. Separation is incompatible with the existence of $V$.\footnote{\label{fn:classicalassume}NB: I assume classical logic throughout.}  More generally, Separation is incompatible with the principle of Complementation (i.e.\ with the principle that every set has an absolute complement). But it does not immediately follow that Complementation is false; only that we must choose between Separation and Complementation. 

Both principles are very natural. Separation, however, has the weight of history behind it; and this might not merely be a historical accident. There is a serious argument in favour of Separation and against Complementation, which runs as follows. The paradoxes of na\"ive set theory forced us to develop a less na\"ive conception of \emph{set}. The best such conception (according to this argument) is the cumulative iterative conception, as articulated by this bare-bones story {(recycled from Pt.\ref{pt:lt}):}
\begin{storytime}{\textbf{The Basic Story.}} Sets are arranged in stages. Every set is found at some stage. At any stage \stage{s}: for any sets found before \stage{s}, we find a set whose {members} are exactly those sets. We find nothing else at \stage{s}.
\end{storytime}\noindent
It is easy to see that this conception of set yields Separation rather than Complementation: any subset of a set $a$ occurs at (or before) any stage at which $a$ itself occurs. So (the argument concludes) we should embrace Separation and reject Complementation.

I take this argument very seriously. However, its success hinges on whether the ordinary cumulative iterative conception really is the `best' conception of \emph{set}. Whatever exactly `best' is supposed to mean, the argument lays down a challenge: produce an equally good or better conception of \emph{set}, which accepts Complementation and rejects Separation.

This paper considers a very specific reply to this challenge, due to Forster's development of work by Church and Oswald.\footnote{\textcites{Church:STUS}{Oswald:PhD}; see also \textcites{Mitchell:PhD}{Sheridan:VCST}. \textcite{Forster:CSTUS} includes a nice summary of the technicalities behind the original Church--Oswald idea.}  Forster's idea is to make a small tweak to the story of the ordinary hierarchy, so that `each time we [find] a new set\ldots we also [find] a companion to it which is to be its complement'.\footnote{\textcite[100]{Forster:ICS}. Note that I speak of `finding' sets, whereas Forster speaks of `creating' them. Talk of `creation' leads Forster to say that the members of $V$ change, stage-by-stage, as more sets are created, so that $V$ is `intensional', in a way that $\emptyset$ is not \parencite*[100]{Forster:ICS}. I think that Forster should regard $\emptyset$ as equally `intensional', since what $\emptyset$ omits changes, stage-by-stage. However, if sets are discovered (rather than created) stage-by-stage, then all issues concerning intensionality can be side-stepped: all that changes, stage-by-stage, is our knowledge about $V$'s members and $\emptyset$'s non-members.
	
If we admit contingently-existing urelements, then the discussion of intensionality becomes much more complicated. In the actual world, $\text{Boudica} \in \Setabs{x}{x = x}$; but in a possible world where she never existed, $\text{Boudica} \notin \Setabs{x}{x = x}$; by contrast, in all possible worlds, $\text{Boudica} \notin \Setabs{x}{x \neq x}$. From this, one might infer that $V$ is intensional whereas $\emptyset$ is not. But this inference is not immediate; it requires two substantial, further, assumptions: (1) that the descriptions `$\Setabs{x}{x \neq x}$' and `$\Setabs{x}{x=x}$' \emph{rigidly} designate $\emptyset$ and $V$ respectively, and (2) that intensionality concerns trans-world variation of \emph{members} rather than trans-world variation of \emph{non-members}. I hope to explore both assumptions elsewhere. (Thanks to James Studd, Timothy Williamson, Stephen Yablo, and an anonymous referee for this journal, for pushing me on this point.)} In slightly more detail, we offer the following bare-bones story:
\begin{storytime}{\textbf{The Complemented Story.}} Sets are arranged in stages. Every set is found at some stage. At any stage \stage{s}: for any sets found {before} \stage{s}, we find both
	\begin{listclean}
		\item[\lowclause] a set whose \emph{members} are exactly those sets, and 
		\item[\highclause] a set whose \emph{non-members} are exactly those sets.
	\end{listclean} 
	We find nothing else at \stage{s}.
\end{storytime}\noindent
According to our new story, we find each set using either clause \lowclause or clause \highclause. Moreover, if we find a set using clause \lowclause, then we find its {absolute} complement using clause \highclause, and vice versa.  {This is the absolute complement since, in clause \highclause, we quantify over all sets that will ever be discovered, not just those discovered before stage \stage{s}.} This story therefore secures Complementation; it describes the bare idea of a \emph{complemented} hierarchy of sets. But it only describes the \emph{bare} idea, since, for example, it says nothing about the height of the hierarchy.

In what follows, I will develop an axiomatic theory of this story, and explore that theory's behaviour. To be clear: I am not claiming that we should reject the ordinary hierarchy in favour of the complemented. My aim is only to provide a coherent (and surprisingly elegant) conception of \emph{set} which allows for Complementation rather than Separation.

In what follows, I will speak of \emph{low} sets and \emph{high} sets.\footnote{Note that every set will be low or high. This terminology departs somewhat from Church's. \textcite[298]{Church:STUS} defined `a \emph{low} set as a set which has a one-to-one relation with a well-founded set' and `a \emph{high} set as a set which is the complement of a low set'. This leaves logical space for sets which are neither low nor high (in Church's terms), and \textcite[305]{Church:STUS} used such sets to provide a Frege--Russell definition of cardinal numbers.}  A set is low iff we find it using clause \lowclause; we characterize low sets by saying `exactly these things, which we found earlier, are this set's members'. The limiting case of a low set is the empty set, $\emptyset$. A set is high iff we find it using clause \highclause; we characterize high sets by saying `exactly these things, which we found earlier, are omitted from this set'. The limiting case of a high set is the universe, $V$. (Note that low sets can have high sets as members, e.g.\ $\{V\}$ would be a low set with a high member.)

\section{Boolean Stage Theory}\label{s:3:bst}
Given a model of \ZF, there are simple methods for constructing models of the complemented hierarchy.\footnote{See \textcites[\S\S1--2]{Forster:CSTUS}[106--8]{Forster:ICS}; and my interpretation $\cointer$ in \S\ref{s:interpret:bltinterpretszf}.} However, if the idea of a complemented hierarchy is genuinely to rival that of the ordinary hierarchy, it cannot remain parasitic upon \ZF; it needs a fully autonomous {theory}. I will provide such a theory over the next two sections.\footnote{The approach in this section follows Scott and Boolos, but in the setting of complemented hierarchies rather than the ordinary hierarchies; see Pt.\ref{pt:lt} \S\S\ref{s:1:st} and \ref{s:1:history}.}

The Complemented Story, which introduces the bare-bones idea of a complemented hierarchy, speaks of both stages and sets. To begin, then, I will present a theory which quantifies distinctly over both sorts of entities. Boolean Stage Theory, or \BST, has two distinct sorts of first-order variable, for \emph{sets} (lower-case italic) and for \textbf{stages} (lower-case bold). It has five primitive predicates:
\begin{listbullet}
	\item[$\in$:] a relation between sets; read `$a \in b$' as `$a$ is in $b$'
	\item[$<$:] a relation between stages; read `$\stage{r} < \stage{s}$' as `$\stage{r}$ is before $\stage{s}$'
	\item[{$\foundat$}:] a relation between a set and a stage; read {`$a \foundat \stage{s}$'} as `$a$ is found at $\stage{s}$'
	\item[$\lowpred$:] a property of sets; read `$\lowpred(a)$' as `$a$ is low', i.e.\ we find $a$ using clause \lowclause
	\item[$\highpred$:] a property of sets; read `$\highpred(a)$' as `$a$ is high', i.e.\ we find $a$ using clause \highclause
\end{listbullet}
For brevity, I write $a \foundby \stage{s}$ for $\exists \stage{r}(a \foundat \stage{r} < \stage{s})$, i.e.\ $a$ is found before $\stage{s}$. Then \BST has {eight} axioms:\footnote{\label{fn:cheapBST}Using classical logic yields `cheap' proofs of the existence of a stage, an empty set, and a universal set, via \ref{bst:stage}, \ref{bst:find0} and \ref{bst:find1}. Those who find such proofs \emph{too} cheap might wish to add some explicit existence axioms. (Cf.\ Pt.\ref{pt:lt} footnote \ref{fn:lt:cheap}.)}
\begin{listaxiom}
	\labitem{\ref{ext}}{bst:ext} 
	$\forall a \forall b (\forall x(x \in a \liff x \in b) \lonlyif a = b)$
	\labitem{\ref{st:ord}}{bst:ord} 
	$\forall \stage{r} \forall \stage{s}\forall \stage{t}(\stage{r} < \stage{s} < \stage{t} \lonlyif \stage{r} < \stage{t})$
	\labitem{\ref{st:stage}}{bst:stage} 
	$\forall a \exists \stage{s}\, \ a \foundat \stage{s}$
	\labitem{Cases}{bst:cases} 
	$\forall a(\lowpred(a) \lor \highpred(a))$
	\labitem{Priority$_{\lowindicator}$}{bst:pri0} 
	$\forall \stage{s} (\forall a : \lowpred)(a \foundat \stage{s} \lonlyif (\forall x \in a)x \foundby \stage{s})$
	\labitem{Priority$_{\highindicator}$}{bst:pri1} $\forall \stage{s}(\forall a: \highpred)(a \foundat \stage{s} \lonlyif (\forall x \notin a)x \foundby \stage{s})$
	\labitem{Specification$_{\lowindicator}$}{bst:find0} 
	$\forall F \forall \stage{s}((\forall x : F)x \foundby \stage{s} \lonlyif (\exists a : \lowpred)(a \foundat \stage{s} \land \forall x(F(x) \liff x\in a)))$
	\labitem{Specification$_{\highindicator}$}{bst:find1} 
	$\forall F \forall \stage{s}((\forall x : F)x \foundby \stage{s} \lonlyif (\exists a : \highpred)(a \foundat \stage{s} \land \forall x(F(x) \liff x \notin a)))$
\end{listaxiom}
I will now explain how to justify each axiom. 

The first two axioms make implicit assumptions explicit. Whilst I did not mention \ref{ext} when I told the story of the complemented hierarchy, I take it as analytic that sets are extensional.\footnote{For brevity of exposition, I am considering hierarchies of pure sets.}  Similarly, \ref{bst:ord} records the analytic fact that `before' is transitive. Note, though, that I do not {explicitly} assume that the stages are {well-ordered},\footnote{Here I part company with \textcite[100]{Forster:ICS}, who {explicitly} stipulates that the stages are well-ordered. Ultimately, \BST proves a well-ordering result (Theorem \ref{thm:ecs:wo}).}  as it is unclear at this point what would justify that assumption. (After all, if we are willing to countenance entities as ill-founded as $V$, then it is not immediately obvious that we should refuse to countenance a hierarchy with infinite descending chains of stages. And the Complemented Story does not {explicitly require} that the stages be well-ordered.)

Informally, \ref{bst:stage} says that every set is discovered at some stage; this claim appears verbatim in the Complemented Story. Likewise, \ref{bst:cases} says that every set is either low or high, and this is immediate from the fact that every set is discovered using either clause \lowclause or clause \highclause. (Note, though, that I do not assume at the outset that this is an exclusive disjunction; initially, we should be open to the thought that one set could be discovered using both clauses.)\footnote{Ultimately, \BST proves that no set is discovered using {both} clauses (Lemma \ref{lem:bst:lowin}).}

Next, \ref{bst:pri0} and \ref{bst:pri1} say that if we find a low set at a stage, then we find all its members earlier, and if we find a high set at a stage, then we find all its non-members earlier; both claims follow from clauses \lowclause and \highclause. Finally, {\ref{bst:find0} and \ref{bst:find1}} say that if every $F$ was found before a certain stage, then at that stage we find both the low set of all $F$s, and the high set of all non-$F$s; again, both claims follow from \lowclause and \highclause. 

Since all eight axioms hold of {the Complemented Story}, any complemented hierarchy satisfies \BST.

\section{Boolean Level Theory}\label{s:3:blt}
Unfortunately, \BST contains rather a lot of primitives. Fortunately, most of them can be eliminated. In this section, I present Boolean Level Theory, or \BLT. This theory's only primitive is $\in$, but it makes exactly the same claims about sets as \BST does.\footnote{The approach in this section mirrors Pt.\ref{pt:lt} \S\S\ref{s:1:lt} and \ref{s:1:LTST}, which builds on work by Montague, Scott, Derrick and Potter; see also Pt.\ref{pt:lt} \S\ref{s:1:history}.} I start with a key definition:\footnote{Compare Montague's and Scott's $\pot$-operation, presented in Pt.\ref{pt:lt} Definition \ref{def:pot}.}
\begin{define}\label{def:bpot}
	For any set $a$, let $a$'s \emph{absolute complement} be $\complement{a} = \Setabs{x}{x \notin a}$, if it exists. Let $\bpot{a} = \Setabs{x}{(\exists c \notin c \in a)(x \subseteq c \lor \complement{x} \subseteq c)}$, if it exists.\footnote{By the notational conventions, $\bpot{a} = \Setabs{x}{\exists c(c \in a \land c \notin c \land (x \subseteq c \lor \complement{x} \subseteq c))}$. \BLT's axiom \ref{blt:comp} guarantees that $\complement{a}$ exists for every $a$. However, we do not initially assume that $\bpot{a}$ exists for every $a$; instead, we initially treat every expression of the form `$b = \bpot{a}$' as shorthand for `$\forall x(x \in b \liff (\exists c \notin c \in a)(x \subseteq c \lor (\exists z \subseteq c)\forall y(y \in z \liff y \notin x)))$', and must double-check whether $\bpot{a}$ exists. Ultimately, though, \BLT proves that $\bpot{a}$ exists for every $a$: if $a \notin a$ then $\bpot{a} \subseteq \bevof{a}$ (see Definition \ref{def:bevof}); if $a \in a$ then $\bpot{a} = V$.}
\end{define}\noindent
The definition of $\complement{a}$ needs no comment, but the definition of $\bpot{a}$ merits explanation. It turns out that \BST proves that $a$ is low iff $a \notin a$, and $a$ is high iff $a \in a$ (see Lemma \ref{lem:bst:lowin}). Seen in this light, $\bpot{a}$ collects together all the subsets of low members of $a$, and all the complements of such subsets. As a specific example, if $b$ is low, then $\bpot{\{b\}} = \Setabs{x}{x \subseteq b \lor \complement{x} \subseteq b}$, i.e.\ it is the result of closing $b$'s powerset under complements. We use this operation in this next definition (where `\boohist' is short for `boolean-history', and `\boolevel' is short for `boolean-level'):\footnote{Compare Pt.\ \ref{pt:lt} Definition \ref{def:level}, which simplifies the Derrick--Potter definition of `level'. Here, `\boohist' is short for `boolean-history'; `\boolevel' is short for `boolean level'.}
\begin{define}\label{def:bevel}
	Say that $h$ is a \emph{\boohist}, written $\bistpred(h)$, iff $h\notin h \land (\forall x \in h)x = \bpot(x \cap h)$. Say that $s$ is a \emph{\boolevel}, written $\bevpred(s)$, iff $(\exists h : \bistpred) s = \bpot{h}$.
\end{define}\noindent
The intuitive idea behind Definition \ref{def:bevel} is that the \boolevel{}s go proxy for the stages of the Complemented Story, and each \boohist is an initial sequence of \boolevel{}s. (It is far from obvious that these definitions work as described, but we will soon see that they do.) Using these definitions, \BLT has just four axioms:\footnote{\label{fn:emptyandV}As in footnote \ref{fn:cheapBST}, classical logic yields a `cheap' proof of the existence of $\emptyset$ and $V$.}
\begin{listaxiom}
	\labitem{\ref{ext}}{yetmoreext} $\forall a \forall b(\forall x(x \in a \liff x \in b) \lonlyif a = b)$
	\labitem{Complements}{blt:comp} {$\forall a(\exists c = \complement{a})(a \notin a \liff c \in c)$}
	\labitem{Separation$_\notin$}{blt:sep}$\forall F(\forall a \notin a)(\exists b \notin b)\forall x(x \in b \liff (F(x) \land x \in a))$
	\labitem{Stratification$_\notin$}{blt:strat} $(\forall a \notin a)(\exists s : \bevpred) a \subseteq s$
\end{listaxiom}\noindent
Intuitively, \ref{blt:comp} tells us that every set has a complement, and a set is low iff its complement is high; \ref{blt:sep} tells us that arbitrary subsets of low sets exist (and are low); and \ref{blt:strat} tells us that every low set is a subset of some \boolevel (which corresponds to the thought that it is found at some stage). These axioms and definitions are vindicated by this next result, which shows that \BLT has exactly the same set-theoretic content as \BST (see \S\ref{s:thm:BLTBST} for the proof):
\begin{thm}\label{thm:BLTBST}
	$\BST \proves \phi$ iff $\BLT \proves \phi$, for any \BLT-sentence $\phi$.
\end{thm}\noindent
Otherwise put: no information about sets is gained or lost by moving between \BST and \BLT. Moreover, since every complemented hierarchy satisfies \BST, every complemented hierarchy satisfies \BLT. In what follows, then, I will treat \BLT as the canonical theory of complemented hierarchies. 

\section{Characteristics and extensions of \BLT}\label{s:3:characteristic}
To give a sense of how \BLT behaves, I will state some of its `characteristic' results (the proofs are in \S\ref{s:proofs:blt:elementary}). The first two results allow us to characterize \BLT with a simple slogan: \emph{a boolean algebra of sets arranged in well-ordered levels.}
\begin{thm}[\BLT]\label{thm:ecs:wo}
	The \boolevel{}s are well-ordered by $\in$.
\end{thm}
\begin{thm}[\BLT]\label{thm:blt:boolean}
	The sets form a boolean algebra under complementation, $\cap$ and $\cup$.
\end{thm}
\noindent
This first result is quite surprising:\footnote{It will be much less surprising for those who have read Pt.\ref{pt:lt} \S\ref{s:1:wodiscussion}.} the Complemented Story does not \emph{explicitly} specify that the stages must be well-ordered (see \S\ref{s:3:bst}); but, since every complemented hierarchy satisfies \BLT (see \S\ref{s:3:blt}), every complemented hierarchy has well-ordered levels.

The well-ordering of the \boolevel{}s yields a powerful tool, which intuitively allows us to consider the \boolevel at which a set is first found:
\begin{define}[\BLT]\label{def:bevof}
	If $a \notin a$, let $\bevof{a}$ be the $\in$-least \boolevel with $a$ as a subset; i.e., $a \subseteq \bevof{a}$ and $\lnot (\exists s : \bevpred)a \subseteq s \in \bevof{a}$. If $a \in a$, let $\bevof{a} = \bevof{\complement{a}}$.
\end{define}\noindent
Note that $\bevof{a}$ exists for any $a$, by \ref{blt:strat}, \ref{blt:comp} and Theorem \ref{thm:ecs:wo}.

A third characteristic result is that there is a \emph{contra-automorphism} on the universe.\footnote{See \textcite[Definition 16 and subsequent comments]{Forster:CSTUS}. This result inspires my epigraph, from Le Guin. I owe the point to Brian King: in 2006, he arrived at an idea like the Complemented Story (independently of Forster) and explained it using Le Guin's image.} Roughly put: replacing membership with non-membership (and vice versa) yields an isomorphic universe. Formally:
\begin{define}\label{def:gameneg}
	We recursively define $a$'s \emph{negative}, written $\gameneg a$, as follows:
	\begin{align*}
		\gameneg a &\coloneq 
		\complement{\Setabs{\gameneg x}{x \in a}} \text{, if }a \notin a &
		\gameneg a&\coloneq 
		\Setabs{\gameneg x}{x \notin a} \text{, if }a \in a
	\end{align*}
\end{define}
\begin{thm}[\BLT]\label{thm:blt:antimorphism}
	$\forall a \forall b(a \in b \liff \gameneg a \notin  \gameneg b)$
\end{thm}\noindent 
This immediately yields a nice duality:
\begin{cor}[\BLT] 
	$\phi \liff \phi\leguin$, for any \BLT-sentence $\phi$, where $\phi\leguin$ is the sentence which results from $\phi$ by replacing every `$\in$' with `$\notin$' and vice versa.
\end{cor}\noindent
These results highlight some of \BLT's deductive strengths. Now let me comment on its (deliberate) weakness. By design, \BLT axiomatizes only the \emph{bare} idea of a complemented hierarchy, and so makes no comment on the hierarchy's height.\footnote{Beyond the fact that classical logic guarantees the existence of at least one stage; see footnotes \ref{fn:cheapBST} and \ref{fn:emptyandV}.} If we want to ensure that our hierarchy is reasonably tall, three axioms suggest themselves (where `$P$' is a second-order function-variable in the statement of \ref{blt:rep}):
\begin{listaxiom}
	\labitem{Endless$_\notin$}{blt:cre} $(\forall s : \bevpred)(\exists t : \bevpred)s \in t$
	\labitem{Infinity$_\notin$}{blt:inf} $(\exists s : \bevpred)((\exists q : \bevpred)q \in s \land 
	(\forall q : \bevpred)(q \in s \lonlyif (\exists r : \bevpred)q \in r \in s))$
	\labitem{Unbounded$_\notin$}{blt:rep} $\forall P(\forall a \notin a)(\exists s : \bevpred)(\forall x \in a) P(x) \in s$
\end{listaxiom}\noindent
\ref{blt:cre} says there is no last \boolevel. \ref{blt:inf} says that there is an infinite \boolevel, i.e.\ a \boolevel with no immediate predecessor. \ref{blt:rep} states that the hierarchy of \boolevel{}s is so tall that no low set can be mapped unboundedly into it (recall that the low sets are precisely the non-self-membered sets). 

To make all of this more familiar, here are some simple facts relating \BLT to \ZF. Let \BLTplus stand for $\BLT + \clearme{\ref{blt:cre}}$, and \BLTzf stand for $\BLT + \clearme{\ref{blt:inf}} + \clearme{\ref{blt:rep}}$; then:\footnote{Since \BLTplus proves Pairing, \BLTplus extends NF$_2$, the sub-theory of Quine's NF whose axioms are \ref{ext}, Pairing, and Theorem \ref{thm:blt:boolean}. However, \BLTplus does not extend NF$_\text{O}$, the theory which adds to NF$_2$ the axiom that $\Setabs{x}{a \in x}$ exists for every $a$; in particular, $\Setabs{x}{\emptyset \in x}$ does not exist; see the proof of Proposition \ref{prop:pstl:mltint}.\ref{blt:nopowersets} in  \S\ref{s:proofs:blt:elementary}. For discussion of NF$_2$ and NF$_\textspaced{O}$, see \textcite[\S2]{Forster:CSTUS}.}
\begin{prop}\label{prop:blt:facts}\textcolor{white}{.}
	\begin{listn-0}
		\item\label{blt:empty} \BLT proves the Axiom of Empty Set, i.e.\ $\exists a \forall x\phantom{(}x \notin a$.
		\item\label{blt:union} \BLT proves Union, i.e.\ $\forall a(\bigcup a\text{ exists})$.
		\item\label{blt:pairing} \BLTplus proves Pairing, i.e.\ $\forall a \forall b(\{a, b\}\text{ exists})$, but \BLT does not.
		\item \BLTplus proves Powersets-restricted-to-low-sets, i.e.\ $(\forall a \notin a)(\powerset{a}\text{ exists})$, but \BLT does not.
		\item\label{blt:nopowersets} \BLT contradicts Powersets, i.e.\ it proves $\exists a \lnot \exists  b \forall x(x \in b \liff x \subseteq a)$.
		\item\label{blt:highfoundation} \BLT proves Foundation-restricted-to-high-sets, i.e.\ $(\forall a \in a)(\exists x \in a)a \cap x = \emptyset$.
		\item\label{blt:nofoundation} \BLTplus contradicts Foundation, i.e.\ it proves $(\exists a \neq \emptyset)(\forall x \in a)a \cap x \neq \emptyset$.
		\item \BLTzf proves \ref{blt:cre}.
	\end{listn-0}
\end{prop}\noindent
If we want to state this result with maximum shock value: of the standard axioms of \ZF, \BLT validates only \ref{ext}, Empty Set, and Union (though \BLT is also consistent with Pairing and standard formulations of Infinity).

\section{The quasi-categoricity of BLT}\label{s:3:quasicat}
We have seen that every complemented hierarchy satisfies \BLT, so that every complemented hierarchy has well-ordered \boolevel{}s. In fact, we can push this point further, by noting that \BLT is quasi-categorical.\footnote{This mirrors the discussion of \LT's quasi-categoricity; see Pt.\ref{pt:lt} \S\ref{s:1:quasicat}.}
	
Informally, we can spell out \BLT's quasi-categoricity as follows: \emph{Any two complemented hierarchies are structurally identical for so far as they both run, but one may be taller than the other}. So, when we set up a complemented hierarchy, our only choice is how tall to make it.
	
In fact, there are at least two ways to explicate the informal idea of quasi-categoricity, and \BLT is quasi-categorical on both explications.\footnote{Both ways make essential use of second-order logic, albeit in different ways.} The first notion of quasi-categoricity should be familiar from Zermelo's results for \ZF, and uses the full semantics for second-order logic:
\begin{thm}\label{thm:blt:externalquasi} Given full second-order logic:
	\begin{listn-0}
		\item The \boolevel{}s of any model of \BLT are well-ordered.\footnote{i.e.\ if $\model{M} \models \BLT$ then $\Setabs{s \in M}{\model{M} \models \bevpred(s)}$ is well-ordered by $\in^\model{M}$.}
		\item For any ordinal $\alpha > 0$, there is a model of \BLT whose \boolevel{}s form an $\alpha$-sequence.\footnote{i.e.\ there is some $\model{M} \models \BLT$ such that $\Setabs{s \in M}{\model{M} \models \bevpred(s)}$ is isomorphic to $\alpha$.}
		\item\label{blt:initialseg} Given any two models of \BLT, one is isomorphic to an initial segment of the other.\footnote{When $\model{A}$ and $\model{M}$ are models of \BLT, say that $\model{A}$ is an \emph{initial segment} of $\model{M}$ iff either $\model{A} = \model{M}$ or there is some $s$ such that $\model{M} \models \bevpred(s)$ and $\model{A}$ is isomorphic to the substructure of $\model{M}$ whose domain is $\Setabs{x \in M}{\model{M} \models \bevof{x} \in s}$.}
	\end{listn-0}
\end{thm}\noindent
Since this result involves semantic ascent, it is an \emph{external} quasi-categoricity result. There is also an \emph{internal} quasi-categoricity result for \BLT, which is a theorem of the (second-order) object language, but this point requires a little more explanation.

In embracing \ref{ext}, \BLT assumes that everything is a pure set. {Here is an easy way to avoid making that assumption. Consider the following formula, which relativises \BLT to a new primitive predicate, \purepred:\footnote{Here, `$\subseteq$' and `$\bevpred$' should be defined in terms of $\varin$ rather than $\in$; similarly for `$\bevof$' in the statement of Theorem \ref{thm:blt:quasicat}.}
	\begin{align*}
		\BLT(\purepred, \varin) \coloneq {}&
		(\forall a : \purepred)(\forall b : \purepred)(\forall x (x \varin a \liff x \varin b)\lonlyif a = b) 
		\land {}\\
		& (\forall a: \purepred)(\exists c: \purepred)((\forall x : \purepred)(x \varin c \liff x \notvarin a) \land (a \notvarin a \liff c \varin c))
		\land {}\\		
		&\forall F (\forall a : \purepred)({a \notvarin a} \lonlyif {}\\
		&\phantom{\forall F (\forall a : \purepred)(}(\exists b : \purepred)(b \notvarin b \land \forall x(x \varin b \liff (F(x) \land x \varin a))))
		\land {}\\	
		& (\forall a : \purepred)(a \notvarin a \lonlyif (\exists s : \bevpred) a \subseteq s)
		\land {}\\
		& \forall x\forall y(y \varin x \lonlyif (\purepred(x) \land \purepred(y)))
	\end{align*}\noindent
The first four conjuncts say that the pure sets satisfy \BLT;\footnote{With one insignificant caveat (see footnote \ref{fn:emptyandV}): whereas classical logic guarantees that any model of \BLT contains an empty set and a universal set,  $\LT(\purepred, \varin)$ allows that there may be no pure sets.} the last says that, when we use `$\varin$', we restrict our attention to membership facts between pure sets. This avoids the assumption that everything is a pure set. Moreover, I can use this formula to state our internal quasi-categoricity result (I have labelled the lines to facilitate its explanation):\footnote{\citepossess[ch.11]{ButtonWalsh:PMT} proofs carry over straightforwardly to \BLT.}
\setcounter{equation}{0}
\begin{thm}\label{thm:blt:quasicat}This is a deductive theorem of impredicative second-order logic:
	\begin{align}
		(\BLT(\purepred_1&, \mathord{\varin_1}) \land \BLT(\purepred_2, \mathord{\varin_2})) \lonlyif {} \nonumber\\
		\exists R(&\label{spi:VtoV} \forall v \forall y(R(v,y) \lonlyif (\purepred_1(v) \land \purepred_2(y))) \land {}\\
		&\label{spi:exhausts} ((\forall v : \purepred_1)\exists y R(v,y) \lor (\forall y : \purepred_2)\exists v R(v,y)) \land{}\\
		&\label{spi:preserves} \forall v\forall y\forall x\forall z ((R(v, y) \land R(x, z)) \lonlyif (v \varin_1 x \liff y \varin_2 z))\land {}\\
		&\label{spi:functional} \forall  v \forall y\forall z ((R(v,y) \land R(v,z)) \lonlyif y = z)\land {}\\
		&\label{spi:injective} \forall v\forall x \forall  y  ((R(v,y) \land R(x, y)) \lonlyif v= x)\land {}\\
		&\label{spi:segmentsa} \forall v \forall x \forall y ((\bevof{}_1{x}  \subseteq_1 \bevof{}_1{v} \land R(v, y)) \lonlyif \exists z R(x, z))\land {}\\
		&\label{spi:segmentsb} \forall v \forall y \forall z((\bevof{}_2{z}  \subseteq_2 \bevof{}_2{y} \land R(v, y)) \lonlyif  \exists x R(x, z)))
	\end{align}
\end{thm}\noindent
Intuitively, the point is this. Suppose two people are using their versions of \BLT, subscripted with `$1$' and `$2$' respectively. Then there is some second-order entity, a relation $R$, which takes us between their sets \eqref{spi:VtoV}, exhausting the sets of one or the other person \eqref{spi:exhausts}; which preserves membership \eqref{spi:preserves}; which is functional \eqref{spi:functional} and injective \eqref{spi:injective}; and whose domain is an initial segment of one \eqref{spi:segmentsa} or the other's \eqref{spi:segmentsb} hierarchy. Otherwise put: \BLT is (internally) quasi-categorical.

As a bonus, this internal quasi-categoricity result can be lifted into an internal \emph{total}-categoricity result. To explain how, consider this abbreviation (where `$P$' is a second-order function-variable): 
	$$\somany x \Phi(x) \coloneq \exists P(\forall x \Phi(P(x)) \land (\forall y : \Phi)\exists ! x\ P(x) = y)$$
This formalizes the idea that there is a bijection between the $\Phi$s and the universe (see Pt.\ref{pt:lt} \S\ref{s:1:quasicat}). Using this notation, we can state our internal total-categoricity result:
\begin{thm}This is a deductive theorem of impredicative second-order logic:
	\begin{align*}
		(\BLT(\purepred_1, \mathord{\varin_1}) \land {} & \somany x\, \purepred_1(x) \land  \BLT(\purepred_2, \mathord{\varin_2}) \land \somany x\,  \purepred_2(x))\lonlyif {} \nonumber\\
		\exists R(&\forall v \forall y\left(R(v,y) \lonlyif {}\left(\purepred_1(v) \land \purepred_2(y)\right)\right) \land {} \\
		&(\forall  v : \purepred_1)\exists !y R(v,y) \land (\forall y : \purepred_2) \exists! v R(v,y)\land {} \\
		&\forall v\forall y\forall x\forall z \left((R(v, y) \land R(x, z)) \lonlyif (v \varin_1 x \liff y \varin_2 z)\right))
	\end{align*}
\end{thm}\noindent
Intuitively, if both \BLT-like hierarchies are as large as the universe, then there is a structure-preserving \emph{bijection} between them.

\section{Ordinary set theory as a proper part of \BLT}\label{s:3:helow}
The Complemented Story provides two clauses for finding sets. Clause \lowclause tells us that, at each stage \stage{s} and for any sets found  before \stage{s}, we find a set whose members are exactly those sets. But this is exactly what we would find according to the Basic Story (see \S\ref{s:3:story}), which deals with ordinary, uncomplemented hierarchies. Intuitively, then, we should be able to recover an ordinary hierarchy by considering a complemented hierarchy whilst ignoring any use of clause \highclause. This intuitive idea is exactly right; the aim of this section is to explain it carefully. 

First, I must formalize the notion of a set which we find without ever using clause \highclause. I call such sets \emph{hereditarily low}, or \emph{helow} for short. So: helow sets are low, their members are low, the members of their members are low, etc. Here is the precise definition:
\begin{define}\label{def:helow}
	Say that $a$ is \emph{helow}, or $\helow(a)$, iff there is some transitive $c \supseteq a$ such that $(\forall x \in c)x \notin x$. 
\end{define}\noindent
To restrict our attention to the ordinary (uncomplemented) hierarchy, we then just restrict our attention to the helow sets. To implement this formally, for any formula $\phi$, let $\phi\helowint$ be the formula which results by restricting all of $\phi$'s quantifiers to helow sets. Using this notation, we can then prove results of this shape: \emph{If some theory of uncomplemented hierarchies proves $\phi$, then some suitable theory of complemented hierarchies proves $\phi\helowint$.} 

To state these results precisely, we need a suitable theory of uncomplemented hierarchies. That theory is \LT, discussed in Pt.\ref{pt:lt}. In a nutshell: \LT stands to uncomplemented hierarchies exactly as \BLT stands to complemented hierarchies. I will now briefly recap \LT's key elements. To formalize the Basic Story, we define a predicate, \levpred, to capture the notion of a \emph{level} of an uncomplemented hierarchy (Pt.\ref{pt:lt} Definition \ref{def:level}); then \LT is the theory whose axioms are \ref{ext}, \ref{sep}, and \ref{lt:strat}, which states that $\forall a (\exists s : \levpred)a \subseteq s$ (see Pt.\ref{pt:lt} \S\ref{s:1:lt}). It transpires that \LT is quasi-categorical, and that every uncomplemented hierarchy satisfies \LT, no matter how tall or short it is (see Pt.\ref{pt:lt} \S\S\ref{s:1:wodiscussion}--\ref{s:1:quasicat}). If we want to secure a tall uncomplemented hierarchy, we can consider the axioms \ref{lt:cre}, \ref{lt:inf} and \ref{lt:rep} (see Pt.\ref{pt:lt} \S\ref{s:1:ltsubzf}); these are exactly like \ref{blt:cre}, \ref{blt:inf} and \ref{blt:rep} (see \S\ref{s:3:blt} of this part), except that they replace `$\bevpred$' with `$\levpred$'. Let \LTplus stand for $\LT + \clearme{\ref{lt:cre}}$; it turns out that \ZF is deductively equivalent to $\LT + \clearme{\ref{lt:inf}} + \clearme{\ref{lt:rep}}$; so \LT, \LTplus, and \ZF are three theories which axiomatize uncomplemented hierarchies, making successively stronger demands on the hierarchy's height. With this background in place, here is the result which intuitively states that the helow part of any complement hierarchy is an ordinary (uncomplemented) hierarchy (see \S\ref{s:blt:helow:appendix} for the proof):
\begin{thm}\label{thm:blt:helowinterpret} For any \LT-sentence $\phi$:
	\begin{listn-0}
		\item\label{helowint:vanilla} If $\LT \proves \phi$, then $\BLT	\proves \phi\helowint$
		\item\label{helowint:inf} If $\LTplus \proves \phi$, then $\BLTplus \proves \phi\helowint$
		\item\label{helowint:zf} If $\ZF \proves \phi$, then $\BLTzf \proves \phi\helowint$
	\end{listn-0}
\end{thm}

\section{Definitional equivalence}\label{s:3:de}
Theorem \ref{thm:blt:helowinterpret}.\ref{helowint:zf} allows us to regard \ZF as the result of restricting attention to the helow-fragment of \BLTzf's universe of sets. But we also have a much deeper interpretative result, as follows (see \S\ref{s:blt:de:appendix}):\footnote{Forster conjectured that a result of this shape should hold.}
\begin{thm}\label{thm:synonymy}
	\ZF and \BLTzf are definitionally equivalent, as are \LTplus and \BLTplus.
\end{thm}\noindent
As an immediate consequence, \ZF and \BLTzf are \emph{equiconsistent}, as are \LTplus and \BLTplus. However, definitional equivalence is much stronger than mere equiconsistency. 

Roughly, to say that two theories are definitionally equivalent is to say that each theory can define all the primitive expressions of the other, such that each theory can simulate the other perfectly, and where combining the two simulations gets you back exactly where you began.\footnote{For a precise statement of what definitional equivalence requires, see \textcite[ch.5]{ButtonWalsh:PMT}.} So, in some purely formal sense, \ZF and \BLTzf can be regarded as notational variants; as wrapping the same deductive content in different notational packaging.

One might be tempted to go further, and suggest that Theorem \ref{thm:synonymy} shows that there is \emph{no} relevant difference between \ZF and \BLTzf. That, however, would require further argument.\footnote{Compare Pt.\ref{pt:pst} \S\ref{s:2:equivalencethesis}.} Precisely because definitional equivalence is a purely formal property, it ignores all non-formal matters, and these may be philosophically significant. There is more philosophical discussion to be had about the significance of Theorem \ref{thm:synonymy}, but that must wait for another time.

\section{Conway games and surreal numbers in BLT}\label{s:3:conway}
Since \ZF and \BLTzf are definitionally equivalent, there is a sense in which each can do anything that the other can. Still, \BLTzf can do some things more easily than \ZF. This is neatly illustrated by considering John Conway's theory of games and surreal numbers.\footnote{Joel David Hamkins suggested this application of \BLT to me; many thanks to him, both for the initial suggestion, and for much subsequent correspondence.}

Consider two-player games in which players move alternately, with no element of chance, where the game must end in a win or loss. (Think of chess, but without the possibility of stalemate.) Abstractly, such games can be thought of as specifications of permissible positions: to make a move in such a game is just to select a new position which is permissible given the current game state; and you lose when it is your turn to move but there is no permissible position. (Think of being checkmated: you must move to a position where your King is not in check, but no such move is available.) Crucially, any position in any such game can be considered as a game in its own right. (Imagine the version of chess which always starts with the pieces arranged as after the Queen's Gambit in regular chess.) So every game can be regarded, abstractly, as nothing other than a specification of which games each player can move to. Otherwise put, if we call the two game-players Low and High, then a game is just a specification of \emph{low options}, i.e.\ games that Low can move to, and \emph{high options}, i.e.\ games that High can move to. 
	
The idea is very natural. However, as Conway remarked, formalizing it `in \ZF destroys a lot of its symmetry.' He therefore suggested that `the proper set theory in which to perform such a formalisation would be one with two kinds of membership': a game would just be a set with `low-members' (low options) and `high-members' (high options).\footnote{\textcite[66]{Conway:ONG}. \textcite{CoxKaye:AZF} take up this suggestion and offer an axiomatic theory with two kinds of membership; they prove it is definitionally equivalent with \ZF. By Theorem \ref{thm:synonymy}, it is definitionally equivalent with \BLTzf.} However, we can easily implement this idea in \BLT, using only \emph{one} kind of membership. We start by saying that the games are the sets, and then stipulate:
\begin{define}[\BLT]
	If $a$ is low, the set of $a$'s \emph{low options} is $\lowopt{a} \coloneq \Setabs{x \in a}{x\notin x}$; the set of $a$'s \emph{high options} is $\highopt{a} \coloneq \Setabs{x \in a}{x \in x}$. If $a$ is high, $\lowopt{a} \coloneq  \lowopt{\complement{a}}$ and $\highopt{a} \coloneq \highopt{\complement{a}}$.
\end{define}\noindent
Intuitively, then, $a$ and $\complement{a}$ represent the same game. Moreover, there is a natural algebra on the games, given as follows (I explain the definitions below):\footnote{The well-ordering of \boolevel{}s guarantees determinacy, and licenses induction and recursive definitions (see footnote \ref{fn:blt:recursion}, below). Definition \ref{def:gamealgebra} and \ref{def:surreal} are \BLT-implementations of \citepossess[chs.0--1]{Conway:ONG} definitions. (As defined, the sum of two low sets is always low; an arbitrary choice was required.) For Theorem \ref{thm:Conway}, see \citepossess[78]{Conway:ONG}; for Theorem \ref{thm:Surreal}, see \textcite[ch.1]{Conway:ONG}. For an accessible presentation, see also \textcite[\S\S2--4]{SchleicherStoll:ICGN}. }
\begin{define}[\BLT]\label{def:gamealgebra}
	With $\gameneg$ as in Definition \ref{def:gameneg}, define $+$ and $\leq$ recursively:
	\begin{align*}
		a + c &\coloneq 
		\Setabs{x + c}{x \in \lowopt{a}} \cup \Setabs{a + x}{x \in \lowopt{c}} \cup \Setabs{\complement{y + c}}{y \in \highopt{a}} \cup \Setabs{\complement{a + y}}{y \in \highopt{c}}\\
		a \gameleq c &\text{ iff } (\forall y \in \highopt{c})y \gamenleq a \land (\forall x \in \lowopt{a})c \gamenleq x
	\end{align*}
	We stipulate that $a \gameeq c\text{ iff }a \gameleq c \gameleq a$, and define $a - c \coloneq a + (\gameneg c)$.
\end{define}\noindent
We can make these algebraic operations intuitive as follows. To take the \emph{negative} of a game is to reverse the players' roles (cf.\ Theorem \ref{thm:blt:antimorphism}). To \emph{add} two games is to place them side-by-side, allowing a player to move in one game without affecting the other. But the \emph{partial-order} requires slightly more explanation. Suppose High plays first on the game $a$; then Low has a winning strategy \emph{iff} whatever move High makes, i.e.\ for all $y \in \highopt{a}$, if Low plays first on $y$ then High has no winning strategy. Similarly, suppose Low plays first on $a$; then High has a winning strategy \emph{iff} for all $x \in \lowopt{a}$, if High plays first on $x$ then Low has no winning strategy. So, if we gloss `$\emptyset \leq z$' as `Low has a winning strategy as second player on $z$' and gloss `$z \leq \emptyset$' as `High has a winning strategy as second player on $z$', this motivates two important special cases of the partial order: 
\begin{align*}
	\emptyset \leq a&\text{ iff }(\forall y \in \highopt{a})y \nleq \emptyset&
	a \leq \emptyset& \text{ iff }(\forall x \in \lowopt{a})\emptyset \nleq x
\end{align*}
The remainder of the definition is then set up so that $a - b \leq \emptyset$ iff $a \leq b$. More generally, we have the following foundational result:
\begin{thm}[\BLT]\label{thm:Conway}
	The sets form a partially-ordered abelian Group, with $\emptyset = 0$ and $+, -, \gameleq$ as in Definition \ref{def:gamealgebra}, all modulo $\gameeq$.\footnote{To quotient by $\gameeq$, define $[a]\coloneq  \Setabs{b \gameeq a}{(\forall x \gameeq a)\bevof{b} \subseteq \bevof{x}}$; cf.\ \textcites{Scott:DAST}[65]{Conway:ONG}.}
\end{thm}\noindent
We can obtain a totally-ordered Field by restricting our attention to surreals:  
\begin{define}[\BLT]\label{def:surreal}
	We specify that $a$ is \emph{surreal} iff: for all $x \in \lowopt{a}$ and all $y \in \highopt{a}$, both $x$ and $y$ are surreal and $x \ngeq y$. We define multiplication on surreals thus: 
	\begin{align*}
		a \cdot c \coloneq {}& \Setabs{x \cdot c + a \cdot y - x \cdot y}{(x \in \lowopt{a} \land y \in \lowopt{c}) \lor (x \in \highopt{a} \land y \in \highopt{c})} \cup {}\\
		& \Setabs{\complement{x\cdot c + a\cdot y - x\cdot y}}{(x \in \lowopt{a} \land y \in \highopt{c}) \lor (x \in \highopt{a} \land y \in \lowopt{c})}
	\end{align*}	
	We say that $a$ is a \emph{surreal-ordinal} iff $a$ is both helow and surreal.
\end{define}
\begin{thm}[\BLT]\label{thm:Surreal}
	The surreals form a totally-ordered Field, modulo $\gameeq$.
\end{thm}\noindent
Summing up: Conway's beautifully rich, nonstandard, theory of surreal numbers is available, essentially off-the-shelf, within \BLT.

\section{Conclusion}
The Complemented Story lays down a conception of \emph{set} which rivals the (ordinary) cumulative notion, but which accepts Complementation and rejects Separation (see \S\ref{s:3:story}). 

I have shown that any complemented hierarchy satisfies \BLT (see \S\S\ref{s:3:bst}--\ref{s:3:blt}). So, given the characteristic results of \BLT, the sets of any complemented hierarchy are arranged into well-ordered \boolevel{}s, and constitute a boolean algebra (see \S\ref{s:3:characteristic}). Moreover, \BLT is quasi-categorical (see \S\ref{s:3:quasicat}); so our only choice, in setting up a complemented hierarchy, is how tall to make it.

The theory \BLTzf arises from \BLT just by adding axioms which state that the complemented hierarchy \emph{is} quite tall (see \S\ref{s:3:characteristic}). And we can regard \ZF as either a proper part of \BLTzf (see \S\ref{s:3:helow}), or as a notational variant (in a purely formal sense) of \BLTzf (see \S\ref{s:3:de}). But both interpretations suggest that there is no obvious \emph{a priori} reason to favour Separation over Complementation. And in some settings, such as the discussion of Conway games, using Complementation is extremely natural (see \S\ref{s:3:conway})

\startappendix

\section{Characteristics of \BLT}\label{s:proofs:blt:elementary}
The remainder of this paper consists of proofs of the results discussed in the main text. Many of the simpler proofs are similar to results for Pt.\ref{pt:lt}; in such cases, I omit the proof and refer interested readers to the appropriate result from Pt.\ref{pt:lt}.

This first appendix deals with the results from \S\ref{s:3:characteristic}. Initially, I will work in \ECS, the subtheory of \BLT whose only axioms are \ref{ext}, \ref{blt:comp} and \ref{blt:sep} (see \S\ref{s:3:blt}). I start with some simple results and definitions:
\begin{lem}[\ECS]\label{lem:ecs:lowhigh}
	If $c \subseteq a \notin a$, then $c \notin c$; if $a \in a \subseteq c$, then $c \in c$.
\end{lem}
\begin{proof}
	If $c \subseteq a \notin a$, then $c \notin c = \Setabs{x \in a}{x \in c}$ by \ref{blt:sep} and \ref{ext}. If $a \in a \subseteq c$, then $\complement{c} \subseteq \complement{a} \notin \complement{a}$ by \ref{blt:comp}, so that $ \complement{c} \notin \complement{c}$ as before, and $c \in c$ by \ref{blt:comp}.
\end{proof}
\begin{define}\label{def:bltpot}\label{def:blttrans}
	Say that $a$ is \emph{\bltpot} iff $\forall x(\exists c(x \subseteq c \notin c \in a) \lonlyif x \in a)$. Say that $a$ is \emph{\blttrans} iff $(\forall x \notin x \in a)x\subseteq a$. Say that $a$ is \emph{complement-closed} iff $\forall x(x \in a \liff \complement{x} \in a)$. 
\end{define}
\begin{lem}[\ECS]\label{lem:ecs:bwiden}
	{If $\bpot{a}$ exists (see Definition \ref{def:bpot})}, then:
	\begin{listn-0}
		\item\label{bwiden:elem} $(\forall x \notin x \in \bpot{a})\exists c(x \subseteq c \notin c \in a)$.
		\item\label{bwiden:sub} $\bpot{a}$ is \bltpot.
		\item\label{bwiden:comp}  $\bpot{a}$ is complement-closed.
	\end{listn-0}
\end{lem}
\begin{proof}
	\emphref{bwiden:elem} Fix $x \notin x \in \bpot{a}$; so for some $c \notin c \in a$, either $x \subseteq c$ or $\complement{x} \subseteq c$ . But $\complement{x} \in \complement{x}$ by \ref{blt:comp}, so $\complement{x} \nsubseteq c$ by Lemma \ref{lem:ecs:lowhigh}.
	
	\emphref{bwiden:sub} Fix $x \subseteq c \notin c \in \bpot{a}$; so $x \subseteq c \subseteq b \notin b \in a$ for some $b$ by \eqref{bwiden:elem}; hence $x \in \bpot{a}$. 
	
	\emphref{bwiden:comp} Fix $x \in \bpot{a}$. If $x \subseteq c$ for some $c \notin c \in a$, then $x= \complement{\complement{x}} \subseteq c$ so that $\complement{x} \in \bpot{a}$; if $\complement{x} \subseteq c$  for some $c \notin c \in a$, then $\complement{x} \in \bpot{a}$ straightforwardly.
\end{proof}\noindent
It follows that \boolevel{}s (see Definition \ref{def:bevel}) have several important closure properties:
\begin{lem}[\ECS]\label{lem:ecs:closed} 
	Every \boolevel is \blttrans, \bltpot, {complement-closed, and non-self-membered}.
\end{lem}
\begin{proof}
	Let $s$ be a \boolevel, i.e.\ $s = \bpot{h}$ for some \boohist $h$. So $s$ is \bltpot and complement-closed by Lemma \ref{lem:ecs:bwiden}. For \blttransity, fix {$a \notin a \in s = \bpot{h}$; so $a \subseteq c \notin c \in h$ for some $c$ by Lemma \ref{lem:ecs:bwiden}.\ref{bwiden:elem}; and $c = \bpot(c \cap h)$ as $h$ is a \boohist; so $a \subseteq \bpot(c \cap h) \subseteq \bpot{h} = s$. To see $s \notin s$, suppose $s\in s$ for reductio. Then $\complement{s} \notin \complement{s} \in s$ by \ref{blt:comp}, so $\complement{s}  \subseteq s$ by \blttransity, so $s = V$. Since $h \notin h$ by definition, and $h \in V = s = \bpot{h}$, by Lemma \ref{lem:ecs:bwiden}.\ref{bwiden:elem} there is some $c$ such that $h \subseteq c \notin c \in h$. Since $h$ is a \boohist, $c = \bpot{(h \cap c)} = \bpot{h} = V$, contradicting the fact that $c \notin c$.}
\end{proof}\noindent
From here, we can prove the well-ordering of the \boolevel{}s, by proving a sequence of results like those from Pt.\ref{pt:lt} \S\ref{s:1:LTwo}; I leave this to the reader:\footnote{For Lemma \ref{lem:ecs:levhist}, first note that if $h$ is a history and $c \in h$, then $c = \bpot{(c \cap h)} \subseteq \bpot{h} \notin \bpot{h}$ by Lemma \ref{lem:ecs:closed}, so $c\notin c$ by Lemma \ref{lem:ecs:lowhigh}. For Lemmas \ref{lem:ecs:acc}--\ref{lem:ecs:comparability}, reason about non-self-membered sets in the first instance, then deal with self-membered sets using \ref{blt:comp} and complement-closure.}
\begin{lem}[\ECS]\label{lem:ecs:regularity}
	If there is an $F$, {and all $F$s are non-self-membered and \bltpot}, then there is an $\in$-minimal $F$. Formally: $\forall F((\exists x F(x) \land (\forall x : F)(x \notin x \land x\text{ is \bltpot{}})) \lonlyif (\exists a : F)(\forall x : F)x \notin a)$
\end{lem}
\begin{lem}[\ECS]\label{lem:ecs:induction}
	If some \boolevel is $F$, then there is an $\in$-minimal \boolevel which is $F$. 
	Formally: $\forall F((\exists s : \bevpred)F(s) \lonlyif (\exists s: \bevpred)(F(s) \land (\forall r : \bevpred)(F(r) \lonlyif r \notin s)))$
\end{lem}
\begin{lem}[\ECS]\label{lem:ecs:levhist}
	Every member of a \boohist is a \boolevel.
\end{lem}
\begin{lem}[\ECS]\label{lem:ecs:acc} $s = \bpot{\Setabs{r \in s}{\bevpred(r)}}$, for any \boolevel $s$.
\end{lem}
\begin{lem}[\ECS]\label{lem:ecs:comparability}
	All \boolevel{}s are comparable, i.e.\ 
	$(\forall s: \bevpred)(\forall t : \bevpred)(s \in t \lor s = t \lor t \in s)$
\end{lem}\noindent
Combining Lemmas \ref{lem:ecs:induction} and \ref{lem:ecs:comparability}, \ECS proves that the \boolevel{}s are well-ordered by $\in$; this is Theorem \ref{thm:ecs:wo}. This licenses our use of the $\bevof{}$-operator (see Definition \ref{def:bevof}). Here are some simple results about that operator, which can be proved by tweaking the proof of Pt.\ref{pt:lt}  Lemma \ref{lem:lt:levof}:
\begin{lem}[\BLT]\label{lem:blt:bevof}
	For any sets $a, c$, and any \boolevel{}s $r, s$: 
	\setcounter{ncounts}{0}
	\begin{listn}
		\item\label{bevofexists} $\bevof{a}$ exists
		\item\label{bevofnotin} $a \notin \bevof{a}$
		\item\label{bevofquick} $r\subseteq s$ iff $s\notin r$
		\item\label{bevofidem} $s = \bevof{s}$
		\item\label{bevofsubs} if $c \subseteq a \notin a$ or $a \in a \subseteq c$, then $\bevof{c} \subseteq \bevof{a}$
		\item\label{bevofin} if $c \in a \notin a$ or $c \notin a \in a$, then $\bevof{c} \in \bevof{a}$
	\end{listn}
\end{lem}\noindent
Moreover, we can now show that sets are closed under arbitrary pairwise intersection:
\begin{lem}[\BLT]\label{lem:bst:capcup}
	For any sets $a$ and $c$, the set $a \cap c = \Setabs{x}{x \in a \land x \in c}$ exists.
\end{lem}
\begin{proof}
	First suppose that either $a \notin a$ or $c\notin c$ (or both); without loss of generality, suppose $a \notin a$; now $a \cap c = \Setabs{x \in a}{x \in c}$ exists by \ref{blt:sep}. Next suppose that both $a \in a$ and $c \in c$. So both $\complement{a} \notin \complement{a}$ and $\complement{c} \notin \complement{c}$ by \ref{blt:comp}. Let $s$ be the maximum of $\bevof{\complement{a}}$ and $\bevof{\complement{c}}$. Since $s$ is \bltpot{}, both $\complement{a} \subseteq s$ and $\complement{c} \subseteq s$, so $\complement{a} \cup \complement{c} = \Setabs{x \in s}{x \in \complement{a} \lor x \in \complement{c}}$ exists by \ref{blt:sep}. Now $a \cap c = \complement{\complement{a} \cup \complement{c}}$ exists by \ref{blt:comp}.
\end{proof}\noindent
This immediately entails that the sets form a boolean algebra, which is Theorem \ref{thm:blt:boolean}. Our next result shows that the universe is contra-automorphic:\footnote{\label{fn:blt:recursion}Theorem \ref{thm:ecs:wo} licenses recursive definitions. We can regard as defining second-order entities. If we are using second-order logic, such definitions yield a second-order entity. If we are using first-order logic, then (as usual) we define a term by considering a strictly increasing sequence of first-order 'bounded approximations' (specifying the behaviour of the term over the last few \boolevel{}s manually, if there is a last \boolevel).}
\begin{namedthm}[\ref{thm:blt:antimorphism} \normalfont{(\BLT)}]
	$\forall a \forall b(a \in b \liff \gameneg a \notin  \gameneg b)$
\end{namedthm}
\begin{proof}
	Recall that negative is given as in Definition \ref{def:gameneg} by 
	\begin{align*}
		\gameneg a &\coloneq 
		\complement{\Setabs{\gameneg x}{x \in a}} \text{, if }a \notin a &
		\gameneg a&\coloneq 
		\Setabs{\gameneg x}{x \notin a} \text{, if }a \in a
	\end{align*}
	Fix a \boolevel $s$ and for induction suppose that, for any $x, y \in s$:
	\begin{listn-0} 
		\item\label{pi:rank} $\gameneg x$ is well-defined and $\bevof{x} = \bevof{(\gameneg x)}$; and
		\item\label{pi:bijective} $x = y$ iff $\gameneg x = \gameneg y$.
	\end{listn-0}
	It suffices to show that both properties hold of $a, b$ when $\bevof{a} = \bevof{b} = s$.
	
	\emph{Concerning (\ref{pi:rank}).} Suppose $a \notin a$. If $x \in a$, then $\bevof{(\gameneg x)} = \bevof{x} \in \bevof{a}$ by induction assumption \eqref{pi:rank} and Lemma \ref{lem:blt:bevof}.\ref{bevofin}. Using \ref{blt:sep}, let $c \notin c = \Setabs{v \in \bevof{a}}{(\exists x \in a)v = \gameneg x} = \Setabs{\gameneg x}{x \in a}$. Moreover, $\bevof{c} = \bevof{a}$, by the well-ordering of \boolevel{}s and since $\bevof{(\gameneg x)} = \bevof{x} \in \bevof{a}$ for all $x \in a$. Now $\complement{c} \in \complement{c} = \gameneg a$ by \ref{blt:comp}; so $\bevof{a} = \bevof{c} = \bevof{\complement{c}} = \bevof{(\gameneg a)}$. The case when $a \in a$ is similar, defining $c \notin c = \Setabs{v \in \bevof{a}}{(\exists x \notin a)v = \gameneg x} = \Setabs{\gameneg x}{x \notin a} = {\gameneg a}$. 
	
	\emph{Concerning (\ref{pi:bijective}).} If $a \in a \liff b \in b$, then $a = b$ iff $\gameneg a = \gameneg b$ by induction assumption \eqref{pi:bijective}. Without loss of generality, suppose that $a \in a$ and $b \notin b$; in establishing \eqref{pi:rank}, we found that $\gameneg a \notin \gameneg a$ and $\gameneg b \in \gameneg b$; so $a \neq b$ and $\gameneg a \neq \gameneg b$.
\end{proof}\noindent
I ended \S\ref{s:3:characteristic} by stating some simple facts about extensions of \BLT. I will prove the distinctively boolean facts, leaving the remainder to the reader:
\begin{namedprop}[\ref{prop:pstl:mltint}\normalfont{, fragment}]
	\textcolor{white}{.}
	\begin{listn-0}
		\item[\eqref{blt:union}] \BLT proves Union, i.e.\ $\forall a(\bigcup a\text{ exists})$
		\item[\eqref{blt:nopowersets}] \BLT contradicts Powersets, i.e.\ it proves $\exists a \lnot \exists  b \forall x(x \in b \liff x \subseteq a)$
		\item[\eqref{blt:highfoundation}] \BLT proves Foundation-restricted-to-high-sets, i.e.\ $(\forall a \in a)(\exists x \in a)a \cap x = \emptyset$. 
		\item[\eqref{blt:nofoundation}] \BLTplus contradicts unrestricted Foundation, i.e.\ it proves $(\exists a \neq \emptyset)(\forall x \in a)a \cap x \neq \emptyset$. 
	\end{listn-0}
\end{namedprop}
\begin{proof}
	\emphref{blt:union} If $a \in a$, then $\bigcup a = \complement{\Setabs{x \in \complement{a}}{(\forall y \in a)x \notin y}}$, which exists by \ref{blt:sep} and \ref{blt:comp}. If $a \notin a$, then using \ref{blt:sep} let $a_0 = \Setabs{x \in a}{x \notin x}$ and let $a_1 = \Setabs{x \in a}{x \in x}$. I will show that $\bigcup a_0$ and $\bigcup a_1$ exist, so that, using \ref{blt:comp} and Lemma \ref{lem:bst:capcup}: 
	$$\bigcup a = \bigcup a_0 \cup \bigcup a_1   = \complement{\complement{\bigcup a_0} \cap \complement{\bigcup a_1}}$$
	Clearly $\bigcup a_0$ exists by \ref{blt:sep} on $\bevof{a}$. If $a_1 = \emptyset$ then $\bigcup a_1 = \emptyset$; otherwise, $\bigcup a_1 = \complement{\bigcap \Setabs{\complement{x}}{x \in a_1}}$, which exists by \ref{blt:comp} and \ref{blt:sep} on $\bevof{a}$.
	
	\emphref{blt:nopowersets} If there is only one \boolevel, then the only sets are $\emptyset$ and $V = \{\emptyset, V\}$, so that $\powerset \emptyset = \{\emptyset\}$ does not exist. Otherwise, we find $\complement{\{\emptyset\}}$ at the second \boolevel, and if $\powerset\complement{\{\emptyset\}}$ existed it would be $\Setabs{x}{\emptyset \notin x}$. So suppose for reductio that $a = \Setabs{x}{\emptyset \notin x}$. Then $\emptyset \notin \emptyset$, so $\emptyset \in a$, so $a \notin a$. Now $\complement{a} \in \complement{a} = \Setabs{x}{\emptyset \in x}$ by \ref{blt:comp}, so that $\emptyset \in \complement{a}$, contradicting that $\emptyset \in a$. 
		
	\emphref{blt:highfoundation} If $a \in a$ then $\complement{a} \in a$ by \ref{blt:comp}, and $a \cap \complement{a} = \emptyset$.
	
	\emphref{blt:nofoundation} We find $\{V\}$ at the second \boolevel, and $\{V\} \cap V \neq \emptyset$. 
\end{proof}

\section{The set-theoretic equivalence of \BST and \BLT}\label{s:thm:BLTBST}
I now want to prove Theorem \ref{thm:BLTBST}, which states that \BLT and \BST say exactly the same things about sets. (This mirrors Pt.\ref{pt:lt} \S\ref{s:1:LTST}.)

To show that \BST says no more about sets than \BLT does, I define a translation $* : \BST \functionto \BLT$, whose non-trivial actions are as follows:\footnote{So the other clauses are: 
			$(\lnot \phi)^* \coloneq \lnot \phi^*$; 
			$(\phi \land \psi)^* \coloneq (\phi^* \land \psi^*)$; 
			$(\forall x\phi)^* \coloneq \forall x \phi^*$; 
			$(\forall F \phi)^* \coloneq \forall F \phi^*$; and $\alpha^* \coloneq \alpha$ for all atomic formulas $\alpha$ which are not of the forms mentioned in the main text.}
	\begin{align*}
		\lowpred(x)&:= x \notin x & \highpred(x) &\coloneq  x \in x\\
		(\stage{s} < \stage{t})^* &\coloneq \stage{s} \in \stage{t} & 
		(x \foundat \stage{s})^* &\coloneq  (x \subseteq \stage{s} \lor \complement{x} \subseteq \stage{s}) &
		(\forall \stage{s} \phi)^* &\coloneq (\forall \stage{s} : \bevpred)(\phi^*)
	\end{align*}\noindent
After translation, I treat all first-order variables as being of the same sort. Fairly trivially, for any \BLT-sentence $\phi$, if $\BST  \proves \phi$ then $\BST^* \proves \phi$. The left-to-right half of Theorem \ref{thm:BLTBST} now follows as $*$ is an interpretation:
\begin{lem}[\BLT]\label{lem:blt:bsttrans} $\BST^*$ holds.
\end{lem}
\begin{proof}
	\ref{bst:ext}$^*$ is \ref{ext}. 
	\ref{bst:ord}$^*$ holds by Lemma \ref{lem:ecs:closed};  \ref{bst:stage}$^*$ holds by {\ref{blt:strat} and \ref{blt:comp}}; and \ref{bst:cases}$^*$ is trivial. Next, by Lemmas \ref{lem:ecs:closed} and \ref{lem:ecs:acc}, we can simplify $(x \foundby \stage{s})^*$ to $x \in \stage{s}$. So, using Lemmas \ref{lem:ecs:lowhigh} and \ref{lem:ecs:closed}, we can simplify \ref{bst:pri0}$^*$ thus:
	\begin{align*}
		(\forall s \in \bevpred)&
		(\forall a\notin a)((a \subseteq s \lor \complement{a} \subseteq s) \lonlyif (\forall x \in a)x \in s) \\
		\text{i.e.\ }(\forall s \in \bevpred)&(\forall a \subseteq s )(\forall x \in a) x \in s 
	\end{align*}
	which is trivial; then \ref{bst:pri1}$^*$ holds similarly, by \ref{blt:comp}. A similar simplification allows us to obtain \ref{bst:find0}$^*$ via \ref{blt:sep}; then \ref{bst:find1}$^*$ holds similarly, by \ref{blt:comp}.\footnote{Note that the $*$-translation of any \BST-Comprehension instance is a \BLT-Comprehension instance.}
\end{proof}\noindent
To obtain the right-to-left half of Theorem \ref{thm:BLTBST}, I will work in \BST. I start by defining \emph{slices}, which will go proxy for stages, and will turn out to be \boolevel{}s, and then stating a few elementary results (for proofs, tweak those of Pt.\ref{pt:lt} \S\ref{s:1:LTST}):
\begin{define}[\BST]
	For each $\stage{s}$, let $\bslice{s} = \Setabs{x}{x \foundby \stage{s}}$. Say that $a$ is a slice iff $a = \bslice{s}$ for some stage $\stage{s}$.
\end{define}
\begin{lem}[\BST]\label{lem:bst:sep:l}$\forall F (\forall a : \lowpred)(\exists b : \lowpred)\forall x(x \in b \liff (F(x) \land x \in a))$
\end{lem}
\begin{lem}[\BST]\label{lem:bst:conversepri}
	$\forall \stage{s}(\forall a : \lowpred)(a \foundat \stage{s} \liff (\forall x \in a)x \foundby \stage{s})$
\end{lem}
\begin{lem}[\BST]\label{lem:bst:basicslice} For any \stage{s}:
	\begin{listn-0}
		\item\label{bslice:low} $\bslice{s}$ exists and is low
		\item\label{bslice:foundup} {$\forall \stage{r}(\forall a : \lowpred)(a \foundat \stage{r}  \leq \stage{s} \lonlyif a \foundat \stage{s})$}	
		\item\label{bslice:foundat} $(\forall a : \lowpred)(a \subseteq \bslice{s} \liff a \foundat \stage{s})$
	\end{listn-0}
\end{lem}
\noindent
We must now part company slightly with the strategy of Pt.\ref{pt:lt} \S\ref{s:1:LTST}, to handle low and high sets, and their relation to (non-)self-membership:
\begin{lem}[\BST]\label{lem:bst:induction}
	If some slice is $F$, then there is an $\in$-minimal slice which is $F$. 
\end{lem}
\begin{proof}
	Every slice is low, by Lemma \ref{lem:bst:basicslice}.\ref{bslice:low}. Subsets of low sets are low, by a result like Lemma \ref{lem:ecs:lowhigh}. From this, and Lemma \ref{lem:bst:basicslice}, it follows that $\forall\slice{s}\forall x((\exists c : \lowpred)x \subseteq c \in \slice{s} \lonlyif x \in \slice{s})$. The result now follows, reasoning as in Pt.\ref{pt:lt} Lemma \ref{lem:es:regularity}. 
\end{proof}
\begin{lem}[\BST]\label{lem:bst:lowin}$a$ is low iff $a \notin a$; and $a$ is high iff $a \in a$. 
\end{lem}
\begin{proof}
	Suppose for reductio that $a \in a$ is low. Using \ref{bst:stage} and Lemma \ref{lem:bst:induction}, let $\bslice{s}$ be an $\in$-minimal slice such that $\exists \stage{t}(a \foundat \stage{t} \land \bslice{t}= \bslice{s})$; let $\stage{t}$ witness this. Since $a \in a \foundat \stage{t}$ and $a$ is low, $a \foundat \stage{r} < \stage{t}$ for some $\stage{r}$ by \ref{bst:pri0}; so $\bslice{r} \in \bslice{t} = \bslice{s}$ by Lemma \ref{lem:bst:basicslice}, contradicting $\bslice{s}$'s minimality. Discharging the reductio: if $a$ is low, then $a \notin a$. Similarly: if $a$ is high, then $a \in a$. The biconditionals follow by \ref{bst:cases}.
\end{proof}
\begin{lem}[\BST]\label{lem:bst:comp} $\complement{a}$ exists; and $a \notin a \liff \complement{a} \in \complement{a}$; and $\forall \stage{s}(a \foundat \stage{s} \liff \complement{a} \foundat \stage{s})$.
\end{lem}
\begin{proof}
	Using \ref{bst:stage}, let $a \foundat \stage{s}$. If $a \notin a$, {then $a$ is low by Lemma \ref{lem:bst:lowin}, so} $(\forall x \in a)x \foundby \stage{s}$ by {\ref{bst:pri0}}, so that by {\ref{bst:find1}} and \ref{ext} $\Setabs{x}{x \notin a} = \complement{a} \foundat \stage{s}$ exists and is high, i.e.\ $\complement{a} \in \complement{a}$ {by Lemma \ref{lem:bst:lowin}}. If $a \in a$, reason similarly {using \ref{bst:pri1} and \ref{bst:find0}}. 
\end{proof}\noindent
Note that $\BST \proves \ECS$ by Lemmas \ref{lem:bst:sep:l}, \ref{lem:bst:lowin}, and \ref{lem:bst:comp}. So Lemmas \ref{lem:ecs:lowhigh}--\ref{lem:ecs:comparability} hold verbatim within \BST. We can now complete our reasoning about slices, by resuming the proof-strategy of Pt.\ref{pt:lt} \S\ref{s:1:LTST}; at this point, I leave the remaining details to the reader:
\begin{lem}[\BST]\label{lem:bst:niceslice} $\bslice{s} \notin \bslice{s}$; and $\bslice{s}$ is \blttrans; and $\bslice{s} = \bpot{\Setabs{\bslice{r}}{\bslice{r} \in \bslice{s}}}$.
\end{lem}
\begin{lem}[\BST]\label{lem:bst:comparability}
	All slices are comparable, i.e.\ 
	$\forall \bslice{s} \forall \bslice{t} (\bslice{s} \in \bslice{t} \lor \bslice{s} = \bslice{t} \lor \bslice{t} \in \bslice{s})$.
\end{lem}
\begin{lem}[\BST]\label{lem:bst:levelsslices} $s$ is a \boolevel iff $s$ is a slice.
\end{lem}\noindent
It follows that \BST proves \ref{blt:strat}, delivering Theorem \ref{thm:BLTBST}.

\section{Helow sets}\label{s:blt:helow:appendix}
In this appendix I prove Theorem \ref{thm:blt:helowinterpret}, which shows how to recover ordinary, uncomplemented hierarchies via helow sets (see Definition \ref{def:helow}). For readability, I refer to non-self-membered sets as \emph{low}, and self-membered sets as \emph{high} (cf.\ Lemma \ref{lem:bst:lowin}). Note that every helow set is low, since all its members are low (i.e.\ non-self-membered). Now: 
\begin{define}[\BLT]\label{def:helowmark}
	If $a$ is low, let $a\helowmark \coloneq \Setabs{x \in a}{x\text{ is helow}}$; by \ref{blt:sep}, $a\helowmark$ exists and is low.
\end{define}
\begin{lem}[\BLT]\label{lem:helow:equiv} $a$ is helow iff every member of $a$ is helow.
\end{lem}
\begin{proof}
	\emph{Left-to-right.} Where $c$ witnesses that $a$ is helow, if $x \in a$, then $x \in c$ and hence $x \subseteq c$, so $c$ also witnesses that $x$ is helow. \emph{Right-to-left}. Let every member of $a$ be helow. Every member of $a$ is low, so $a$ itself is low; hence $a \subseteq (\bevof{a})\helowmark$. Now $(\bevof{a})\helowmark$ witnesses that $a$ is helow: if $x \in c \in (\bevof{a})\helowmark$ then $c$ is helow so $x$ is helow (by left-to-right), so $x \in  (\bevof{a})\helowmark$ as $\bevof{a}$ is \blttrans.
\end{proof}\noindent
I can now begin to show that $\helowobj: \LT \functionto \BLT$, which simply restricts all quantifiers to helow sets (see \S\ref{s:3:helow}), is an interpretation of \LT:
\begin{lem}[\BLT]\label{lem:blt:helow:es} 
Both \ref{lt:ext}$\helowint$ and \ref{sep}$\helowint$ hold.
\end{lem}
\begin{proof}
	For \ref{lt:ext}$\helowint$, fix helow $a$ and $b$ and suppose that $(\forall x : \helow)(x \in a \liff x \in b)$; then $\forall x(x \in a \liff x \in b)$ by Lemma \ref{lem:helow:equiv}, so $a = b$ by \ref{ext}. Similarly, repeated use of Lemma \ref{lem:helow:equiv} shows that \ref{sep}$\helowint$ follows from \ref{blt:sep}.
\end{proof}\noindent
The next task is to connect \boolevel{}s with levels$\helowint$. (See Pt.\ref{pt:lt} Definitions \ref{def:pot}--\ref{def:potent} for the definitions of \emph{potent}, $\pot$, $\histpred$ and $\levpred$.)
\begin{lem}[\BLT]\label{lem:helow:bevlevtrick} For any \boolevel{}s $r, s$:
	 \begin{listn-0}
	 	\item\label{bevlev:helevtrans}	$s\helowmark$ is helow, potent and transitive
		\item\label{bevlev:sin} 
		$r \in s $ iff $r\helowmark \in s\helowmark$		
		\item\label{bevlev:nochange} 
		$s = \bevof{(s\helowmark)}$ 
		\item\label{bevlev:bevpot} 
		$s\helowmark = \pot{h} = \pot\helowint(h)$, where $h = \Setabs{r\helowmark \in s\helowmark}{\bevpred(r)}$.
		\item\label{bevlev:levels} 
		$s\helowmark$ is a level$\helowint$
	\end{listn-0}
\end{lem}
\begin{proof}
	\emphref{bevlev:helevtrans} By Lemma \ref{lem:helow:equiv}, $s\helowmark$ is helow; then $s\helowmark$ is potent and transitive as $s$ is \bltpot{} and \blttrans{}.

	\emphref{bevlev:sin} \emph{Left-to-right}. By \eqref{bevlev:helevtrans}. \emph{Right-to-left}. Let $r\helowmark \in s\helowmark$. So $r \neq s$, since $r\helowmark \notin r\helowmark$. Similarly, $s\helowmark \notin r\helowmark$, since $s\helowmark$ is transitive; so $s \notin r$ by \emph{left-to-right}. So $r \in s$, by Lemma \ref{lem:ecs:comparability}.
	
	\emphref{bevlev:nochange} Induction on \boolevel{}s, using \eqref{bevlev:sin}.
	
	\emphref{bevlev:bevpot} By \eqref{bevlev:helevtrans} and Lemma \ref{lem:helow:equiv}, $h$ is helow. If $a \in \pot{h}$, then $a \in s\helowmark$ as $s\helowmark$ is potent by \eqref{bevlev:helevtrans}. Conversely, if $a \in s\helowmark$, then $a \subseteq r \in s$ for some \boolevel $r$ by Lemma \ref{lem:ecs:acc}, and $a \subseteq r\helowmark \in s\helowmark$ by \eqref{bevlev:sin} and Lemma \ref{lem:helow:equiv}, so $a \in \pot{h}$. So $s\helowmark = \pot{h}$. Repeated use of Lemma \ref{lem:helow:equiv}, as in Lemma \ref{lem:blt:helow:es}, now yields that $\pot{h} = \pot\helowint(h)$.
	
	\emphref{bevlev:levels} With $h$ as in \eqref{bevlev:bevpot}, since $s = \pot\helowint(h)$ it suffices to show that $\histpred\helowint(h)$. If $r\helowmark \in h$, then $r\helowmark \cap h = \Setabs{q\helowmark \in r\helowmark}{\bevpred(q)}$, by \eqref{bevlev:helevtrans}; so $r\helowmark = \pot\helowint{(r \helowmark \cap h)}$ by \eqref{bevlev:bevpot}.
\end{proof}
\begin{lem}[\BLT]\label{lem:blt:levhelow} The levels$\helowint$ are the \boolevel{}s$\helowmark$, i.e.: $\levpred\helowint(a)$ iff $(\exists s : \bevpred)a = s\helowmark$.
\end{lem}
\begin{proof}
	By Lemma \ref{lem:helow:bevlevtrick}, if $s$ is a \boolevel then both $\levpred\helowint(s\helowmark)$ and $\bevof(s\helowmark) = s$. To complete the proof, it suffices to note that if $p$ and $q$ are distinct levels$\helowint$, then $\bevof{p} \neq \bevof{q}$; this follows from  Lemma \ref{lem:blt:bevof}.\ref{bevofin} and the fact that the levels$\helowint$ are well-ordered by $\in$. (The well-ordering of levels$\helowint$ is Pt.\ref{pt:lt} Theorem \ref{thm:es:wo}$\helowint$, which holds via Lemma \ref{lem:blt:helow:es}.)
\end{proof}
\begin{cor}[\BLT]\label{lem:blt:helowinterpret}$\clearme{\ref{lt:strat}}\helowint$ holds; \ref{blt:cre} proves $\clearme{\ref{lt:cre}}\helowint$; \ref{blt:inf} proves $\clearme{\ref{lt:inf}}\helowint$; and \ref{blt:rep} proves $\clearme{\ref{lt:rep}}\helowint$.
\end{cor}\noindent
Recalling that $\LT + \clearme{\ref{lt:inf}} + \clearme{\ref{lt:rep}}$ is equivalent to \ZF (see \S\ref{s:3:helow}), Lemmas \ref{lem:blt:helow:es} and \ref{lem:blt:helowinterpret} yield Theorem \ref{thm:blt:helowinterpret}.

\section{Definitional equivalence}\label{s:blt:de:appendix}
In this appendix, I prove the definitional equivalence discussed in \S\ref{s:3:de}.\footnote{Recall: both \LT and \BLT (and their extensions) are formulated as second-order theories. I continue to frame my discussion in second-order terms in this appendix. However, the theories can easily be reformulated as first-order formulations, and the definitional equivalences hold for these first-orderisations (only the quasi-categoricity results of \S\ref{s:3:quasicat} require second-order resources).}

\subsection{Interpreting \BLTzf in \ZF}\label{s:interpret:bltinterpretszf}
I first define an interpretation, $\cointer$, to simulate (extensions of) \BLT within (extensions of) \LT. The key idea is to use $\emptyset$ as a flag to indicate whether to treat a set as low or high. To allow $\emptyset$ to play this role, I define a bijection $\coinj : V \functionto V \setminus \{\emptyset\}$:\footnote{Many thanks to Randall Holmes for discussion of this construction (and other constructions); the proof in this section is much more self-contained than it would have been, had it not been for his input. Thanks also to Thomas Forster, for encouraging me to consider the question of definitional equivalence. The proof-strategy is similar to \textcite{Lowe:STWWU}.} 
\begin{align*}
	\coinj(a) &\coloneq 
	\begin{cases}
		\{a\}&\text{if }a\text{ is a Zermelo number}\\
		a&\text{otherwise}
	\end{cases}
\end{align*}
where the Zermelo numbers are $0 = \emptyset$ and $n+1 = \{n\}$. I then interpret membership thus:
\begin{align*}
	x  \coin a &\text{ iff }
	(\coinj(x) \in a \liff \emptyset \notin a) 
\end{align*}
Since $\coinj(a) \notin a$ for all $a$, it follows that $a \conotin a$ iff $\emptyset \notin a$ (i.e.\ $a$ is treated as low), and $a \coin a$ iff $\emptyset \in a$ (i.e.\ $a$ is treated as high). I will now prove a sequence of results which establish that $\cointer$ is an interpretation of \BLT. The first few are straightforward:
\begin{lem}[\LTplus]\label{lem:co:subcosub}Where $a \subseteq\cotrans b$ abbreviates $(\forall x \coin a)x \coin b$:
	\begin{listn-0}
		\item\label{subcosub:low} If $\emptyset \notin a$ and $\emptyset \notin b$, then: $a \subseteq b$ iff $a \subseteq\cotrans b$
		\item\label{subcosub:high} If $\emptyset \in a$ and $\emptyset \in b$, then: $a \supseteq b$ iff $a \subseteq\cotrans b$.
	\end{listn-0}
\end{lem}
\begin{proof}
	\emphref{subcosub:low} Since $\coinj$ is a bijection $V \functionto V \setminus \{\emptyset\}$, $a \subseteq b$ iff $\forall x(\coinj(x) \in a \lonlyif \coinj(x) \in b)$ iff $a \subseteq\cotrans b$. 
	
	\emphref{subcosub:high} Similarly, $a \supseteq b$ iff $\forall x(\coinj(x) \notin a \lonlyif \coinj(x) \notin b)$ iff $a \subseteq\cotrans b$.
\end{proof}
\begin{lem}[\LTplus]\label{lem:ltco:ext}
	\ref{ext}$\cotrans$ holds. 
\end{lem}
\begin{proof}
	Suppose $\forall x(x \coin a \liff x \coin b)$. If $a \conotin a$ but $b \coin b$, then $\forall x(\coinj(x) \in a \liff \coinj(x) \notin b)$, so that $a \cup b = V$, which is impossible. Generalising, $a \coin a$ iff $b\coin b$. Now apply Lemma \ref{lem:co:subcosub}.
\end{proof}
\begin{lem}[\LTplus]\label{lem:ltco:sep}
	\ref{blt:sep}$\cotrans$ holds.
\end{lem}
\begin{proof}
	Fix $F$ and $a \conotin a$, i.e.\ $\emptyset \notin a$. Using \ref{sep}, let $b = \Setabs{\coinj(x) \in a}{F(x)}$. Since $\emptyset \notin b$ we have $\forall x(x \coin b \liff (F(x) \land x \coin a))$.
\end{proof}\noindent
The interpretation of complementation is obvious: $\cocomp{a}= a \cup \{\emptyset\}$ if $a \conotin a$, and $\cocomp{a} = a \setminus \{\emptyset\}$ if $a \coin a$. The next result follows trivially:
\begin{lem}[\LTplus]\label{lem:ltco:comp} $\forall a\forall x(x \coin a \liff x \conotin \cocomp{a})$, and \ref{blt:comp}$\cotrans$ holds.
\end{lem}\noindent
The only intricate part of this interpretation concerns the treatment of \boolevel{}s. Within \LTplus, we can define the von Neumann ordinals, and recursively define the following:
	$$W_\gamma = \Setabs{\coinj(x)}{(\exists \beta < \gamma)x \subseteq W_\beta \cup \{\emptyset\}}$$ 
Now \LTplus proves that $W_\gamma$ exists for each $\gamma$, and that these are the \boolevel{}s$\cotrans$:
\begin{lem}[\LTplus]\label{lem:ltco:newlev} $\bevpred\cotrans(s)$ iff $s = W_\gamma$ for some $\gamma$.
\end{lem}
\begin{proof}
	Lemmas \ref{lem:ltco:ext}--\ref{lem:ltco:comp} show that \LTplus proves $\ECS\cotrans$. Hence \LTplus proves Theorem \ref{thm:ecs:wo}$\cotrans$, i.e.\ that the \boolevel{}s$\cotrans$ are well-ordered by $\coin$. For induction on $\gamma$, suppose that if $\beta < \gamma$ then $W_\beta$ is the $\beta^{\text{th}}$ \boolevel{}$\cotrans$. Let $s$ be the $\gamma^{\text{th}}$ \boolevel{}$\cotrans$. By Lemma \ref{lem:co:subcosub}:
	\begin{align*}
		W_\gamma &= \Setabs{\coinj(x)}{(\exists \beta < \gamma)x \subseteq W_\beta \cup \{\emptyset\}}\\
		&= \Setabs{\coinj(x)}{(\exists \beta < \gamma)(x \subseteq\cotrans W_\beta \lor \complement{x}\cotrans \subseteq\cotrans W_\beta)}\\
		&= \Setabs{\coinj(x)}{(\exists W_\beta \conotin W_\beta \coin s)(x \subseteq\cotrans W_\beta \lor \complement{x}\cotrans \subseteq\cotrans W_\beta)}\\
		&=(\bpot\Setabs{w \in s}{\bevpred(w)})\cotrans
	\end{align*}
	So $W_\gamma = s$ by Lemma \ref{lem:ecs:acc}$\cotrans$. By induction, the \boolevel{}s$\cotrans$ are the $W_\gamma$s.
\end{proof}\noindent
I can now prove the crucial proposition: 
\begin{lem}[\LTplus]\label{lem:ltco:rank}\ref{blt:strat}$\cotrans$ holds.
\end{lem}
\begin{proof}
	By Lemma \ref{lem:ltco:newlev}, it suffices to show that $(\forall a \conotin a)\exists \gamma\, \ a \subseteq\cotrans W_\gamma$. Since the levels are well-ordered by $\in$ (Pt.\ref{pt:lt} Theorem \ref{thm:es:wo}), we can write $V_\gamma$ for the $\gamma^{\text{th}}$ level. I claim: if $a \conotin a \subseteq V_\gamma$, then $a \subseteq W_\gamma$. For induction, suppose this holds for all ordinals $\beta < \gamma$. Fix $a \conotin a \subseteq V_\gamma$. If $\gamma = 0$, then $a = \emptyset \subseteq\cotrans W_0 = \emptyset$. Otherwise, fix $x \coin a$, i.e.\ $\coinj(x) \in a \subseteq V_\gamma$; now $x \subseteq V_\beta$ for some $\beta < \gamma$, by Pt.\ref{pt:lt} Lemma \ref{lem:lt:levof}, so that $x \subseteq W_\beta \cup \{\emptyset\}$ by the induction hypothesis; so $\coinj(x) \in W_\gamma$, i.e.\ $x \coin W_\gamma$. Generalising, $a \cosubseteq W_\gamma$.
\end{proof}
\begin{lem}\label{lem:lt:cotransworks}$\LTplus \proves \BLTplus\cotrans$ and $\ZF \proves \BLTzf\cotrans$.
\end{lem}
\begin{proof}
	Lemmas \ref{lem:ltco:ext}--\ref{lem:ltco:rank} establish that $\LTplus \proves \BLT\cotrans$. And $\LTplus \proves \clearme{\ref{blt:cre}}\cotrans$, using \ref{lt:cre} and our explicitly defined \boolevel{}s$\cotrans$, the $W_\gamma$s. Evidently, \ref{lt:inf} yields \ref{blt:inf}$\cotrans$. For \ref{blt:rep}$\cotrans$, fix $P$ and $a \conotin a$; by \ref{lt:rep}, the set $c = \Setabs{\coinj(P(x))}{\coinj(x) \in a}$ exists; by construction, $\emptyset \notin c$ and $(\forall x \coin a)P(x) \coin c$. The result follows, since \ZF is equivalent to $\LT + \clearme{\ref{lt:inf}} + \clearme{\ref{lt:rep}}$ (see \S\ref{s:3:helow}).
\end{proof}

\subsection{Interpreting \ZF in \BLTzf}
I now switch to working in \BLTplus. Using $\coinjoc$---i.e.\ using verbatim the same definitions of `Zermelo number' and of $\coinj$ in \BLTplus as we used in \LTplus---consider this function:
\begin{align*}
	\helowbij(a) &= 
	\begin{cases}
		\Setabs{\coinjoc(\helowbij(x))}{x \in a}&\text{if }a\notin a \\
		\Setabs{\coinjoc(\helowbij(x))}{x \notin a} \cup \{\emptyset\}&\text{if }a \in a
	\end{cases}
\end{align*}\noindent
I will prove that $\helowbij$ is a bijection $V \functionto \helow$. I then define a translation, $\ocinter$, by stipulating: 
	$$x \ocin a\text{ iff }\helowbij(x) \in \helowbij(a)$$
It will follow that $\ocinter$ is an interpretation of \LTplus in \BLTplus. 
\begin{lem}[\BLTplus]\label{lem:oc:ocinjinj} If $\helowbij(a) = \helowbij(b)$, then $a = b$. 	
\end{lem}
\begin{proof}
	Let $\helowbij(a) = \helowbij(b)$, so that $a \notin a \liff b \notin b$. For induction, suppose that $\helowbij(x) = \helowbij(y) \lonlyif x = y$ for all $x, y$ with $\bevof{x}, \bevof{y} \in \bevof{a} \cup \bevof{b}$. If $a \notin a$ and $b \notin b$, then $\Setabs{\coinjoc(\helowbij(x))}{x \in a}  = \Setabs{\coinjoc(\helowbij(x))}{x \in b}$, so that $a = b$ by the induction hypothesis and the injectivity of $\coinjoc$. The case when $a \in a$ is similar. 
\end{proof}
\begin{lem}[\BLTplus]\label{lem:oc:ocallhelow} $\helowbij(a)$ is helow, for any $a$.
\end{lem}
\begin{proof}
	For induction, suppose that $\helowbij(x)$ is helow for all $x$ with $\bevof{x} \in \bevof{a}$. Suppose $a \notin a$; since $\coinjoc(\helowbij(x))$ is helow iff $\helowbij(x)$ is helow, every member of $\helowbij(a)$ is helow; so $\helowbij(a)$ is helow by Lemma \ref{lem:helow:equiv}. The case when $a \in a$ is similar. 
\end{proof}
\begin{lem}[\BLTplus]\label{lem:oc:ocsurj}
	If $a$ is helow, then $a = \helowbij(c)$ for some $c$. 
\end{lem}
\begin{proof}             
	By Lemma \ref{lem:oc:ocinjinj}, $\helowbij^{-1}$ is functional. For induction, suppose that for all helow $z \in \bevof{a}$, we have that $\helowbij^{-1}(z)$ is defined and $\bevof{(\helowbij^{-1}(z))} \subseteq \bevof{z}$.
	
	If $\emptyset \notin a$, let $c \notin c =  \Setabs{\helowbij^{-1}(\coinjoc^{-1}(x)) \in \bevof{a}}{x \in a}$ using \ref{blt:sep}. Fix $x \in a$; then $\coinjoc^{-1}(x) \in \bevof{a}$ and $\coinjoc^{-1}(x)$ is helow, recalling that $a$ is helow and using Lemma \ref{lem:helow:equiv}). Now $\bevof(\helowbij^{-1}(\coinjoc^{-1}(x))) \subseteq \bevof(\coinjoc^{-1}(x)) \in \bevof{a}$ by the induction hypothesis, i.e.\ $\helowbij^{-1}(\coinjoc^{-1}(x)) \in \bevof{a}$. So $c = \Setabs{\helowbij^{-1}(\coinjoc^{-1}(x))}{x \in a}$, so that $a = \helowbij(c)$ and $\bevof{c} \subseteq \bevof{a}$.
	
	If $\emptyset \in a$, then instead let $c = \Setabs{\helowbij^{-1}(\coinjoc^{-1}(x))}{\emptyset \neq x \in a}$; now $a = \helowbij(\complement{c})$. 
\end{proof}
\begin{lem}\label{lem:blt:octransworks}$\BLTplus \vdash \LTplus\octrans$ and $\BLTzf \proves \ZF\octrans\textsl{}$.
\end{lem}
\begin{proof}
	By Lemmas \ref{lem:oc:ocinjinj}--\ref{lem:oc:ocsurj},  $\helowbij : V\functionto \helow$ is a bijection; now use Theorem \ref{thm:blt:helowinterpret}. 
\end{proof}
\subsection{The interpretations are inverse}\label{s:interpret:inverse}
It remains to show that $\cointer$ and $\ocinter$ are mutually inverse, in the sense required for definitional equivalence.\footnote{Via \textcite[Corollary 5.5]{FriedmanVisser:WBIS}, to establish Theorem \ref{thm:synonymy} we could instead verify that $\cointer$ and $\helowobj$ (from \S\ref{s:blt:helow:appendix}) are bi-interpretations.} The key lies in their treatments of the Zermelo numbers. Working informally, let $\zero_n$ be the $n^{\text{th}}$ Zermelo number, and let $\univ_n$ be defined similarly, but starting from $V$ rather than $\emptyset$, i.e.:
\begin{align*}
	\zero_n &= \overbrace{\{\ldots \{}^{n\text{ times}}\emptyset\,\}\ldots\} & \univ_n= \overbrace{\{\ldots \{}^{n\text{ times}}V\,\}\ldots\}
\end{align*}
We can now consider two sequences: 
\begin{align*}
	\begin{array}{llllllll}
		\zero_0, & \zero_1, & \zero_2, & \zero_3, & \ldots, & \zero_{2n}, & \zero_{2n+1}, & \ldots\\
		\zero_0, & \univ_0, & \zero_1, & \univ_1, & \ldots, & \zero_n, & \univ_n, & \ldots
	\end{array}
\end{align*}
Inutitively, $\cointer$ treats the former sequence as the latter, and $\ocinter$ treats the latter as the former. The proof that $\cointer$ and $\ocinter$ are mutually inverse simply builds on this intuitive thought.

Here are two facts which make the intuitive thought precise:
\begin{lem}[\LTplus]\label{lem:co:znvn}
		 $\forall x\, \ x \conotin \emptyset$, and $\forall x\, \ x \coin \{\emptyset\}$, and $\forall x(x \coin \zero_{n+2} \liff x = \zero_n)$ for all $n$.
\end{lem}
\begin{lem}[\BLTplus]\label{lem:oc:znvn}
	$\helowbij(\zero_n) = \zero_{2n}$ and $\helowbij(\univ_n) = \zero_{2n+1}$, for all $n$.
\end{lem}\noindent
The proofs of both facts are trivial. Using the second fact, though, I can build up to the proof in \BLTplus that $x \in a$ iff $(x\coin a)\octrans$:
\begin{lem}[\BLTplus]\label{lem:fJ:effect}
	The function $\coinj\octrans$, i.e.\ the $\ocinter$-interpretation of \LT's definition of $\coinj$, maps $\zero_n \mapsto \univ_n \mapsto \zero_{n+1}$, and $x \mapsto x$ otherwise.
\end{lem}
\begin{proof}
	Note that $\zero_{2n} \in \zero_{2n+1} \in \zero_{2n+2}$, with these membership facts unique. So $\helowbij(\zero_{n}) \in \helowbij(\univ_{n}) \in \helowbij(\zero_{n+1})$, by Lemma \ref{lem:oc:znvn}, i.e.\ $\zero_{n} \ocin \univ_{n} \ocin \zero_{n+1}$.
\end{proof}
\begin{lem}[\BLTplus]\label{lem:hfeh}
	$\helowbij(\coinj\octrans(a)) = \coinjoc(\helowbij(a))$, for all $a$.
\end{lem}
\begin{proof}
	By Lemmas \ref{lem:oc:znvn}--\ref{lem:fJ:effect}, we have $\helowbij(\coinj\octrans(\zero_n)) = \helowbij(\univ_n) = \zero_{2n+1} = \coinjoc(\zero_{2n}) = \coinjoc(\helowbij(\zero_n))$ and 
	$\helowbij(\coinj\octrans(\univ_n)) = \helowbij(\zero_{n+1}) = \zero_{2n+2} = \coinjoc(\zero_{2n+1}) = \coinjoc(\helowbij(\univ_n))$. 
	Now suppose $a \neq \zero_n$ and $a \neq \univ_n$ for any $n$, so that $\coinj\octrans(a) = a$ and hence $\helowbij(\coinj\octrans(a)) = \helowbij(a)$; moreover, $\helowbij(a) \neq \zero_n$ for any $n$ by  Lemma \ref{lem:oc:znvn}; so $\helowbij(\coinj\octrans(a)) = \helowbij(a) = \coinjoc(\helowbij(a))$. 
\end{proof}
\begin{lem}[\BLTplus]\label{lem:coinjoctrans}
	$\helowbij(\coinj\octrans(x)) \in \helowbij(a) \liff a \notin a$ iff $x \in a$
\end{lem}
\begin{proof}
	If $a \notin a$ then $\helowbij(a) = \Setabs{\helowbij(\coinj\octrans(x))}{x \in a}$ by Lemma \ref{lem:hfeh}. If $a \in a$ then $\helowbij(a) = \Setabs{\helowbij(\coinj\octrans(x))}{x \notin a} \cup \{\emptyset\}$, and note that $\emptyset \neq \helowbij(\coinj\octrans(x)) = \coinjoc(\helowbij(x))$ for all $x$.
\end{proof}
\begin{lem}[\BLTplus]\label{lem:inv1} $x \in a$ iff $(x \coin a)\octrans$
\end{lem}
\begin{proof}
	Using Lemma \ref{lem:coinjoctrans} and the fact that $a \notin a$ iff $\helowbij(\emptyset) = \emptyset \notin \helowbij(a)$, note the following chain of equivalent formulas:
	\begin{listn-0}
		\item\label{twowandh} $x \in a$ 
		\item\label{twowandb3} $\helowbij(\coinj\octrans(x)) \in \helowbij(a) \liff \helowbij(\emptyset) \notin \helowbij(a)$
		\item\label{twowandb} $(\coinj(x) \in a \liff \emptyset \notin a )\octrans$
		\item\label{twowanda} $(x \coin a)\octrans$
		\qedhere
	\end{listn-0}
\end{proof}\noindent
It remains to show in \LTplus that $x \in a$ iff $(x \ocin a)\cotrans$. Working in \BLTplus, define $\ocinj$ as a map sending $\zero_{n+1} \mapsto \univ_n \mapsto \zero_n$ and $x \mapsto x$ otherwise; by Lemma \ref{lem:fJ:effect}, if $x \neq \emptyset$ then $\ocinj^{-1}(x) = \coinj\octrans(x)$. We then have two quick results:
\begin{lem}[\BLTplus]\label{lem:helowotherwiseput}
	$\helowbij(x) \in \helowbij(a)$ iff $(x = \emptyset \land a \in a) \lor (x \neq \emptyset \land (\ocinj(x) \in a \liff a \notin a))$
\end{lem} 
\begin{proof}
	If $x = \emptyset$, then $\helowbij(\emptyset) = \emptyset \in \helowbij(a)$ iff $a \in a$. If $x \neq \emptyset$; use Lemma \ref{lem:coinjoctrans}.	
\end{proof}
\begin{lem}[\LTplus]\label{lem:helowocinjcotrans}
	If $x \neq \emptyset$, then $\coinj(\ocinj\cotrans(x)) = x$.
\end{lem}
\begin{proof}
	By Lemma \ref{lem:co:znvn}, $\ocinj\cotrans$ maps $\zero_{n+2} \mapsto \zero_{n+1} \mapsto \zero_n$, and $x \mapsto x$ otherwise.
\end{proof}
\begin{lem}[\LTplus]\label{lem:inv2} $x \in a$ iff $(x \ocin a)\cotrans$
\end{lem}
\begin{proof}
	Using Lemmas \ref{lem:helowocinjcotrans} and \ref{lem:helowotherwiseput}$\cotrans$, note the following chain of equivalent formulas:
	\begin{listn-0}
		\item $x \in a$
		\item\label{twowand:h:2} $(\emptyset = x \land x \in a) \lor (\emptyset \neq x \land x \in a)$
		\item\label{twowand:h:3} $(\emptyset = x \land x \in a) \lor (\emptyset \neq x \land \coinj(\ocinj\cotrans(x)) \in a)$
		\item\label{twowand:h:4} $(\emptyset = x \land a \coin a) \lor (\emptyset \neq x\land (\ocinj\cotrans(x) \coin a \liff \emptyset \notin a))$
		\item\label{twowand:h:5} $(\emptyset = x \land a \coin a) \lor (\emptyset \neq x \land (\ocinj\cotrans(x) \coin a \liff a \conotin a)) $
		\item\label{twowand:h:6} $((\emptyset = x \land a \in a) \lor (\emptyset \neq x \land (\ocinj(x) \in a \liff a \notin a)))\cotrans$
		\item\label{twowand:h:7} $(\helowbij(x) \in \helowbij(a))\cotrans$	
		\item $(x \ocin a)\cotrans$ \qedhere
	\end{listn-0}
\end{proof}\noindent
Theorem \ref{thm:synonymy} now follows from Lemmas \ref{lem:lt:cotransworks}, \ref{lem:blt:octransworks}, \ref{lem:inv1}, and \ref{lem:inv2}.

\subsection{Finitary cases of definitional equivalences}
The base theories, \LT and \BLT, are not definitionally equivalent. To see this, consider:
\begin{align*}
	\clearme{lt}(1) &\coloneq  1 & \clearme{blt}(1)&\coloneq  2\\
	\clearme{lt}(n+1) &\coloneq  2^{\clearme{lt}(n)} & \clearme{blt}(n+1) &\coloneq  2^{\clearme{blt}(n)+1}
\end{align*}
Any model of \LT with $n$ levels has $\text{lt}(n)$ sets, and any model of \BLT with $n$ \boolevel{}s has $\text{blt}(n)$ sets. In particular, there is a model of \LT with four sets, but no model of \BLT has four sets. So \LT and \BLT are not definitionally equivalent.

There is, though, a nice definitional equivalence when we insist that there are infinitely many sets but that every set is finite. Concretely: let \LTfin be $\LTplus + \lnot\clearme{\ref{lt:inf}}$, and let \BLTfin be $\BLTplus + \lnot\clearme{\ref{blt:inf}}$. Our earlier results immediately entail that $\LTfin$ and $\BLTfin$ are definitionally equivalent. Moreover, as noted in Pt.\ref{pt:lt} \S\ref{s:1:ltsubzf},  \LTfin is equivalent to \ZFfin. Finally, \ZFfin and \PA are definitionally equivalent.\footnote{\textcite[Theorems 3.3, 6.5, 6.6]{KayeWong:IAST}. \ZFfin is the theory with all of \ZF's axioms except that: (i) Zermelo's axiom of infinity is replaced with its negation; and (ii) it has a new axiom, $\forall a(\exists t \supseteq a)(t\text{ is transitive})$.} So:
\begin{lem}
	\PA, \ZFfin, \LTfin, and \BLTfin are definitionally equivalent. 
\end{lem}

\section*{Acknowledgements}
Special thanks to Thomas Forster, Joel David Hamkins, Randall Holmes, and Brian King, for extensive discussion. Thanks also to Neil Barton, Sharon Berry, Luca Incurvati, Juliette Kennedy, \O{}ystein Linnebo, Michael Potter, Chris Scambler, James Studd, Rob Trueman, Sean Walsh, Will Stafford, audiences at MIT, Oxford, and Paris, and anonymous referees for \emph{Bulletin of Symbolic Logic}.
\stopappendix

\printbibliography
\end{document}